\documentstyle[11pt,leqno, amsfonts]{article}
\begin{document}
\newtheorem{tm}{Theorem}
\newtheorem{la}{Lemma}
\newtheorem{cy}{Corollary}
\newtheorem{pn}{Proposition}
\newcommand{\bx}{{\bf x}}
\newcommand{\by}{{\bf y}}
\newcommand{\bw}{{\bf w}}
\newcommand{\ekp}{\enskip}
\newcommand{\wtd}{\widetilde}
\newcommand{\noi}{\noindent}

\pagestyle{myheadings}
\markright{SPIKE TRAIN MODELS}
\begin{center}
{\large {\bf SOME THEORETICAL RESULTS ON NEURAL}} \\
{\large {\bf SPIKE TRAIN PROBABILITY MODELS
}}\footnote{
{\em AMS} 2000 {\em subject classifications}. Primary 62E20; secondary 62G20, 62M20.
}
 \\
\vspace{0.3cm}
{\sc By Hock Peng Chan and Wei-Liem Loh} \\
\vspace{0.3cm}
{\em National University of Singapore}  \\
\end{center}
\begin{quote}
This article contains two main theoretical results on neural spike train models.
The first assumes that the spike train is modeled as a counting or point process on the real line
where the conditional intensity function
is a product of a free firing rate function $s$, which depends only on the stimulus,
and a recovery function $r$, which depends only on the time since the last spike. 
If $s$ and $r$ belong to a $q$-smooth class of functions,
it is proved that sieve maximum likelihood estimators for $s$ and $r$ achieve essentially the optimal convergence rate
(except for a logarithmic factor) under $L_1$ loss. 

The second part of this article considers template matching of multiple spike trains.
$P$-values for the occurrences of a given template or
pattern in a set of spike trains are computed using a general scoring system. 
By identifying the pattern with an experimental stimulus, multiple spike
trains can be deciphered to provide useful information. 
\end{quote}

\section{Introduction}

In the field of neuroscience, it is generally  acknowledged that neurons are the basic units of information processing in the brain.
They do this by generating characteristic 
very short duration  and highly peaked electric
action potentials, or more simply spikes, in its body [see, for example, Dayan and Abbott (2001)]. 
These spikes can travel along nerve fibers that extend over relatively long
distances to other cells. The temporal pattern of these spikes depends dynamically on the stimuli of the neuron or the biochemicals
induced by the spikes of other neurons. The collection of such spikes generated by a neuron over a time period is called a spike train.
In this way, information is transmitted via spike trains.
Because the spikes are of very short duration and are highly peaked, point processes or counting processes are the most commonly used
probability models for neural spike trains, with points on the time axis representing the temporal location of the spikes
[see, for example, Brillinger (1992)].

Let $N(t)$ denote the number of spikes on the interval $[0, t)$. $N(.)$ counts the number of spikes and hence is a counting process.
Let $w_1< w_2 <\cdots < w_{N(t)}$ be all the spike times occurring in  $[0, t)$. We assume that the following limit exists:
\begin{displaymath}
\lambda (t | w_1,\cdots, w_{N(t)} ) = \lim_{\delta \downarrow 0} \frac{1}{\delta} E[ N( t+\delta) - N(t) | w_1,\cdots, w_{N(t)} ], \hspace{0.5cm} \mbox{a.s.}
\end{displaymath}
$\lambda(.| .)$ is known as the conditional intensity of $N(.)$.

In the neuroscience literature, a number of probability models for $\lambda(.| .)$ have been proposed. 
One of the simplest is when $\lambda(t| w_1,\cdots, w_{N(t)} )$ depends only on $t$.
This leads to an inhomogeneous Poisson process [see, for example, Ventura {\em et al.} (2002)]. 
It is well known that for a short period of time after a spike has been discharged, it is more difficult or even impossible for a neuron
to fire another spike [see, for example, Dayan and Abbott (2001), page 4]. 
Such a time interval is called the refractory period. The main drawback with the inhomogeneous
Poisson process model is that it does not incorporate
the refractory period of the neuron.
To account for this,  a number of researchers [for example, Johnson and Swami (1983) and Kass and Ventura (2001)] 
have proposed modeling $\lambda(.|.)$ by
\begin{displaymath}
\lambda_0 (t| w_1,\cdots, w_{N(t)} ) = f( t, t- w_{N(t)}),
\end{displaymath}
where $f$ is a nonnegative function.
This model is Markovian in that it only depends on the present time $t$ and the duration $t-w_{N(t)}$ since the last spike.

A simpler alternative model for the conditional intensity of $N(.)$
that has been proposed in the literature [see, for example, Johnson and Swami (1983), Miller (1985)
and Berry and Meister (1998)]
is
\begin{equation}
\lambda_1 (t| w_1,\cdots, w_{N(t)} ) = \left\{ \begin{array}{ll}
s(t) & \mbox{if $N(t)=0$,} \\
s(t) r( t- w_{N(t)}) & \mbox{if $N(t)\geq 1$},
\end{array}
\right.
\label{eq:1.1}
\end{equation}
where $s, r$ are nonnegative functions. 
$s$ and $r$ are known as the free firing rate function and the recovery function respectively.
This model has the added attractiveness of easy interpretability.

This article consists of two parts. The first part considers sieve maximum
likelihood estimation of $s$ and $r$
in (\ref{eq:1.1})
based on $n$ independent realizations of $N(t), t\in [0, T),$ where $0< T <\infty$
and $N(.)$ is a counting process with conditional intensity
$\lambda_1 (.|.)$.
Here we assume that the true free firing rate function $s$ and recovery function $r$ both lie in the
class of $q$-smooth functions $\Theta_{\tilde{\kappa}, q}$ where 
$\Theta_{\tilde{\kappa}, q}$ is defined as in (\ref{eq:2.62}).
Section 2 computes upper bounds on the metric entropy of $\Theta_{\tilde{\kappa}, q}$ as well as other function spaces
induced by $\Theta_{\tilde{\kappa}, q}$. 
These results are needed in Section 3.

Section 3 focuses on sieve maximum likelihood estimators $\hat{s}_n$ and $\hat{r}_n$ for $s$ and $r$
respectively. Assuming that there
exists an absolute refractory period (that is, there exists a constant $\theta>0$ such that
$r(u) = 0$, for all $u\in [0, \theta]$), it is proved in Theorem  \ref{tm:3.2} that
for $q>1/2$,
\begin{displaymath}
E_{s,r} [ \int_0^T | \hat{s}_n (t) - s(t) | dt] = O(n^{-q/(2q+1)} \log^{1/2} n), \hspace{0.5cm}\mbox{as $n\rightarrow\infty$.}
\end{displaymath}
If, in addition, $s(t)>0$ for $t\in [0, T]$, then Theorem \ref{tm:3.2} shows that
\begin{displaymath}
E_{s,r} [ \int_0^{T^*} | \hat{r}_n (u) - r(u) | du] = O(n^{-q/(2q+1)} \log^{1/2} n), \hspace{0.5cm}\mbox{as $n\rightarrow\infty$},
\end{displaymath}
where $T^*$ is an arbitrary but fixed constant satisfying $0< T^* < T$.

In Section 4, corresponding lower bounds for the convergence rate are established. 
In particular under the assumptions of Section 3, Theorems \ref{tm:4.1} and 
\ref{tm:4.2} prove that it is not possible to achieve a faster convergence rate
than $n^{-q/(2q+1)}$. Thus we conclude that sieve maximum likelihood estimators  for $s$ and $r$
achieved essentially the optimal convergence rate (except for a logarithmic factor).

The second part of the article deals with the detection of multiple spike train
patterns. Let $\bw^{(i)} = \{ w^{(i)}_1,\ldots,w^{(i)}_{N_i(T)} \}$ be the
spike times of the $i$th neuron for $1 \leq i \leq d$, and let
$\bw = (\bw^{(1)},\ldots,\bw^{(d)})$. Loosely speaking, the pattern or
template $\bw$ is said to have
occurred at time $t$ in the spike trains $\by = (\by^{(1)},\ldots,\by^{(d)})$ if
for most $y \in \by^{(i)} \cap [t,t+T)$, $1 \leq i \leq d$, there exists
$w \in \bw^{(i)}$ such that $y-w-t$ is close to 0. A more rigorous definition
of a match, via a user-chosen score function, is given in Section 5.
When the number of matches is significantly large, we can identify the
onset of the patterns $\bw$ in $\by$ with the stimulus provided to the subjects
when $\bw$ is recorded. For example, $\bw$ can be the spike times of an
assembly of neurons of a zebra finch when its own song is played while
awake and
$\by$ is the spike trains of the same assembly when it is sleeping.  
The replaying of these patterns
during sleep has been hypothesized to play an important role in bird song
learning [cf.
Dave and Margoliasch (2000) and Mooney (2000)].

In Brown, Kass and Mitra (2004), it was stated that ``research in statistics
and signal processing on multivariate point process models has not
been nearly as extensive as research on models of multivariate 
continuous-valued processes'' in the section titled ``future 
challenges for multiple spike train data analysis''. We develop in Sections
5 and 6 a theory for computing the distribution of 
scan statistics in multivariate point processes and apply it to obtain
$p$-values for template matching. The accuracy of these computations
is then verified independently via computer experiments.

\section{Metric entropy} 

In this section, suppose $N(.)$ is a counting process
with conditional intensity $\lambda_1 (.|.)$ as given by (\ref{eq:1.1}).
We assume that a realization
of $N$ is observed on the interval $[0, T)$, $0<T<\infty$, and that
the spike times are $0 < w_1 < w_2 <\cdots < w_{N(T)} < T$.
It is convenient to let $\{w_1,\cdots, w_{N(T)}\}$  denote the point process corresponding to $N(t), t\in [0, T)$, 
and ${\cal N}$ be the set of all possible realizations of $\{w_1,\cdots, w_{N(T)}\}$.
It follows from Chapter 7 of Daley and Vere-Jones (2002) that the likelihood
is the local Janossy density given by
\begin{equation}
p_{s,r} ( \{w_1,\cdots, w_{N(T)}\}) = e^{-\int_0^T s(t) r(t-w_{N(t)}) dt} \prod_{j=1}^{N(T)} s(w_j) r(w_j- w_{j-1} ),
\label{eq:4.1}
\end{equation}
and hence its log-likelihood function is
\begin{eqnarray*}
l( s, r| \{w_1,\cdots, w_{N(T)}\})
&=& - \int_0^T s(t) r( t- w_{ N(t)} ) dt 
+ \sum_{j=1}^{ N (T) } \log[ s( w_j ) r (w_j - w_{ j-1} ) ]
\nonumber \\
&=& 
- \int_0^{w_1 \wedge T}  s(t) dt 
- \int_{w_1 \wedge T}^T  s(t) r( t- w_{ N(t)} ) dt 
\nonumber \\
&& + \sum_{j=1}^{ N (T) } \log[ s( w_j ) r (w_j - w_{ j-1} ) ],
\end{eqnarray*}
where $r( t - w_0)= 1$ for all $t\in [0, T)$.

Next let $q, q_0, q_1$ be constants satisfying $q>0$, $q = q_0+q_1$, $q_0$ a nonnegative integer and $0< q_1\leq 1$.
Furthermore we write $\tilde{\kappa} = (\kappa_0,\cdots, \kappa_{q_0+1} )$ to be a vector of strictly positive constants.
In this section, we assume that the true free firing rate function $s$ and the recovery function $r$ lie in the
$q$-smooth function class $\Theta_{\tilde{\kappa}, q}$ where
\begin{eqnarray}
\Theta_{\tilde{\kappa}, q} &=& \Big\{ f=g^2: g\in {\cal C}^{q_0} [0, T), \min_{t\in [0, T)} g(t) \geq 0, \max_{t\in [0, T)} |\frac{ d^j}{dt^j} g(t) | < \kappa_j,
j=0,\cdots, q_0,
\nonumber \\
&&\hspace{0.5cm}
|\frac{d^{q_0}}{ dt^{q_0}} g (t_1) - \frac{d^{q_0} }{dt^{q_0}} g (t_2)| < \kappa_{q_0+1} |t_1-t_2|^{q_1}, \forall t_1, t_2 \in [0, T) \Big\}.
\label{eq:2.62}
\end{eqnarray}
Let $\{0< \delta_n\leq 1: n= 1, 2, \cdots\}$ be a sequence of constants (to be suitably chosen later and $\delta_n$
depends only on $n$) such that
$\delta_n \rightarrow 0$ as $n\rightarrow \infty$. We define a sieve for the parameter space of $\Theta_{\tilde{\kappa}, q}$ by
\begin{eqnarray*}
\Theta_{\tilde{\kappa}, q, n} &=& \Big\{ f=g^2: g\in {\cal C}^{q_0} [0, T), \min_{t\in [0, T)} g(t) \geq \delta_n, 
\max_{t\in [0, T)} |\frac{d^j }{ dt^j} g (t)| < \kappa_j,
j=0,\cdots, q_0,
\nonumber \\
&&\hspace{0.5cm}
|\frac{ d^{q_0} }{ dt^{q_0}} g (t_1) - \frac{ d^{q_0} }{ dt^{q_0}} g (t_2)| < \kappa_{q_0+1} |t_1-t_2|^{q_1}, \forall t_1, t_2 \in [0, T) \Big\}.
\end{eqnarray*}
Let $\Theta_{\tilde{\kappa}, q}$ and $\Theta_{\tilde{\kappa}, q, n}$ be endowed with the metrics 
$\rho_{\Theta_{\tilde{\kappa}, q}}$ and
$\rho_{\Theta_{\tilde{\kappa}, q, n}}$ respectively where
\begin{eqnarray*}
\rho_{\Theta_{\tilde{\kappa}, q}} (f_1, f_2) &=& \sup_{t\in [0, T)} | f_1^{1/2}(t) - f_2^{1/2}(t) |, 
\hspace{0.5cm}\forall f_1, f_2 \in \Theta_{\tilde{\kappa}, q},
\nonumber \\
\rho_{\Theta_{\tilde{\kappa}, q, n}} (f_1, f_2) &=& \sup_{t\in [0, T)} | f_1^{1/2}(t) - f_2^{1/2}(t) |, 
\hspace{0.5cm}\forall f_1, f_2 \in \Theta_{\tilde{\kappa}, q, n}.
\end{eqnarray*}
We observe that any $f\in \Theta_{\tilde{\kappa}, q}$ can be approximated arbitrarily closely by $ (f^{1/2}+\delta_n)^2 \in \Theta_{\tilde{\kappa}, q, n}$ by choosing 
$n$ sufficiently large.
Consequently a sieve for the parameter space of $(s, r)$  can now be expressed as
$\Theta_{\tilde{\kappa}, q, n}^2 = \Theta_{\tilde{\kappa}, q, n} \times \Theta_{\tilde{\kappa}, q, n}$ 
with metric $\rho_{\Theta_{\tilde{\kappa}, q, n}^2}$ where
\begin{displaymath}
\rho_{\Theta_{\tilde{\kappa}, q, n}^2} ((f_1, g_1), (f_2, g_2)) = \rho_{\Theta_{\tilde{\kappa}, q, n}}( f_1, f_2) 
+ \rho_{\Theta_{\tilde{\kappa}, q, n}} (g_1, g_2), \hspace{0.5cm}
\forall (f_1,g_1 ), (f_2, g_2) \in \Theta_{\tilde{\kappa}, q, n}^2.
\end{displaymath}
Next let
\begin{displaymath}
{\cal F}_{\tilde{\kappa}, q, n} = \mbox{ \{$p_{s_1,r_1}$ is as in (\ref{eq:4.1}): $(s_1,r_1) \in \Theta_{\tilde{\kappa}, q, n}^2$\}},
\end{displaymath}
be endowed with the Hellinger metric $\rho_{{\cal F}_{\tilde{\kappa}, q, n}}$ where
\begin{eqnarray*}
&& \rho_{{\cal F}_{\tilde{\kappa}, q, n}} (p_{s_1, r_1}, p_{s_2, r_2}) 
\nonumber \\
&=& \| p_{s_1, r_1}^{1/2} - p_{s_2,r_2}^{1/2} \|_2
\nonumber \\
&=& \Big\{ \sum_{j=0}^\infty \int_{0<w_1<\cdots < w_j<T} [ p_{s_1, r_1}^{1/2} ( \{w_1,\cdots, w_j\}) - p_{s_2,r_2}^{1/2} ( \{w_1,\cdots, w_j\}) 
]^2 dw_1\cdots dw_j \Big\}^{1/2}.
\end{eqnarray*}
For $\varepsilon>0$, let $\Theta_{\tilde{\kappa}, q, n}^2 (\varepsilon ) \subseteq \Theta_{\tilde{\kappa}, q, n}^2$ denote 
a finite $\varepsilon$-net for $\Theta_{\tilde{\kappa}, q, n}^2$ with
respect to the metric $\rho_{\Theta_{\tilde{\kappa}, q, n}^2}$. 
This implies that for each $(s_1, r_1) \in \Theta_{\tilde{\kappa}, q, n}^2$, 
there exists a $(s_2, r_2) \in \Theta_{\tilde{\kappa}, q, n}^2 (\varepsilon)$ such that
$\rho_{\Theta_{\tilde{\kappa}, q, n}^2} ((s_1,r_1), (s_2, r_2) ) \leq \varepsilon$.

Now suppose that for each $\varepsilon >0$, there exist measurable nonnegative functions $f_{l,\varepsilon}$ and $f_{u,\varepsilon}$ 
on $\Theta_{\tilde{\kappa}, q, n}^2 (\varepsilon) \times {\cal N}$ such that for each $(s_1,r_1) \in \Theta_{\tilde{\kappa}, q, n}^2$, there
is some $(s_2, r_2) \in \Theta_{\tilde{\kappa}, q, n}^2 (\varepsilon)$ satisfying
\begin{equation}
\rho_{\Theta_{\tilde{\kappa}, q, n}^2 } ((s_1,r_1), (s_2, r_2) ) \leq \varepsilon,
\label{eq:4.5}
\end{equation}
with
\begin{eqnarray}
f_{l,\varepsilon} ((s_2, r_2), \{w_1,\cdots, w_{N(T)}\})  
&\leq & p_{s_1, r_1}( \{w_1,\cdots, w_{N(T)}\}) 
\nonumber \\
&\leq &
f_{u,\varepsilon} ((s_2, r_2), \{w_1,\cdots, w_{N(T)}\}), \hspace{0.5cm}\mbox{a.s.,}
\label{eq:4.6}
\end{eqnarray}
and
\begin{eqnarray}
&& \Big\{ \sum_{j=0}^\infty \int_{0<w_1<\cdots < w_j<T} [ f_{u,\varepsilon}^{1/2} ((s_2, r_2), \{w_1,\cdots, w_j\}) 
\nonumber \\
&&\hspace{0.5cm} - f_{l,\varepsilon}^{1/2} ((s_2,r_2), \{w_1,\cdots, w_j\}) ]^2 dw_1 \cdots dw_j \Big\}^{1/2}
\leq \varepsilon.
\label{eq:4.7}
\end{eqnarray}

{\sc Definition.} For $\varepsilon>0$, the $\varepsilon$-entropy of $\Theta_{\tilde{\kappa}, q, n}^2$ with respect to 
the metric $\rho_{\Theta_{\tilde{\kappa}, q, n}^2}$ is defined to be
\begin{eqnarray*}
H (\varepsilon, \Theta_{\tilde{\kappa}, q, n}^2, \rho_{\Theta_{\tilde{\kappa}, q, n}^2} ) 
&=& \log [ \min\{ \mbox{card $\Theta_{\tilde{\kappa}, q, n}^2 (\varepsilon)$: 
$\Theta_{\tilde{\kappa}, q, n}^2 (\varepsilon)$
is a $\varepsilon$-net for $\Theta_{\tilde{\kappa}, q, n}^2$}
\nonumber \\
&& \hspace{0.5cm}\mbox{with respect to the metric $\rho_{\Theta_{\tilde{\kappa}, q, n}^2}$} \} ].
\end{eqnarray*}
The $\varepsilon$-entropies of 
$\Theta_{\tilde{\kappa}, q}$ with respect to $\rho_{\Theta_{\tilde{\kappa}, q}}$ and
$\Theta_{\tilde{\kappa}, q, n}$ with respect to $\rho_{\Theta_{\tilde{\kappa}, q, n}}$ are defined
in a similar manner.

{\sc Definition.}
The $\varepsilon$-entropy of ${\cal F}_{\tilde{\kappa}, q, n}$ with bracketing
with respect to the metric $\rho_{{\cal F}_{\tilde{\kappa}, q, n}}$ is defined to be
\begin{eqnarray*}
H^B (\varepsilon, {\cal F}_{\tilde{\kappa}, q, n}, \rho_{{\cal F}_{\tilde{\kappa}, q, n}}) 
&=& \log [ \min\{ \mbox{card $\Theta_{\tilde{\kappa}, q, n}^2 (\varepsilon)$:
(\ref{eq:4.5}), (\ref{eq:4.6}) and (\ref{eq:4.7}) are satisfied} \} ].
\end{eqnarray*}

We observe from Kolmogorov and Tihomirov (1961), page 308, and Dudley (1999), page 11, that
the $\varepsilon$-entropy of $\Theta_{\tilde{\kappa}, q, n}$ satisfies
\begin{displaymath}
H( \varepsilon, \Theta_{\tilde{\kappa}, q, n}, \rho_{\Theta_{\tilde{\kappa}, q, n}} ) \leq 
H( \varepsilon, \Theta_{\tilde{\kappa}, q}, \rho_{\Theta_{\tilde{\kappa}, q}} ) \leq
\frac{C_{\tilde{\kappa}, q} }{\varepsilon^{1/q} },
\end{displaymath}
and hence
\begin{equation}
H (\varepsilon, \Theta_{\tilde{\kappa}, q, n}^2, \rho_{\Theta_{\tilde{\kappa}, q, n}^2} ) 
\leq 2 H(\frac{\varepsilon}{2}, \Theta_{\tilde{\kappa}, q, n}, \rho_{\Theta_{\tilde{\kappa}, q, n}})
\leq \frac{ 2^{(q+1)/q} C_{\tilde{\kappa}, q} }{\varepsilon^{1/q} },
\label{eq:4.20}
\end{equation}
where $C_{\tilde{\kappa}, q}$ is a constant depending only on $\tilde{\kappa}$ and $q$.
Thus we conclude from Lemma \ref{la:a.3} in Appendix A that
\begin{displaymath}
H^B (\varepsilon, {\cal F}_{\tilde{\kappa}, q, n}, \rho_{{\cal F}_{\tilde{\kappa}, q, n}} )
\leq \frac{ 2^{(q+ 2)/q} C_{\tilde{\kappa}}^{1/q} C_{\tilde{\kappa}, q} }{ \varepsilon^{1/q}},
\end{displaymath}
where $C_{\tilde{\kappa}}$ is a constant that depends only on
$\tilde{\kappa}$.

Next let $f: {\cal N} \rightarrow R$ be a nonnegative function such that
\begin{displaymath}
\sum_{j=0}^\infty \int_{0\leq w_1<\cdots < w_j <T} f(\{w_1,\cdots, w_j\}) dw_1 \cdots dw_j < \infty.
\end{displaymath} 
We follow Wong and Shen (1995) by defining
\begin{displaymath}
Z_f (\{w_1,\cdots, w_{N(T)} \}) = \log [ \frac{ f (\{w_1,\cdots, w_{N(T)} \}) }{
p_{s,r} (\{w_1,\cdots, w_{N(T)} \}) } ],
\end{displaymath}
where $s$ is the true free firing rate function and $r$ the true recovery function.
For $\tau>0$, we write
\begin{eqnarray*}
&& \tilde{f} (\{w_1,\cdots, w_{N(T)} \}) 
\nonumber \\
&=&
\left\{ \begin{array}{ll}
f (\{ w_1,\cdots, w_{N(T)} \}), & \mbox{if $f (\{ w_1,\cdots, w_{N(T)} \}) \geq e^{-\tau} p_{s,r}(\{w_1,\cdots, w_{N(T)}\} )$}, \\
e^{-\tau} p_{s,r}(\{ w_1,\cdots, w_{N(T)} \}), & \mbox{if $f (\{ w_1,\cdots, w_{N(T)} \}) < e^{-\tau} p_{s,r}(\{w_1,\cdots, w_{N(T)}\} )$}, \\
\end{array}
\right.
\nonumber
\end{eqnarray*} 
and
\begin{equation}
\tilde{Z}_f 
= Z_{\tilde{f}}
= \left\{ \begin{array}{ll}
Z_f, & \mbox{if $Z_f \geq -\tau$,}
\\
-\tau, &
\mbox{if $Z_f < -\tau$}.
\end{array}
\right.
\label{eq:4.13}
\end{equation}
Let $\tilde{\cal Z}_{\tilde{\kappa}, q, n} = \{ \tilde{Z}_{p_{s_1, r_1}}: p_{s_1, r_1} \in {\cal F}_{\tilde{\kappa}, q, n} \}$ be the space of
truncated log-likelihood ratios (based on one observation).
Define $H^B (\varepsilon, \tilde{\cal Z}_{\tilde{\kappa}, q, n}, \rho_{\tilde{\cal Z}_{\tilde{\kappa}, q, n}})$ 
to be the $\varepsilon$-entropy of $\tilde{\cal Z}_{\tilde{\kappa}, q, n}$ with
bracketing with respect to the metric
\begin{eqnarray*}
&& \rho_{\tilde{\cal Z}_{\tilde{\kappa}, q, n} } (\tilde{Z}_{p_{s_1, r_1}}, \tilde{Z}_{p_{s_2, r_2}} )
\nonumber \\
&=&
\Big\{ E_{s,r} [ \tilde{Z}_{p_{s_1, r_1}} (\{ w_1,\cdots, w_{N(T)} \}) - \tilde{Z}_{p_{s_2, r_2}} (\{w_1,\cdots, w_{N(T)} \}) ]^2 \Big\}^{1/2}
\nonumber \\
&=& \Big\{ \sum_{j=0}^\infty \int_{0<w_1<\cdots < w_j<T} [ \tilde{Z}_{p_{s_1, r_1}} (\{w_1,\cdots, w_j\}) - 
\tilde{Z}_{p_{s_2, r_2}} (\{ w_1,\cdots, w_j\}) ]^2
\nonumber \\
&&\hspace{0.5cm} \times
e^{-\int_0^T s(t) r(t- w_{\zeta(t)}) dt} [ \prod_{i=1}^j s(w_i) r(w_i - w_{i-1}) ] dw_1\cdots dw_j \Big\}^{1/2},
\end{eqnarray*}
where $\zeta(t) = \max\{ k\geq 0: w_k <t\}$ and
$E_{s,r}$ denotes expectation when the true free firing rate function is $s$ and the recovery function is $r$.
We observe from Lemma \ref{la:a.6} in Appendix A that
\begin{equation}
H^B (\varepsilon, \tilde{\cal Z}_{\tilde{\kappa}, q, n}, \rho_{\tilde{Z}_{\tilde{\kappa}, q, n}} )
\leq 2^{(q+ 2)/q} C_{\tilde{\kappa}}^{1/q} C_{\tilde{\kappa}, q} (\frac{ 2 e^{\tau/2} }{ \varepsilon})^{1/q}.
\label{eq:2.10}
\end{equation}

\section{Sieve maximum likelihood estimation}

In this section,  we assume that we have $n$ independent identically distributed copies of $N(t), t\in [0, T)$,
with conditional intensity as given by (\ref{eq:1.1}).
Let these $n$ copies be denoted by $N_i(t), t\in [0, T)$, and 
the spike times be written as $0< w_{i,1} <\cdots < w_{i, N_i (T)}< T$, $i=1,\cdots, n$.
Inspired by Wong and Shen (1995), we shall first establish a number of 
likelihood ratio probability inequalities.

\begin{pn} \label{pn:4.1}
Let $0< \varepsilon <1$ and $C_{\tilde{\kappa}}, C_{\tilde{\kappa}, q}$ be as in (\ref{eq:2.10}).
Suppose that
\begin{equation}
\int_{\varepsilon^2/2^8}^{\sqrt{2} \varepsilon} 
[ 2^{(q+2)/q} C_{\tilde{\kappa}}^{1/q} C_{\tilde{\kappa}, q} (\frac{ 10}{ x})^{1/q} ]^{1/2} dx
\leq \frac{ n^{1/2} \varepsilon^2 }{ 2^{13} \sqrt{2} }.
\label{eq:4.90}
\end{equation}
Then
\begin{displaymath}
P_{s,r}^* \{ \sup_{\| p_{s_1,r_1}^{1/2} - p_{s,r}^{1/2} \|_2 \geq \varepsilon, p_{s_1, r_1}\in {\cal F}_{\tilde{\kappa}, q, n} }
\prod_{i=1}^n \frac{p_{s_1, r_1} (\{ w_{i,1},\cdots, w_{i, N_i (T)} \}) }{ p_{s, r} (\{ w_{i,1},\cdots, w_{i, N_i(T)} \} )}
\geq e^{- n \varepsilon^2/8} \}
\leq 4 \exp[ -\frac{ n\varepsilon^2}{ 2^7 (250)}].
\end{displaymath}
$P^*_{s,r}$ is the outer measure corresponding to the density $p_{s,r}$.
\end{pn}
We refer the reader to Appendix A for a proof of Proposition \ref{pn:4.1}.

Next we define nonnegative functions $s^\dag_n$ and $r^\dag_n$ on $t\in [0, T)$ by
\begin{equation}
\sqrt{ s^\dag_n (t)} = \sqrt{ s (t)} + \delta_n, \hspace{0.5cm}
\sqrt{ r^\dag_n (t)} = \sqrt{ r (t)} + \delta_n,
\label{eq:3.23}
\end{equation}
where $\{0< \delta_n \leq 1: n=1, 2, \cdots\}$ is as in Section 2.
Since $s, r\in \Theta_{\tilde{\kappa}, q}$,  we have $s^\dag_n, r^\dag_n \in \Theta_{\tilde{\kappa}, q, n}$ 
and $p_{s^\dag_n, r^\dag_n} \in {\cal F}_{\tilde{\kappa}, q, n}$ for sufficiently large $n$.
We further observe from Lemma 8 of Wong and Shen (1995) and Lemma \ref{la:a.43} in Appendix A that
\begin{displaymath}
0 \leq \delta^\dag_n := E_{s,r} (\frac{ p_{s,r} }{p_{s^\dag_n, r^\dag_n } } -1) \leq C_{\tilde{\kappa}, 1} \delta_n,
\end{displaymath}
where 
$C_{\tilde{\kappa}, 1}$ is a constant depending only on $\tilde{\kappa}$.

\begin{pn} \label{pn:4.2}
Let $0< \varepsilon <1, \delta^\dag_n \leq 1$ and that (\ref{eq:4.90}) holds.
Then
\begin{eqnarray*}
&& P_{s,r}^* \{ \sup_{ \|p_{s_1, r_1}^{1/2} - p_{s,r}^{1/2} \|_2 \geq \varepsilon, p_{s_1, r_1} \in {\cal F}_{\tilde{\kappa}, q, n} }
\prod_{i=1}^n \frac{ p_{s_1, r_1} (\{ w_{i,1},\cdots, w_{i, N_i(T)} \} ) }{ p_{s^\dag_n, r^\dag_n } (\{ w_{i,1},\cdots, w_{i, N_i(T)} \}) }
\geq e^{ -n \varepsilon^2/16} \}
\nonumber \\
&\leq & 4 \exp[ - \frac{ n \varepsilon^2 }{ 2^7 (250) } ]
+ \exp[ - n( \frac{ \varepsilon^2 }{16} - \delta^\dag_n )].
\end{eqnarray*}
\end{pn}
{\sc Proof.} First we observe from Proposition \ref{pn:4.1} that
\begin{eqnarray*}
&& P_{s,r}^* \{ \sup_{ \|p_{s_1, r_1}^{1/2} - p_{s,r}^{1/2} \|_2 \geq \varepsilon, p_{s_1, r_1} \in {\cal F}_{\tilde{\kappa}, q, n} }
\prod_{i=1}^n \frac{ p_{s_1, r_1} (\{ w_{i,1},\cdots, w_{i, N_i(T)} \} ) }{ p_{s^\dag_n, r^\dag_n } (\{ w_{i,1},\cdots, w_{i, N_i(T)} \}) }
\geq e^{ - n \varepsilon^2/16} \}
\nonumber \\
&\leq &
P_{s,r}^* \{ \sup_{ \|p_{s_1, r_1}^{1/2} - p_{s,r}^{1/2} \|_2 \geq \varepsilon, p_{s_1, r_1} \in {\cal F}_{\tilde{\kappa}, q, n} }
\prod_{i=1}^n \frac{ p_{s_1, r_1} (\{ w_{i,1},\cdots, w_{i, N_i(T)} \} ) }{ p_{s, r } (\{ w_{i,1},\cdots, w_{i, N_i(T)} \}) }
\geq e^{ - n \varepsilon^2/8} \}
\nonumber \\
&& + P_{s,r} \{ 
\prod_{i=1}^n \frac{ p_{s, r} (\{ w_{i,1},\cdots, w_{i, N_i(T)} \} ) }{ p_{s^\dag_n, r^\dag_n } (\{ w_{i,1},\cdots, w_{i, N_i(T)} \}) }
\geq e^{ n \varepsilon^2/16} \}
\nonumber \\
&\leq & 4 \exp[ - \frac{ n \varepsilon^2 }{ 2^7 (250) } ]
+ P_{s,r} \{ 
\prod_{i=1}^n \frac{ p_{s, r} (\{ w_{i,1},\cdots, w_{i, N_i(T)} \} ) }{ p_{s^\dag_n, r^\dag_n } (\{ w_{i,1},\cdots, w_{i, N_i(T)} \}) }
\geq e^{ n \varepsilon^2/16} \}.
\end{eqnarray*}
Now using Markov's inequality,
\begin{eqnarray*}
&& P_{s,r} \{ 
\prod_{i=1}^n \frac{ p_{s, r} (\{ w_{i,1},\cdots, w_{i, N_i(T)} \} ) }{ p_{s^\dag_n, r^\dag_n } (\{ w_{i,1},\cdots, w_{i, N_i(T)} \}) }
\geq e^{ n \varepsilon^2/16} \}
\nonumber \\
&\leq & 
e^{ - n \varepsilon^2/16}
\prod_{i=1}^n E_{s,r} [ \frac{ p_{s, r} (\{ w_{i,1},\cdots, w_{i, N_i(T)} \} ) }{ p_{s^\dag_n, r^\dag_n } (\{ w_{i,1},\cdots, w_{i, N_i(T)} \}) }]
\nonumber \\
&=& (1 + \delta^\dag_n )^n
e^{ - n \varepsilon^2/16}
\nonumber \\
&\leq &
\exp( -\frac{ n \varepsilon^2 }{16} + n \delta^\dag_n ).
\end{eqnarray*}
This proves Proposition \ref{pn:4.2}. \hfill $\Box$

{\sc Definition.}
Let $\eta_n$ be a sequence of positive numbers converging to 0. We call an estimator $p_{\hat{s}_n, \hat{r}_n}:
{\cal N}^n \rightarrow {\mathbb{R}}^+$ a $\eta_n$-sieve MLE of $p_{s,r}$ if
$(\hat{s}_n, \hat{r}_n) \in \Theta_{\tilde{\kappa}, q, n}^2$ and
\begin{displaymath}
\frac{1}{n} \sum_{i=1}^n \log [ p_{\hat{s}_n, \hat{r}_n} (\{ w_{i,1},\cdots, w_{i, N_i(T)} \}) ]
\geq \sup_{p_{s_1, r_1}\in {\cal F}_{\tilde{\kappa}, q, n} } \frac{1}{n} \sum_{i=1}^n \log[ p_{s_1, r_1} 
(\{ w_{i,1},\cdots, w_{i, N_i(T)} \}) ] - \eta_n.
\end{displaymath}
The corresponding $\hat{s}_n$ and $\hat{r}_n$ are called $\eta_n$-sieve MLEs of $s$ and $r$ respectively.

{\sc Definition.} If $\delta_n$ in Section 2 satisfies $\delta_n =0$ for $n=1,2, \cdots$, then
a $\eta_n$-sieve MLE $p_{\hat{s}_n, \hat{r}_n}$ is more simply called a $\eta_n$-MLE of $p_{s,r}$.

\begin{pn} \label{pn:4.50}
Let $\varepsilon_n>0$ be the smallest value of $\varepsilon$
satisfying (\ref{eq:4.90}) and $0< \eta_n < \varepsilon_n^2/16$. 
If $p_{\hat{s}_n,\hat{r}_n}$ is a $\eta_n$-sieve MLE of $p_{s,r}$, then
\begin{displaymath}
P_{s, r} ( \| p_{\hat{s}_n, \hat{r}_n}^{1/2} - p_{s,r}^{1/2} \|_2 \geq \varepsilon_n) 
\leq 4 \exp[ - \frac{ n \varepsilon_n^2 }{ 2^7 (250) } ]
+ \exp[ - n( \frac{ \varepsilon_n^2 }{16} - \delta^\dag_n )],
\end{displaymath}
where $P_{s,r}$ denotes probability when the true free firing rate function is $s$ and the recovery function is $r$.
\end{pn}
{\sc Proof.} We observe from Proposition \ref{pn:4.2} that
\begin{eqnarray*}
&& P_{s,r} ( \| p_{\hat{s}_n, \hat{r}_n}^{1/2} - p_{s,r}^{1/2} \|_2 \geq \varepsilon_n)
\nonumber \\
&\leq & 
P_{s,r}^* \{ \sup_{ \| p_{s_1, r_1}^{1/2} - p_{s,r}^{1/2} \|_2 \geq \varepsilon_n, p_{s_1, r_1} \in {\cal F}_{\tilde{\kappa}, q, n} }
\prod_{i=1}^n \frac{ p_{s_1, r_1} (\{w_{i,1},\cdots, w_{i, N_i(T)} \}) }{
p_{s^\dag_n, r^\dag_n } (\{w_{i,1},\cdots, w_{i, N_i(T)} \}) }
\geq e^{-n \eta_n} \}
\nonumber \\
&\leq & 
P^*_{s,r} \{ \sup_{ \| p_{s_1, r_1}^{1/2} - p_{s,r}^{1/2} \|_2 \geq \varepsilon_n, p_{s_1, r_1} \in {\cal F}_{\tilde{\kappa}, q, n} }
\prod_{i=1}^n \frac{ p_{s_1, r_1} (\{w_{i,1},\cdots, w_{i, N_i(T)} \}) }{
p_{s^\dag_n, r^\dag_n } (\{w_{i,1},\cdots, w_{i, N_i(T)} \}) }
\geq e^{-n \varepsilon_n^2/16 } \}
\nonumber \\
&\leq & 4 \exp[ - \frac{ n \varepsilon_n^2 }{ 2^7 (250) } ]
+ \exp[ - n( \frac{ \varepsilon_n^2 }{16} - \delta^\dag_n )].
\end{eqnarray*}
This proves Proposition \ref{pn:4.50}. \hfill $\Box$

\begin{tm} \label{tm:3.1}
Let $\varepsilon_n>0$ be the smallest value of $\varepsilon$
satisfying (\ref{eq:4.90}), $q>1/2$ and $0< \eta_n < \varepsilon_n^2/16$. 
If $p_{\hat{s}_n,\hat{r}_n}$ is a $\eta_n$-MLE of $p_{s,r}$, then
\begin{displaymath}
E_{s,r}  \| p_{\hat{s}_n, \hat{r}_n}^{1/2} - p_{s,r}^{1/2} \|_2 = O(n^{-q/(2q+1)}), 
\hspace{0.5cm} \mbox{as $n\rightarrow\infty$.}
\end{displaymath}
\end{tm}
{\sc Proof.}
We observe that $\delta^\dag_n =0$ for $n=1,2,\cdots$ (from the definition of $\eta_n$-MLE) 
and $\varepsilon_n$ is exactly of order $n^{-q/(2q+1)}$
as $n\rightarrow \infty$.
Hence we observe from Proposition \ref{pn:4.50} that
\begin{eqnarray*}
E_{s,r}  \| p_{\hat{s}_n, \hat{r}_n}^{1/2} - p_{s,r}^{1/2} \|_2  
&\leq & \varepsilon_n
+ 8 \exp[ - \frac{ n \varepsilon_n^2 }{ 2^7 (250) } ]
+2  \exp[ - n( \frac{ \varepsilon_n^2 }{16} - \delta^\dag_n )]
\nonumber \\
&=& O(n^{-q/(2q+1)}),
\end{eqnarray*}
as $n\rightarrow\infty$.
\hfill $\Box$

Now we assume that there exists a refractory period in which the neuron cannot discharge another spike after 
a spike had been fired [see, for example, Brillinger (1992) and Johnson and Swami (1983)]. More precisely, we suppose that there exists a constant $\theta>0$ such that
\begin{equation}
r( u) = 0, \hspace{0.5cm}\forall u\in [0, \theta].
\label{eq:4.86}
\end{equation}
Then the number of spikes on the interval $[0, T)$ can be at most $n_\theta = \lceil T/\theta \rceil$.

\begin{pn} \label{pn:4.6}
Let $\varepsilon_n>0$ be the smallest value of $\varepsilon$
satisfying (\ref{eq:4.90}) and $0< \eta_n < \varepsilon_n^2/16 \leq (1 - e^{-1})^2/32$. 
If (\ref{eq:4.86}) holds and $p_{\hat{s}_n,\hat{r}_n}$ is a $\eta_n$-sieve MLE of $p_{s,r}$, then
\begin{eqnarray*}
&& P_{s,r} \Big\{ \sum_{j=0}^{n_\theta } \int_{0<w_1<\cdots < w_j <T} p_{s,r} (\{w_1,\cdots, w_j\}) \log [\frac{ p_{s,r} (\{w_1,\cdots, w_j\}) }{
p_{\hat{s}_n, \hat{r}_n} (\{w_1,\cdots, w_j\}) } ]
dw_1\cdots dw_j
\nonumber \\
&& \hspace{0.5cm} > [ 6 + \frac{ 2 \log (2) }{(1- e^{-1})^2} +
8 \max\{1, \log( \frac{ e^{ \bar{\kappa}^4 (\bar{\kappa}^4 +1) T/2}  }{ \varepsilon_n \delta_n^{ 2 n_\theta} } )   \} ] \varepsilon_n^2 \Big\}
\nonumber \\
&\leq & 4 \exp[ - \frac{ n \varepsilon_n^2 }{ 2^7 (250) } ]
+ \exp[ - n( \frac{ \varepsilon_n^2 }{16} - \delta^\dag_n )],
\end{eqnarray*}
where $\bar{\kappa} = \kappa_0 \vee 1$.
\end{pn}
{\sc Proof.} First we observe that
\begin{eqnarray*}
&& \sum_{j=0}^{n_\theta } \int_{0<w_1<\cdots < w_j<T}  \frac{ p_{s,r}^2 (\{w_1,\cdots, w_j\}) }{ p_{\hat{s}_n, \hat{r}_n}(\{w_1,\cdots, w_j\}) } 
dw_1 \cdots dw_j
\nonumber \\
&= & 
\sum_{j=0}^{n_\theta } \int_{0<w_1<\cdots < w_j<T}  \frac{ e^{-2 \int_0^T s(t) r(t- w_{\zeta(t)}) dt} \prod_{i=1}^j s^2 (w_i)
r^2 (w_i-w_{i-1}) }{ 
e^{- \int_0^T \hat{s}_n (t) \hat{r}_n (t- w_{\zeta(t)}) dt} \prod_{i=1}^j \hat{s}_n (w_i)
\hat{r}_n (w_i-w_{i-1}) } 
dw_1 \cdots dw_j
\nonumber \\
&\leq & \frac{ e^{\bar{\kappa}^4 (\bar{\kappa}^4 +1) T} }{\delta_n^{4 n_\theta }}.
\end{eqnarray*}
Now we observe from Theorem 5 of Wong and Shen (1995) that
\begin{displaymath}
\| p_{\hat{s}_n, \hat{r}_n}^{1/2} - p_{s,r}^{1/2} \|_2^2 \leq  \varepsilon_n^2 \leq (1-e^{-1})^2/2
\end{displaymath}
implies that
\begin{eqnarray*}
&& \sum_{j=0}^{n_\theta } \int_{0<w_1<\cdots < w_j<T}  p_{s,r} (\{w_1,\cdots, w_j\}) 
\log[ \frac{ p_{s,r} (\{w_1,\cdots, w_j\}) }{ p_{\hat{s}_n, \hat{r}_n}(\{w_1,\cdots, w_j\}) }] 
dw_1 \cdots dw_j
\nonumber \\
&\leq & [ 6 + \frac{ 2 \log (2) }{(1- e^{-1})^2} +
8 \max\{1, \log( \frac{ e^{\bar{\kappa}^4 (\bar{\kappa}^4 +1) T/2} }{ \varepsilon_n \delta_n^{ 2 n_\theta} } )   \} ] \varepsilon_n^2.
\end{eqnarray*}
Hence it follows from Proposition \ref{pn:4.50} that
\begin{eqnarray*}
&& P_{s,r} \Big\{ \sum_{j=0}^{n_\theta } \int_{0<w_1<\cdots < w_j<T}  p_{s,r} (\{w_1,\cdots, w_j\}) 
\log[ \frac{ p_{s,r} (\{w_1,\cdots, w_j\}) }{ p_{\hat{s}_n, \hat{r}_n}(\{w_1,\cdots, w_j\}) }] 
dw_1 \cdots dw_j
\nonumber \\
&& \hspace{0.5cm} > [ 6 + \frac{ 2 \log (2) }{(1- e^{-1})^2} +
8 \max\{1, \log( \frac{ e^{ \bar{\kappa}^4 (\bar{\kappa}^4 +1) T/2}  }{ \varepsilon_n \delta_n^{2 n_\theta } } )  \} ] \varepsilon_n^2 \Big\}
\nonumber \\
&\leq &
4 \exp[ - \frac{ n \varepsilon_n^2 }{ 2^7 (250) } ]
+ \exp[ - n( \frac{ \varepsilon_n^2 }{16} - \delta^\dag_n )].
\end{eqnarray*}
This proves Proposition \ref{pn:4.6}. \hfill $\Box$

The following theorem is the main result of this section.

\begin{tm} \label{tm:3.2}
Let $\varepsilon_n>0$ be the smallest value of $\varepsilon$
satisfying (\ref{eq:4.90}), $q>1/2$, $0< \eta_n < \varepsilon_n^2/16$
and $\delta_n = n^{-\alpha}$ for some constant $\alpha\in (2q/(2q+1), 1)$.
Suppose that (\ref{eq:4.86}) holds and $\hat{s}_n,\hat{r}_n$ are $\eta_n$-sieve MLEs of $s,r$ respectively. 
Then
\begin{displaymath}
E_{s,r} [ \int_0^T |\hat{s}_n (t) - s(t)| dt ] = O( n^{-q/(2q +1)} \log^{1/2} n ),
\hspace{0.5cm}\mbox{as $n\rightarrow\infty$.}
\end{displaymath}
If, in addition, $s(t) >0$ for all $t\in [0, T]$, then
\begin{displaymath}
E_{s,r} [ \int_0^{T^*} |\hat{r}_n (u) - r(u)| du ] = O( n^{-q/(2q+1)} \log^{1/2} n),
\hspace{0.5cm}\mbox{as $n\rightarrow\infty,$}
\end{displaymath}
where $T^*$ is any constant satisfying $0<T^*< T$.
\end{tm}
{\sc Proof.} 
First we observe from (\ref{eq:4.90}) that 
$\varepsilon_n$ is exactly of order $n^{-q/(2q+1)}$ as $n\rightarrow\infty$.
We observe from Proposition \ref{pn:4.6} and  Lemma \ref{la:a.1} (in Appendix A) that
\begin{eqnarray*}
&& P_{s, r} \Big\{ 
\min\{ 
\frac{1}{20 \int_0^T s(t) e^{-\int_0^t s(u) du} dt} , \frac{1}{200} \}
[ \int_0^T  | \hat{s}_n (t) - s(t) | 
e^{-\int_0^t s(u) du } dt ]^2
\nonumber \\
&& \hspace{0.5cm} \leq 
[ 6 + \frac{ 2 \log (2) }{(1- e^{-1})^2} +
8 \max\{1, \log( \frac{ e^{ \bar{\kappa}^4 (\bar{\kappa}^4 +1) T/2} }{ \varepsilon_n \delta_n^{2 n_\theta } } )   \} ] \varepsilon_n^2 \Big\}
\nonumber \\
&\geq &
1 - 4 \exp[ - \frac{ n \varepsilon_n^2 }{ 2^7 (250) } ]
- \exp[ - n( \frac{ \varepsilon_n^2 }{16} - \delta^\dag_n )],
\end{eqnarray*}
or equivalently,
\begin{eqnarray*}
&& P_{s,r} \Big\{ 
\int_0^T  | \hat{s}_n (t) - s(t) | 
e^{-\int_0^t s(u) du } dt 
\leq 
\Big\{ \max\{ 
20 \int_0^T s(t) e^{-\int_0^t s(u) du} dt, 200 \}
\nonumber \\
&&\hspace{0.5cm}\times
[ 6 + \frac{ 2 \log (2) }{(1- e^{-1})^2} +
8 \max\{1, \log( \frac{ e^{\bar{\kappa}^4 (\bar{\kappa}^4 +1) T/2}  }{ \varepsilon_n \delta_n^{2 n_\theta } } )   \} ] \varepsilon_n^2 
\Big\}^{1/2} \Big\}
\nonumber \\
&\geq &
1 - 4 \exp[ - \frac{ n \varepsilon_n^2 }{ 2^7 (250) } ]
- \exp[ - n( \frac{ \varepsilon_n^2 }{16} - \delta^\dag_n )].
\end{eqnarray*}
This implies that
\begin{eqnarray*}
&& E_{s,r} [ \int_0^T  | \hat{s}_n (t) - s(t) | 
e^{-\int_0^t s(u) du } dt ]
\nonumber \\
&\leq & 
\Big\{ \max\{ 
20 \int_0^T s(t) e^{-\int_0^t s(u) du} dt, 200 \}
\nonumber \\
&&\hspace{0.5cm}\times
[ 6 + \frac{ 2 \log (2) }{(1- e^{-1})^2} +
8 \max\{1, \log( \frac{ e^{\bar{\kappa}^4 (\bar{\kappa}^4 +1) T/2}  }{ \varepsilon_n \delta_n^{2 n_\theta } } )   \} ] \varepsilon_n^2 
\Big\}^{1/2}
\nonumber \\
&&
+ \kappa_0^2 T
 \{ 4 \exp[ - \frac{ n \varepsilon_n^2 }{ 2^7 (250) } ]
+ \exp[ - n( \frac{ \varepsilon_n^2 }{16} - \delta^\dag_n )] \},
\end{eqnarray*}
and consequently
\begin{eqnarray}
&& E_{s,r} [ \int_0^T  | \hat{s}_n (t) - s(t) | dt ]
\nonumber \\
&\leq & 
\Big\{
\max\{ 20 \kappa_0^2 T e^{2 \kappa_0^2 T}, 200 e^{2 \kappa_0^2 T} \}
[ 6 + \frac{ 2 \log (2) }{(1- e^{-1})^2} +
8 \max\{1, \log[ \frac{ e^{ \bar{\kappa}^4 (\bar{\kappa}^4 +1) T/2}  }{ \varepsilon_n \delta_n^{2 n_\theta } } ]   \} ] \varepsilon_n^2 
\Big\}^{1/2}
\nonumber \\
&&
+ \kappa_0^2 T e^{\kappa_0^2 T} 
 \{ 4 \exp[ - \frac{ n \varepsilon_n^2 }{ 2^7 (250) } ]
+ \exp[ - n( \frac{ \varepsilon_n^2 }{16} - \delta^\dag_n )] \}
\label{eq:3.77} \\
&=& O( n^{-q/(2q+1)} \log^{1/2} n ),
\nonumber
\end{eqnarray}
as $n\rightarrow \infty$.
Next we assume, in addition, that $s(t) >0$ for all $t\in [0, T]$.
Let $\xi(t), t\in [0, T)$ be as in (\ref{eq:a.63}).
Since  $s\in \Theta_{\tilde{\kappa}, q}$ and
\begin{displaymath}
s(t) e^{-\int_0^t s(u) du} \leq \xi(t) \leq \max\{ s(t), s(t) r(u): u \in [0, T)\},
\end{displaymath}
we have $0 < \min_{0\leq t< T} \xi(t) \leq \max_{0\leq t< T} \xi(t) \leq \bar{\kappa}^4$.
Thus as in the previous case,
\begin{eqnarray*}
&& P_{s,r} \Big\{ 
\min\{ 
\frac{1}{20 \int_0^T \int_0^t \xi( t-u) s(t) r(u) e^{-\int_{t-u}^t s(v) r( v-t+u) dv} du dt} , \frac{1}{200} \}
\nonumber \\
&& \hspace{0.5cm}\times
[ \int_0^T \int_0^t | \hat{s}_n (t) \hat{r}_n (u) - s(t) r(u) |
\xi( t-u) e^{- \int_{t-u }^t s(v) r(v-t+u ) dv } du dt ]^2
\nonumber \\
&& \hspace{0.5cm} \leq 
[ 6 + \frac{ 2 \log (2) }{(1- e^{-1})^2} +
8 \max\{1, \log( \frac{ e^{ \bar{\kappa}^4 (\bar{\kappa}^4 +1) T/2} }{ \varepsilon_n \delta_n^{2 n_\theta } } )  \} ] \varepsilon_n^2
\Big\}
\nonumber \\
&\geq &
1 - 4 \exp[ - \frac{ n \varepsilon_n^2 }{ 2^7 (250) } ]
- \exp[ - n( \frac{ \varepsilon_n^2 }{16} - \delta^\dag_n )],
\end{eqnarray*}
or equivalently,
\begin{eqnarray*}
&& P_{s, r} \Big\{ 
\int_0^T \int_0^t | \hat{s}_n (t) \hat{r}_n (u) - s(t) r(u) |
\xi( t-u) e^{- \int_{t-u }^t s(v) r(v-t+u ) dv } du dt
\nonumber \\
&& \hspace{0.5cm} \leq 
\Big\{ \max \{ 
20 \int_0^T \int_0^t \xi( t-u) s(t) r(u) e^{-\int_{t-u}^t s(v) r( v-t+u) dv} du dt, 200 \}
\nonumber \\
&&\hspace{0.5cm}\times
[ 6 + \frac{ 2 \log (2) }{(1- e^{-1})^2} +
8 \max\{1, \log( \frac{ e^{ \bar{\kappa}^4 (\bar{\kappa}^4 +1) T/2} }{ \varepsilon_n \delta_n^{2 n_\theta } } ) \} ] \varepsilon_n^2
\Big\}^{1/2} \Big\}
\nonumber \\
&\geq &
1 - 4 \exp[ - \frac{ n \varepsilon_n^2 }{ 2^7 (250) } ]
- \exp[ - n( \frac{ \varepsilon_n^2 }{16} - \delta^\dag_n )].
\end{eqnarray*}
This implies that
\begin{eqnarray*}
&& E_{s,r} [ \int_0^T \int_0^t | \hat{s}_n (t) \hat{r}_n (u) - s(t) r(u) |
\xi( t-u) e^{- \int_{t-u }^t s(v) r(v-t+u ) dv } du dt ]
\nonumber \\
&\leq &
\Big\{ \max \{ 
20 \int_0^T \int_0^t \xi( t-u) s(t) r(u) e^{-\int_{t-u}^t s(v) r( v-t+u) dv} du dt, 200 \}
\nonumber \\
&&\hspace{0.5cm}\times
[ 6 + \frac{ 2 \log (2) }{(1- e^{-1})^2} +
8 \max\{1, \log( \frac{ e^{ \bar{\kappa}^4 (\bar{\kappa}^4 +1) T/2} }{ \varepsilon_n \delta_n^{2 n_\theta } } ) \} ] \varepsilon_n^2
\Big\}^{1/2}
\nonumber \\
&&
+ \bar{\kappa}^8 T^2
\{ 4 \exp[ - \frac{ n \varepsilon_n^2 }{ 2^7 (250) } ]
+ \exp[ - n( \frac{ \varepsilon_n^2 }{16} - \delta^\dag_n )] \}
\nonumber \\
&\leq &
\Big\{ \max \{ 
20 \bar{\kappa}^8 T^2, 200 \}
[ 6 + \frac{ 2 \log (2) }{(1- e^{-1})^2} +
8 \max\{1, \log( \frac{ e^{ \bar{\kappa}^4 (\bar{\kappa}^4 +1) T/2} }{ \varepsilon_n \delta_n^{2 n_\theta } } ) \} ] \varepsilon_n^2
\Big\}^{1/2}
\nonumber \\
&&
+ \bar{\kappa}^8 T^2
\{ 4 \exp[ - \frac{ n \varepsilon_n^2 }{ 2^7 (250) } ]
+ \exp[ - n( \frac{ \varepsilon_n^2 }{16} - \delta^\dag_n )] \},
\end{eqnarray*}
and
\begin{eqnarray}
&& E_{s,r} [ \int_0^T \int_0^t | \hat{s}_n (t) \hat{r}_n (u) - s(t) r(u) |
du dt ]
\nonumber \\
&\leq &
\Big\{ \frac{ \max \{ 
20 \bar{\kappa}^8 T^2 e^{ 2 \bar{\kappa}^4 T}, 200 e^{ 2 \bar{\kappa}^4 T} \} }{ \min_{0\leq t< T} \xi^2 (t) }
[ 6 + \frac{ 2 \log (2) }{(1- e^{-1})^2} +
8 \max\{1, \log( \frac{ e^{ \bar{\kappa}^4 (\bar{\kappa}^4 +1) T/2} }{ \varepsilon_n \delta_n^{2 n_\theta } } ) \} ] \varepsilon_n^2
\Big\}^{1/2}
\nonumber \\
&&
+\frac{ \bar{\kappa}^8 T^2 e^{ \bar{\kappa}^4 T} }{ \min_{0\leq t <T} \xi(t) }
\{ 4 \exp[ - \frac{ n \varepsilon_n^2 }{ 2^7 (250) } ]
+ \exp[ - n( \frac{ \varepsilon_n^2 }{16} - \delta^\dag_n )] \}
\label{eq:3.78} \\
&=& O( n^{-q/(2q+1)} \log^{1/2} n),
\nonumber
\end{eqnarray}
as $n\rightarrow\infty$.
Hence
\begin{eqnarray*}
&& [\min_{0\leq t< T} s(t) ] E_{s,r} [  \int_0^T |\hat{r}_n (u) - r(u)| (T-u) du]
\nonumber \\
&\leq &
E_{s,r} [ \int_0^T \int_u^T s(t) |\hat{r}_n (u) - r(u) | dt du ]
\nonumber \\
&\leq &
E_{s,r} [ \int_0^T |\hat{s}_n (t) - s(t) | \int_0^t \hat{r}_n (u) du dt ]
+ E_{s,r} [ \int_0^T \int_0^t | \hat{s}_n (t) \hat{r}_n (u) - s(t) r(u) |
du dt ]
\nonumber \\
&\leq &
\bar{\kappa}^2 T E_{s,r} [ \int_0^T |\hat{s}_n (t) - s(t) | dt ]
+ E_{s,r} [ \int_0^T \int_0^t | \hat{s}_n (t) \hat{r}_n (u) - s(t) r(u) |
du dt ].
\end{eqnarray*}
Thus we conclude from (\ref{eq:3.77}) and (\ref{eq:3.78}) that
\begin{eqnarray*}
&& 
E_{s,r} [\int_0^{T^*} | \hat{r}_n (u) - r(u) | du ]
\nonumber \\
&\leq & 
\frac{ \bar{\kappa}^2 T}{ (T- T^*) \min_{0\leq t< T} s(t) } E_{s,r} [ \int_0^T |\hat{s}_n (t) - s(t) | dt ]
\nonumber \\
&&
+ \frac{1}{(T-T^*) \min_{0\leq t< T} s(t) } E_{s,r} [ \int_0^T \int_0^t | \hat{s}_n (t) \hat{r}_n (u) - s(t) r(u) |
du dt ]
\nonumber \\
&=& O( n^{-q/(2q +1)} \log^{1/2} n),
\end{eqnarray*}
as $n\rightarrow\infty$.
This proves the theorem.\hfill $\Box$

\section{Lower bounds}

Suppose $N(.)$ is a counting process with conditional intensity $\lambda_1 (.|.)$ as given by (\ref{eq:1.1}).
Let $N_1(.), \cdots, N_n(t), t\in [0, T),$  be independent identically distributed copies of $N(t), t\in [0, T)$.
In this section we shall compute lower bounds on the rate of convergence of an estimator for $s$ and
an estimator for $r$ based on $N_1(t), \cdots, N_n(t), t\in [0, T)$.

Let $0< \theta < T$  be as in (\ref{eq:4.86}) and 
$\Theta_{\tilde{\kappa}, q}$ be as in Section 2. Define 
\begin{eqnarray*}
\Theta_{\theta, \tilde{\kappa}, q }
&=& \Big\{ f= g^2:
g\in {\cal C}^{q_0} [0, T), \min_{t\in [0, T)} g(t) \geq 0, \max_{t\in [0, T)} |\frac{ d^j}{dt^j} g (t)| < \kappa_j, j=0,\cdots, q_0, 
\nonumber \\
&&\hspace{0.5cm}
|\frac{ d^{q_0}}{ dt^{q_0}} g (t_1) - \frac{ d^{q_0}}{ dt^{q_0}} g (t_2) | \leq \kappa_{q_0+1} |t_1-t_2|^{q_1},
\forall t_1,t_2\in [0, T),
\mbox{$g(t) = 0$ if $t\in [0, \theta]$} \Big\}.
\end{eqnarray*}

\begin{la} \label{la:3.1}
Let $\tilde{\Theta}_{\tilde{\kappa}, q, n}\subseteq \Theta_{\tilde{\kappa}, q}$
such that {\rm card}$(\tilde{\Theta}_{\tilde{\kappa}, q, n}) < \infty$.
Suppose that $\tilde{s}_n$ is an estimator for $s$ based on 
$N_1(t), \cdots, N_n(t), t\in [0, T)$.
Then
\begin{eqnarray}
&& \sup \{ E_{s, r} [ \int_0^T | \tilde{s}_n (t) - s(t) | dt]: s \in \Theta_{\tilde{\kappa}, q},
r \in \Theta_{\theta, \tilde{\kappa}, q}  \} 
\nonumber \\
&\geq &
\frac{1}{2} \inf \{ \int_0^T |s_1(t) - s_2(t) | dt: s_1\neq s_2,
s_1, s_2 \in \tilde{\Theta}_{ \tilde{\kappa}, q, n} \}
\Big\{ 1 - \frac{1}{ 
\log[ {\rm card} (\tilde{\Theta}_{\tilde{\kappa}, q, n} )-1 ] }
\Big[ \log 2 
\nonumber \\
&&\hspace{1.0cm}
+ \frac{1}{[{\rm card} (\tilde{\Theta}_{\tilde{\kappa}, q, n} ) ]^2 } \sum_{s_1, s_2 \in \tilde{\Theta}_{\tilde{\kappa}, q, n} }
\sum_{i=1}^n E_{s_1, r_1} \log \frac{ p_{s_1, r_1} (\{ w_{i,1},\cdots, w_{i, N_i(T)} \}) }{
p_{s_2, r_1} (\{ w_{i,1},\cdots, w_{i, N_i(T)} \}) } \Big] \Big\},
\label{eq:3.5}
\end{eqnarray}
for any $r_1\in \Theta_{\theta, \tilde{\kappa}, q}$.
Next suppose that $\tilde{r}_n$ is an estimator for $r$ based on 
$N_1(t), \cdots, N_n(t), t\in [0, T)$.
Let $T^*$ be a constant satisfying $\theta< T^* < T$ and $\tilde{\Theta}_{\theta, T^*, \tilde{\kappa}, q, n}
\subset \Theta_{\theta, \tilde{\kappa}, q}$ such that {\rm card}$(\tilde{\Theta}_{\theta, T^*, \tilde{\kappa}, q, n}) <\infty$
and $r_1(u) = r_2(u), u\in [T^*, T)$ $\forall r_1, r_2 \in \tilde{\Theta}_{\theta, T^*, \tilde{\kappa}, q, n}$.
Then
\begin{eqnarray}
&& \sup \{ E_{s, r} [ \int_0^{T^*} | \tilde{r}_n (t) - r(t) | dt]: s \in \Theta_{\tilde{\kappa}, q},
r \in \Theta_{\theta, \tilde{\kappa}, q}  \} 
\nonumber \\
&\geq &
\frac{1}{2} \inf \{ \int_0^{T^*} |r_1(t) - r_2(t) | dt: r_1\neq r_2,
r_1, r_2 \in \tilde{\Theta}_{\theta, T^*, \tilde{\kappa}, q, n} \}
\nonumber \\
&&\hspace{0.5cm}\times
\Big\{ 1 - \frac{1}{ 
\log [ {\rm card} (\tilde{\Theta}_{\theta, T^*, \tilde{\kappa}, q, n} )-1 ] }
\Big[ \log 2 
\nonumber \\
&&\hspace{1.0cm}
+ \frac{1}{[ {\rm card} (\tilde{\Theta}_{\theta, T^*, \tilde{\kappa}, q, n} ) ]^2 } 
\sum_{r_1, r_2 \in \tilde{\Theta}_{\theta, T^*, \tilde{\kappa}, q, n} }
\sum_{i=1}^n E_{s_1, r_1} \log \frac{ p_{s_1, r_1} (\{ w_{i,1},\cdots, w_{i, N_i(T)} \}) }{
p_{s_1, r_2} (\{ w_{i,1},\cdots, w_{i, N_i(T)} \}) } \Big] \Big\},
\label{eq:3.6}
\end{eqnarray}
for any $s_1\in \Theta_{\tilde{\kappa}, q}$.
\end{la}
We refer the reader to Appendix A for a proof of Lemma \ref{la:3.1}.
Theorems \ref{tm:4.1} and \ref{tm:4.2} (below) are the main results of this section.
They are motivated by the lower bound results in Yatracos (1988).

\begin{tm} \label{tm:4.1}
Let $q>0$. Suppose that $\tilde{s}_n$ is an estimator for $s$ based on 
$N_1(t), \cdots, N_n(t), t\in [0, T)$.
Then there exists a constant $C_{\tilde{\kappa},q}>0$ (depending only on $\tilde{\kappa}$ and $q$) such that
\begin{displaymath}
\sup \{ E_{s, r} [ \int_0^T | \tilde{s}_n (t) - s(t) | dt]: s \in \Theta_{\tilde{\kappa}, q},
r \in \Theta_{\theta, \tilde{\kappa}, q}  \} 
\geq 
C_{\tilde{\kappa}, q} n^{-q/(2q+1)}.
\end{displaymath}
\end{tm}
{\sc Proof.}
Let $\{b_n >0: n=1, 2, \cdots\}$ be a sequence of constants that tend to 0 as $n\rightarrow \infty$
and that $b_n^{-1}$ is an integer.
For $i=1,\cdots, b_n^{-1}$, define $\phi_{i,n}: [0, T) \rightarrow \mathbb{R}$ by 
\begin{displaymath}
\phi_{i,n} (t) = \left\{ \begin{array}{ll}
(b_n T)^q [ 1 - (\frac{ 2 t - (2i-1) b_n T}{ b_n T} )^2 ]^q, & \mbox{if $(i-1) b_n T \leq t< i b_n T$,}
\\
0, & \mbox{otherwise.}
\end{array}
\right.
\end{displaymath}
Writing $q=q_0+q_1$, $q_0$ a nonnegative integer and $0<q_1\leq 1$, we have
\begin{eqnarray*}
\lim_{n\rightarrow\infty} \max_{t\in [0, T)} |\frac{ d^j }{dt^j} \phi_{i,n} (t) | &<& \infty, \hspace{0.5cm} \forall j=0, \cdots, q_0,
\nonumber \\
\lim_{n\rightarrow\infty} \max_{t_1\neq t_2\in [0, T)} |\frac{ d^{q_0} }{dt^{q_0}} \phi_{i,n} (t_1) 
- \frac{ d^{q_0} }{dt^{q_0}} \phi_{i,n} (t_2)|/|t_1-t_2|^{q_1} &<& \infty.
\end{eqnarray*}
Let $\Xi_{a,n}$ denote functions of the form 
\begin{displaymath}
a [ 1 + \sum_{i=1}^{b_n^{-1}} \gamma_i \phi_{i,n} (t) ]^2, \hspace{0.5cm}\forall t\in [0, T),
\end{displaymath}
where
$\gamma_i = 0$ or $1$ and $a>0$ is a suitably small constant such that 
$\Xi_{a,n} \subset \Theta_{\tilde{\kappa}, q}$.
If $s_1, s_2 \in \Xi_{a,n}$ where $s_1 \neq s_2$, then
writing
\begin{equation}
s_1 (t) =
a [ 1 + \sum_{i=1}^{b_n^{-1}} \gamma_{1,i} \phi_{i,n} (t) ]^2,
\hspace{0.5cm}
s_2 (t) =
a [ 1 + \sum_{i=1}^{b_n^{-1}} \gamma_{2,i} \phi_{i,n} (t) ]^2,
\label{eq:3.1}
\end{equation}
with $\gamma_{1,i}, \gamma_{2,i}$ taking values $0$ or $1$, we have
\begin{eqnarray*}
\int_0^T | s_1 (t) - s_2 (t) | dt
&=& 
a \int_0^T | 2 \sum_{i=1}^{b_n^{-1}} (\gamma_{1,i} -\gamma_{2, i} ) \phi_{i,n}(t) 
+ \sum_{i=1}^{b_n^{-1}} (\gamma_{1,i} -\gamma_{2, i}) \phi^2_{i,n} (t) | dt
\nonumber \\
&\geq & 
a \int_0^{b_n T} [ 2 \phi_{1,n}(t) 
+ \phi^2_{1,n} (t) ] dt
\nonumber \\
&=& a (b_n T)^{q+1} J_q + \frac{a (b_n T)^{2q+1} J_{2q} }{2},
\end{eqnarray*}
where 
\begin{equation}
J_l = \int_{-1}^1 (1- y^2)^l dy>0, \hspace{0.5cm}\forall l>0. 
\label{eq:4.77}
\end{equation}
Also it follows from (\ref{eq:3.1}) that
\begin{eqnarray*}
|\frac{ s_1 (t) - s_2(t)}{s_1(t)} |
&\leq &
2 a \phi_{i,n}(t) 
+ a \phi^2_{i,n} (t) 
\nonumber \\
&\leq & 2 a (b_n T)^q + a  (b_n T)^{2 q}, \hspace{0.5cm}\forall t\in [0, T).
\end{eqnarray*}
Let $r_1 \in \Theta_{\theta, \tilde{\kappa}, q}$.
Now using Lemma \ref{la:a.1} in Appendix A, 
\begin{eqnarray*}
&& 
E_{s_1, r_1} \log [ \frac{ p_{s_1, r_1} (\{ w_{i,1},\cdots, w_{i, N_i(T)}\} )
}{
p_{s_2, r_1} (\{ w_{i,1},\cdots, w_{i, N_i(T)}\} ) } ] 
\nonumber \\
&=&
\int_0^T  \{ \frac{s_2 (t) }{s_1 (t)} -1 - \log[ \frac{s_2 (t) }{ s_1 (t)} ] \}
s_1 (t) e^{-\int_0^t s_1 (u) du } dt
\nonumber \\
&& + \int_0^T \int_0^t \{ \frac{ s_2 (t) }{ s_1 (t) } -1 - \log[ \frac{s_2 (t) }{ s_1 (t) } ] \}
\xi( t-u) s_1 (t) r_1 (u) e^{- \int_{t-u }^t s_1 (v) r_1 (v-t+u ) dv } du dt
\nonumber \\
&\leq &
\frac{1}{2} \int_0^T (\frac{ s_1(t)-s_2(t) }{s_1(t)} )^2
s_1 (t) e^{-\int_0^t s_1 (u) du } dt
\nonumber \\
&& + \frac{1}{2} \int_0^T \int_0^t (\frac{ s_1(t) - s_2 (t) }{ s_1 (t) })^2 
\xi( t-u) s_1 (t) r_1 (u) e^{- \int_{t-u }^t s_1 (v) r_1 (v-t+u ) dv } du dt
\nonumber \\
&\leq &
 \frac{ a^2 (b_n T)^{2 q} }{2} [ 2 + (b_n T)^q ]^2 (1 + \bar{\kappa}^8 T^2),
\end{eqnarray*}
where $\bar{\kappa}= \kappa_0 \vee 1$.
Finally we observe from Proposition 3.8 of Birg\'{e} (1983) that
there exists a subset $\tilde{\Theta}_{\tilde{\kappa}, q, n}$ of $\Xi_{a, n}$ such that
\begin{displaymath}
\int_0^T |s_1(t) -s_2(t)| dt \geq \frac{ 1}{8 b_n} 
[ a (b_n T)^{q+1} J_q + \frac{a (b_n T)^{2q+1} J_{2q} }{2} ],
\hspace{0.5cm}
\forall s_1\neq s_2 \in \tilde{\Theta}_{\tilde{\kappa}, q, n},
\end{displaymath}
and $\log[ {\rm card}(\tilde{\Theta}_{\tilde{\kappa}, q, n}) -1] > 0.316/b_n$.
Consequently we conclude from 
(\ref{eq:3.5}) that
\begin{eqnarray*}
&& \sup \{ E_{s, r} [ \int_0^T | \tilde{s}_n (t) - s(t) | dt]: s \in \Theta_{\tilde{\kappa}, q},
r \in \Theta_{\theta, \tilde{\kappa}, q }  \} 
\nonumber \\
&\geq &
\frac{ 1}{16 b_n} 
[ a (b_n T)^{q+1} J_q + \frac{a (b_n T)^{2q+1} J_{2q} }{2} ]
\Big\{ 1 
\nonumber \\
&&\hspace{0.5cm}
- \frac{b_n}{0.316} \Big[ \log 2
+ \frac{ a^2 n (b_n T)^{2 q} }{2} [ 2 + (b_n T)^q ]^2 (1 + \bar{\kappa}^8 T^2) \Big] \Big\}
\nonumber \\
&=&
\frac{ a b_n^q T^{q+1} }{16 } 
[ J_q + \frac{ (b_n T)^q J_{2q} }{2} ]
\Big\{ 1 
- \frac{b_n}{0.316} \Big[ \log 2
+ \frac{ a^2 n (b_n T)^{2 q} }{2} [ 2 + (b_n T)^q ]^2 (1 + \bar{\kappa}^8 T^2) \Big] \Big\}.
\end{eqnarray*}
Thus we conclude that
there exist strictly positive constants $C_0$ and $C_{ \tilde{\kappa}, q}$ 
(depending only on $\tilde{\kappa}$ and $q$) such that
by taking $b_n = 1/\lceil C_0 n^{1/(2 q+1)} \rceil$, we have
\begin{displaymath}
\sup \{ E_{s, r} [ \int_0^T | \tilde{s}_n (t) - s(t) | dt]: s \in \Theta_{\tilde{\kappa}, q},
r \in \Theta_{\theta, \tilde{\kappa}, q}  \} 
\geq  C_{\tilde{\kappa}, q} n^{-q/(2q+1)}.
\end{displaymath}
This proves Theorem \ref{tm:4.1}.\hfill $\Box$

\begin{tm} \label{tm:4.2}
Let $q>0$ and $\theta, T^*$ be constants satisfying $0< \theta< T^* < T$.
Suppose that $\tilde{r}_n$ is an estimator for $r$ based on 
$N_1(t), \cdots, N_n(t), t\in [0, T)$.
Then there exists a constant $C_{\theta, \tilde{\kappa}, q}>0$ (depending only on $\theta, \tilde{\kappa}$ and $q$)
such that
\begin{displaymath}
\sup \{ E_{s, r} [ \int_0^{T^*} | \tilde{r}_n (u) - r(u) | du]: s \in \Theta_{\tilde{\kappa}, q},
r \in \Theta_{\theta, \tilde{\kappa}, q }  \} 
\geq 
C_{\theta, \tilde{\kappa}, q} n^{-q/(2q +1)}.
\end{displaymath}
\end{tm}
{\sc Proof.}
Let $f\in \Theta_{\theta, \tilde{\kappa}, q}$ such that $f(u) >0$ if $u\in [\theta^*, T^*]$  for
constants $\theta^*, T^*$ satisfying
$0< \theta < \theta^* < T^* < T$.
Then $0< \underline{f} := \min_{u\in [\theta^*, T^*] } f(u) \leq \bar{f} := \max_{u\in [\theta^*, T^* ]} f(u) <\infty$.
Let $\{ b_n >0: n=1, 2, \cdots\}$ be a sequence of constants that tend to 0 as $n\rightarrow \infty$
such that $b_n^{-1}$ is an integer.
Then for $i=1,\cdots, b_n^{-1}$, define $\phi_{i,n}: [0, T) \rightarrow \mathbb{R}$ 
by
\begin{eqnarray*}
\phi_{i,n} (u) 
&= & \left\{ \begin{array}{ll}
[ b_n (T^* - \theta^*) ]^q
[ 1 - (\frac{ 2 u - 2 \theta^* - (2i-1) b_n (T^* - \theta^*) }{ b_n (T^* - \theta^*) } )^2 ]^q,  & \\
\hspace{0.5cm}  \mbox{if $\theta^* + (i-1) b_n (T^* -\theta^*) \leq u < \theta^* + i b_n (T^* - \theta^*)$,}
\\
0,\hspace{0.5cm} \mbox{elsewhere}. &
\end{array}
\right.
\end{eqnarray*}
Let $\Xi_{a,n}$ denote functions of the form 
\begin{displaymath}
a [ f^{1/2} (u) + \sum_{i=1}^{b_n^{-1}} \gamma_i \phi_{i,n} (u) ]^2,
\hspace{0.5cm}\forall u\in [0, T),
\end{displaymath}
where
$\gamma_i = 0$ or $1$ and $a>0$ is a sufficiently small constant such that
$\Xi_{a, n} \subset \Theta_{\theta, \tilde{\kappa}, q }$.
If $r_1, r_2 \in \Xi_{a, n}$ where $r_1 \neq r_2$, then
writing
\begin{eqnarray}
r_1 (u)&=& 
a [ f^{1/2} (u) + \sum_{i=1}^{b_n^{-1}} \gamma_{1,i} \phi_{i,n} (u) ]^2,
\nonumber \\
r_2 (u)&=& 
a [ f^{1/2} (u) + \sum_{i=1}^{b_n^{-1}} \gamma_{2,i} \phi_{i,n} (u) ]^2,
\label{eq:3.11}
\end{eqnarray}
with $\gamma_{1,i}, \gamma_{2,i}$ taking values $0$ or $1$, we have
\begin{eqnarray*}
&& \int_0^{T^*} | r_1 (u) - r_2 (u) | du
\nonumber \\
&=& 
a \int_{\theta^*}^{T^*} | 2 f^{1/2} (u) \sum_{i=1}^{b_n^{-1}} (\gamma_{1,i} -\gamma_{2,i}) \phi_{i,n} (u)
+ \sum_{i=1}^{b_n^{-1}} (\gamma_{1,i} -\gamma_{2,i}) \phi_{i,n}^2 (u) | du
\nonumber \\
&\geq & 
a \int_{\theta^*}^{\theta^*+ b_n (T^* -\theta^*)} [ 2 \underline{f}^{1/2} \phi_{1,n} (u)
+ \phi_{1,n}^2 (u) | du
\nonumber \\
&= & 
a \underline{f}^{1/2} [b_n (T^* -\theta^*)]^{q+1} J_q 
+ \frac{ a [ b_n (T^* - \theta^*)]^{2q+1} J_{2q} }{ 2},
\end{eqnarray*}
where $J_q$ is as in (\ref{eq:4.77}).
Also it follows from (\ref{eq:3.11}) that
\begin{eqnarray*}
|\frac{ r_1 (u) - r_2(u)}{r_1(u)} |
&\leq &
\frac{
2 a \bar{f}^{1/2} \phi_{i,n}(u) 
+ a \phi^2_{i,n} (u)
}{ \underline{f}^{1/2} }
\nonumber \\
&\leq & \frac{ 2 a \bar{f}^{1/2} [b_n (T^* - \theta^*)]^q }{\underline{f}^{1/2} } 
+ \frac{ a [ b_n (T^* - \theta^*)]^{2 q} }{\underline{f}^{1/2} }, \hspace{0.5cm}\forall u\in [\theta^*, T^*].
\end{eqnarray*}
Let $s_1 \in \Theta_{\tilde{\kappa}, q}$. Now using Lemma \ref{la:a.1}, 
\begin{eqnarray*}
&& 
E_{s_1, r_1} \log [ \frac{ p_{s_1, r_1} (\{ w_{i,1},\cdots, w_{i, N_i(T)}\} )
}{
p_{s_1, r_2} (\{ w_{i,1},\cdots, w_{i, N_i(T)}\} ) } ] 
\nonumber \\
&=&
\int_0^T \int_0^t \{ \frac{ r_2 (u) }{ r_1 (u) } -1 - \log[ \frac{r_2 (u) }{ r_1 (u) } ] \}
\xi( t-u) s_1 (t) r_1 (u) e^{- \int_{t-u }^t s_1 (v) r_1 (v-t+u ) dv } du dt
\nonumber \\
&= &
\int_0^T \{ \frac{ r_2 (u) }{ r_1 (u) } -1 - \log[ \frac{r_2 (u) }{ r_1 (u) } ] \}
\int_u^T \xi( t-u) s_1 (t) r_1 (u) e^{- \int_{t-u }^t s_1 (v) r_1 (v-t+u ) dv } dt du
\nonumber \\
&\leq &
\frac{1}{2} \int_{\theta^*}^{T^*} (\frac{ r_1(u) - r_2(u) }{ r_1(u)})^2
\int_u^T \xi( t-u) s_1 (t) r_1 (u) e^{- \int_{t-u }^t s_1 (v) r_1 (v-t+u ) dv } dt du
\nonumber \\
&\leq &
\frac{\bar{\kappa}^8 b_n^{2 q} (T - \theta^*)^2 
[ 2 a \bar{f}^{1/2} (T^* - \theta^*)^q 
+ a b_n^q (T^* - \theta^*)^{2 q}]^2 }{2 \underline{f} },
\end{eqnarray*}
where $\bar{\kappa} = \kappa_0 \vee 1$.
Finally we observe from Proposition 3.8 of Birg\'{e} (1983) that
there exists a subset $\tilde{\Theta}_{\theta, T^*, \tilde{\kappa}, q, n}$ of $\Xi_{a, n}$ such that
\begin{eqnarray*}
&& \int_0^{T^*} |r_1(u) -r_2(u)| du 
\nonumber \\
& \geq &
\frac{ 1}{8 b_n} 
\Big[ a \underline{f}^{1/2} [b_n (T^* -\theta^*)]^{q+1} J_q 
+ \frac{ a [ b_n (T^* - \theta^*)]^{2q+1} J_{2q} }{ 2} \Big],
\hspace{0.5cm}
\forall r_1\neq r_2 \in \tilde{\Theta}_{\theta, T^*, \tilde{\kappa}, n},
\end{eqnarray*}
and $\log[ {\rm card}(\tilde{\Theta}_{\theta, T^*, \tilde{\kappa}, n}) -1] > 0.316/b_n$.
Consequently we conclude from 
(\ref{eq:3.6}) that
\begin{eqnarray*}
&& \sup \{ E_{s, r} [ \int_0^{T^*} | \tilde{r}_n (u) - r(u) | du]: s \in \Theta_{\tilde{\kappa}, q},
r \in \Theta_{\theta, \tilde{\kappa}, q }  \} 
\nonumber \\
&\geq &
\frac{ a b_n^q }{16 } 
[ \underline{f}^{1/2} (T^* -\theta^*)^{q+1} J_q 
+ \frac{ b_n^q (T^* - \theta^*)^{2q+1} J_{2q} }{ 2} ]
\Big\{ 1 
\nonumber \\
&&\hspace{0.5cm}
- \frac{b_n}{0.316} \Big[ \log 2
+ \frac{ a \bar{\kappa}^8 n b_n^{2 q} (T - \theta^*)^2 
[ 2 \bar{f}^{1/2} (T^* - \theta^*)^q
+ b_n^q (T^* - \theta^*)^{2 q} ]^2 }{2 \underline{f} }
\Big] \Big\}.
\end{eqnarray*}
Hence there exist strictly positive constants $C_1$ and $C_{\theta, \tilde{\kappa}, q}$ 
(depending only on $\theta, \tilde{\kappa}$ and $q$) such that
by taking $b_n = 1/ \lceil C_1 n^{1/(2q+1)} \rceil$, we have
\begin{displaymath}
\sup \{ E_{s, r} [ \int_0^{T^*} | \tilde{r}_n (u) - r(u) | du]: s \in \Theta_{\tilde{\kappa}, q},
r \in \Theta_{\theta, \tilde{\kappa}, q}  \} 
\geq  C_{\theta, \tilde{\kappa}, q} n^{-q/(2q+1)}.
\end{displaymath}
This proves Theorem \ref{tm:4.2}.\hfill $\Box$

\section{Template matching with continuous kernels}  
 
In the second part of this article, let $\bw =(\bw^{(1)},\ldots,\bw^{(d)})$ be the spike train pattern of an 
assembly of $d$ neurons recorded when an experimental stimulus
is provided to a subject, where $\bw^{(i)} = \{ w^{(i)}_1,\ldots,w^{(i)}_{N_i(T)} \}$ 
are the spike times of the $i$th neuron over the period 
$[0,T)$. The same neurons of the subject are subsequently 
observed for a longer period of time when it is engaged in other activities
and the corresponding spike trains $\by = (\by^{(1)},\ldots,\by^{(d)})$
are checked for occurrences of the template $\bw$. 

For $t \geq 0$, let $\by_t = (\by_t^{(1)},
\ldots,\by_t^{(d)})$, where $\by_t^{(i)} = \{ y-t: y \in \by^{(i)} \cap
[t,t+T) \}$. There are various algorithms in the neuroscience literature
that have been used to determine if there
is a close match between $\by_t$ and $\bw$. In Gr\"{u}n, Diesmann and Aertsen
(2002), $T$ is chosen small and a match is declared if 
$$\{ 1\leq i\leq d : \bw^{(i)} = \emptyset \} = \{ 1\leq i \leq d : \by_t^{(i)} =\emptyset \}.
$$

In the sliding sweeps algorithm [cf.\ Dayhoff and Gerstein (1983)  
and N\'{a}dasdy {\em et al.} (1999)],
a match is declared if 
$$\sup_{1 \leq i \leq d} \sup_{w \in \bw^{(i)}}
\inf_{y-t \in \by_t^{(i)}} |y-t-w| \leq \Delta,
$$
where $\Delta > 0$ is a pre-determined constant.
We shall study in this section the pattern filtering
algorithm [cf. Chi, Rauske and 
Margoliasch (2003)], which uses a scoring system to measure the proximity 
between $\bw$ and $\by_t$. 

Let $f$ be a non-increasing and non-constant
function on $[0,\infty)$ with $f(0) > 0$. The score between 
$\bw$ and $\by_t$ is given by 
\begin{equation}
S_t = \sum_{i=1}^d S_t^{(i)}, \quad {\rm where} \ekp
S_t^{(i)} = T^{-1} \sum_{y-t \in \by_t^{(i)}} \max_{w \in \bw^{(i)}}
f(|y-t-w|). 
\label{1}
\end{equation}
For a given template $\bw$, define the kernel functions
\begin{equation}
g_\bw^{(i)}(u) = \Big[ \max_{w \in \bw^{(i)}} f(|u-w|) \Big] {\bf 1}_{\{ 0
\leq u < T \}}, \hspace{0.5cm} \forall i=1,\cdots, d.
\label{2}
\end{equation}
Then we can also express
$$
S_t^{(i)} = T^{-1} \sum_{y \in \by^{(i)}} g_\bw^{(i)}(y-t).
%\label{3}
$$
The graph of $S_t^{(i)}$ against $t$ is thus a normalized sum of the
kernels $g_\bw^{(i)}(y - \cdot)$ over all $y \in \by^{(i)}$. We declare
a match between $\by_t$ and $\bw$ to be present when the proximity score $S_t$
exceeds a pre-determined threshold
level $c$. To prevent overcounting, a match at time $t$ 
is declared to be new only if the overlap
between the time interval $[t,t+T)$ and the time interval of the previous new
match is less than $\alpha T$ for 
some constant $0 < \alpha < 1$. More specifically,
let $\sigma_1 = \inf \{ t: S_t \geq c \}$ and $\sigma_{j+1} = \inf \{ t>
\sigma_j+(1-\alpha)T: S_t \geq c \}$ for $j \geq 1$. Then the number of
new matches between the spike trains $\by$ over the time interval $[0,a+T)$ and the
template $\bw$ is $U_a := \sup \{ j: \sigma_j \leq a \}$;
with the convention $U_a=0$ if $\sigma_1 > a$. 

To prevent the occurrences of too many (false) matches when $\by$ is pure noise, 
the threshold level $c$ has to be
chosen reasonably large. For $a$ large, there can be
on the average more than one new (false) match between $\by$ and
$\bw$. The Poisson distribution is often used
for modeling $U_a$ to compute the $p$-value under such circumstances.
For small $a$, the occurence of a single match would
itself be rare and we can use the probability of having at least
one match as the $p$-value. For this purpose, we study the
scan statistics 
\begin{displaymath}
M_a := \sup_{0 \leq t \leq a} S_t,
\end{displaymath}
and its dual,
the time to detection 
\begin{displaymath}
V_c := \inf \{ t: S_t \geq c \},
\end{displaymath}
which we
shall show to have asymptotic Gumbel and exponential distributions respectively.
In this section, we consider $f$ to be continuous on $[0,\infty)$ while
in Section 6, we will consider the case in which $f$ is not continuous.

\subsection{Main results}

Let $\by^{(i)}$, $i=1,\cdots,d,$ be independent Poisson
processes with constant intensity $\lambda_i > 0$. Consider
the following regularity conditions on $\bw$ and $f$. 

\smallskip
\noindent (A1) Let $\bw_*^{(1)},\cdots,\bw_*^{(d)}$ be point processes
on $[0,\infty)$ with each $\bw_*^{(i)}$ ergodic, stationary and having
non-constant inter-arrival times and let $\bw^{(i)} = \bw_*^{(i)} \cap
[0,T)$.

\smallskip
\noindent (A2) Let $f$ be continuous and let there be a possibly empty
finite set $H$ such that the second derivative of $f$ exists and is uniformly
continuous and bounded over any interval inside ${\mathbb{R}}^+ \setminus H$.
Moreover,
\begin{equation}
0 < \sup_{x \in {\mathbb{R}}^+ \setminus H} \Big| \frac{d}{dx} f(x) \Big| <
\infty \quad {\rm and} \quad \lim_{x \rightarrow \infty} f(x) > -\infty.
\label{4}
\end{equation}
Let $\mu_\bw = T^{-1} \sum_{i=1}^d \lambda_i \int_0^T g_\bw^{(i)}(u) \; du$
be the expected value of $S_t$ conditioned on $\bw$ known. 
Let the large deviation rate function of $S_t$ be
\begin{equation}
\phi_\bw(c) = \sup_{\theta > 0} \Big[ \theta c- T^{-1} \sum_{i=1}^d 
\lambda_i \int_0^T (e^{\theta g_\bw^{(i)}(u)} -1) \; du \Big] \quad {\rm for} \ c >
\mu_\bw. 
\label{5}
\end{equation}
We shall denote by $\theta_\bw$ ($=\theta_{\bw,c}$) the unique value of
$\theta > 0$ that attains the supremum on the right hand side of (\ref{5}).
By the stationarity of $\bw_*^{(i)}$ in (A1), for all $y \in
{\mathbb{R}}$, the
distribution of $\max_{w \in \bw_*^{(i)}} f(|y-w|)$ is equal to the
distribution of $Z_i := \max_{w \in \bw_*^{(i)}} f(|w|)$. Hence by the
ergodicity of $\bw^{(i)}$ in (A1) and the bounded property of $f$ in (A2), 
\begin{eqnarray}
& & T^{-1} \int_0^T e^{\theta g_\bw^{(i)}(u)} du \rightarrow Ee^{\theta Z_i}
\hspace{0.5cm} \mbox{a.s.\  $\forall \theta > 0$} \cr
& {\rm and} & T^{-1} \int_0^T g_\bw^{(i)}
(u) \ du \rightarrow E Z_i \hspace{0.5cm} \mbox{a.s.\ as $T \rightarrow \infty$}.
\label{6}
\end{eqnarray}
Let $\mu = \sum_{i=1}^d \lambda_i EZ_i$. Define the 
limiting large deviation rate function
\begin{equation}
\phi(c) = \sup_{\theta>0} \Big[ \theta c - \sum_{i=1}^d \lambda_i  
(Ee^{\theta Z_i}-1) \Big] \quad {\rm for} \ c > \mu.
\label{7}
\end{equation}
Let $\theta_*$ ($=\theta_{*,c}$) be the unique value of $\theta > 0$ attaining
the supremum on the right hand side of (\ref{7}). Then by (\ref{5}), (\ref{6}) and (\ref{7}), 
$\mu_\bw \rightarrow \mu, \phi_\bw \rightarrow \phi$ pointwise on 
$(\mu,\infty)$ and $\theta_\bw
\rightarrow \theta_*$ a.s.\  as  $T \rightarrow \infty$. 
Similarly, 
\begin{eqnarray}
v_\bw &:=& T^{-1} \sum_{i=1}^d \lambda_i \int_0^T [g_\bw^{(i)}(u)]^2
e^{\theta_\bw g_\bw^{(i)}(u)} du, 
\nonumber \\
\tau_\bw &:=& T^{-1} \sum_{i=1}^d
\lambda_i \int_0^T \Big[ \frac{d}{du} g_\bw^{(i)}(u) \Big]^2
e^{\theta_\bw g_\bw^{(i)}(u)} du
\label{8}
\end{eqnarray}
both converge almost surely to positive constants as $T \rightarrow \infty$.
Let $P_\bw$ denote the probability measure conditioned on a known $\bw$.
 
\begin{pn}
\label{l1}
Assume (A1)-(A2). Then for any $t \geq 0$, $\Delta > 0$ 
and $c > \mu$,
\begin{equation}
P_\bw \Big\{ \sup_{t < u \leq t+\Delta} S_u \geq c \Big\}
\sim \Delta \zeta_\bw e^{-T \phi_\bw(c)} \hspace{0.5cm} \mbox{a.s.\ as $T \rightarrow
\infty$},
\label{9}
\end{equation}
where $\zeta_\bw = (2 \pi)^{-1} (\tau_\bw/v_\bw)^{1/2}$.   
\end{pn}

By pieceing together the local boundary crossing 
probabilities in (\ref{9}), we are able to obtain the following results.

\begin{tm}
\label{t1}
Assume (A1)-(A2).  

{\rm (a)} Let $c > \mu$. Then the distribution (conditional on $\bw$) of 
$\zeta_\bw e^{-T \phi_\bw(c)} V_c$
converges to the exponential distribution with mean 1 almost surely as
$T \rightarrow \infty$.

{\rm (b)} Let $a \rightarrow \infty$ as $T \rightarrow \infty$ such that
$(\log a)/T$ converges to a positive constant. 
Let $c_\bw > \mu_\bw$ satisfy $\phi_\bw(c_\bw) = (\log a)/T$. Then
for any $z \in {\mathbb{R}}$, 
$$
P_\bw \{ \theta_\bw T (M_a -c_\bw) - \log \zeta_\bw \geq z \} 
\rightarrow 1-\exp(-e^{-z}) \hspace{0.5cm} \mbox{a.s.\ as $T \rightarrow \infty$}.
$$

{\rm (c)} Let $a \rightarrow \infty$ as $T \rightarrow \infty$ such that
$(\log a)/T$ converges to a positive constant. Let $c$ ($=c_T$)
be such that $\eta_\bw := a \zeta_\bw e^{-T \phi_\bw(c)}$ converges to
a constant $\eta > 0$ almost surely. Then
\begin{equation}
P_\bw \{ U_a = k \} - e^{-\eta_\bw} \frac{\eta_\bw^k}{k!}
\rightarrow 0 \hspace{0.5cm} \mbox{a.s. $\forall k=0,1,\cdots.$}
\label{11}
\end{equation}
\end{tm}

{\sc Remark.} Theorem \ref{t1} can be extended to deal with the
situation in which $m > 1$ trials are conducted, giving rise to
$m$ spike train vectors. If we
hypothesize that the times
of recurrence of the template are the same for  
the $m$ trials, then the pattern filtering algorithm is most effectively
applied by comparing $\bw$ against 
a union of the $m$ spike train vectors. If we hypothesize that 
the times of recurrence of the
template is different for each trial, then we can compare $\bw$ against
each spike train vector separately and sum up the number of new
matches for the $m$ trials. The Poisson distribution can again be used to
compute the $p$-values. Note that we do not require the vectors $\by$ to have
the same intensity for the $m$ trials. This has implications 
when one has a spike train $\by$
that is nonstationary but can be broken up into $m$ segments such that each
part has almost constant intensities as we can apply Theorem
\ref{t1}(c) separately on each of the $m$ segments.

{\sc Remark.} In Theorem 1 of Chi (2004), it was shown [without the regularity condition (A2)] that  
$$
\lim_{T \rightarrow \infty} T^{-1} \log V_c \rightarrow \phi(c) \quad {\rm a.s.
\ for \ all} \ c > \mu.
%\label{12}
$$
The question of whether $\log V_c = T \phi_\bw(c) + o(T^{1/2})$ was also
raised in a remark on page 157. Theorem \ref{t1}(a) provides a more precise answer; that
$\log V_c = T \phi_\bw(c)+O_P(1)$.

\subsection{Implementation}

We conduct a small scale simulation study in this subsection to 
test the finite sample accuracy of the analytic approximations 
in Theorem \ref{t1}. An alternative to analytic approximations
is to compute the $p$-values $p_\bw := P_\bw \{ M_a \geq c \}$
via direct Monte Carlo. However,
as $p$-values of interest are often small, a large number of simulation runs is
required for these estimations to be accurate. The 
computational cost is compounded when the time period $[0,a+T)$ of $\by^{(i)}$ is
large. 

We introduce here an importance sampling alternative for the simulation
of $p$-values. We use a change of measure argument, which is also
used in the proof of Proposition \ref{l1}, by generating 
$\by^{(i)}$ from an inhomogeneous Poisson process. Analogous change
of measures for computing $p$-values have been used in  
sequential analysis [cf.\ Siegmund (1976)], 
change-point detection [cf.\ Lai and Shan (1999)] and DNA sequence alignments
[cf.\ Chan (2003)].  

Let $P_{\theta,t}$ denote the probability measure under which $\by^{(i)}$
is generated as a Poisson point process with intensity $\eta_i(v) = 
\lambda_i e^{\theta g_\bw^{(i)}(v-t)}$ for each $1 \leq i \leq d$. Note that
$g_\bw^{(i)}(v-t)=0$ for $v \not\in [t,t+T)$ and hence the
change of measure occurs only for the generation of spikes in the
interval $[t,t+T)$. 
The likelihood of $\by_t^{(i)}$ under $P_{\theta,t}$ is given by 
$$
L_{\theta,t}(\by_t^{(i)}) = \exp \Big( - 
\lambda_i \int_0^T e^{\theta g_\bw^{(i)}(u)} du \Big) \prod_{y \in 
\by_t^{(i)}} \lambda_i e^{\theta g_\bw^{(i)}(y-t)}.
%\label{13}
$$
Hence the likelihood ratio
\begin{eqnarray}
\frac{dP_{\theta,t}}{dP_\bw}(\by) & = & \prod_{i=1}^d 
\frac{L_{\theta,t}(\by_t^{(i)})}{L_{0,t}(\by_t^{(i)})} =
\prod_{i=1}^d \exp \Big[ \theta T S_t^{(i)} - \lambda_i \int_0^T (e^{\theta
g_\bw^{(i)}(u)}-1) \; du \Big] \cr
& = & \exp \Big[ \theta T S_t 
-\sum_{i=1}^d \lambda_i \int_0^T (e^{\theta g_\bw^{(i)}(u)}-1) \ du \Big].
\label{14}
\end{eqnarray}

In our importance sampling algorithm, we first select a small $\Delta > 0$ 
such that $J:=a/\Delta$ is a positive integer. For each simulation run, we 
generate $j$ randomly from $\{ 0,\ldots, J \}$ followed by $\by$ from 
$P_{\theta_\bw,j \Delta}$. The estimate
\begin{eqnarray}
\widehat p & = & (J+1) \Big[ \sum_{j=0}^J  
\frac{dP_{\theta_\bw,j \Delta}}{dP_\bw}(\by) \Big]^{-1}
{\bf 1}_{\{ M_a \geq c \}} \cr
& = & (J+1)  \exp \Big[ \sum_{i=1}^d \lambda_i \int_0^T (e^{\theta g_\bw^{(i)}(u)}-1) \; du \Big]
\Big( \sum_{j=0}^J e^{\theta T S_{j \Delta}} \Big)^{-1} {\bf 1}_{\{ M_a \geq c \}}
\label{15}
\end{eqnarray}
is then 
unbiased for $p_\bw$. The averages of (\ref{15}) over all
the simulation runs is then the importance sampling estimate of $p_\bw$. 

{\sc Example 1.} Consider the Hamming window function
\begin{equation}
f(t) = \left\{ \begin{array}{cl} {1 \over 2} (1-\beta)+{1 \over 2} (1+\beta) \cos \left(
\frac{\pi t}{\varepsilon} \right) & {\rm if} \  0 \leq t < \varepsilon, \cr
- \beta & {\rm if} \ t \geq \varepsilon, \cr \end{array} \right.
\label{16}
\end{equation}
with $\varepsilon = 5$ ms and $\beta=0.4$ [see, for example, Chi, Rauske and Margoliash (2003)]. 

We generate a template $\bw$ over the time interval from 0 to $T=500$ ms on
$d=4$ spike trains, with interarrival distance $X$ ms between two spikes on each
spike train satisfying
\begin{equation}
P \{ X \leq x \} = 1- e^{-(x-1)^+/24}.
\label{17}
\end{equation}
This corresponds to an absolute refractory period or ``dead time'' of $1$ ms after each spike in $\bw^{(i)}$
before the next spike can be generated. In our computer experiment,
a total of 80 spikes were generated on the four spike trains using (\ref{17}).

To compute the $p$-values using direct Monte Carlo, we generated 2000 realizations of $\by$ by
using Poisson point processes with constant intensity $\lambda_i = 0.04$ ms$^{-1}$
on the interval from 0 to $a+T=20$ s. The proportion of times that
$\{ M_a \geq c \}$ occurs is taken as the estimate of $p_\bw$.
For importance sampling, 2000 simulation runs were also executed using the
algorithm described earlier by choosing
$\Delta=0.2$ ms. The following thinning method is used to generate the spike times 
in the interval $[j \Delta,
j \Delta+T)$, where the intensity of $\by^{(i)}$ is not constant under
$P_{\theta_\bw,j \Delta}$:

\smallskip
1. Let $\wtd \by_t^{(i)} = \{ u_1,\ldots,u_N \}$ be generated on $[0,T)$
as a Poisson process with constant intensity $\lambda_i e^{\theta_\bw f(0)}$.

\smallskip
2. Generate independent uniform random variables $R_1,\cdots,R_N$ on $[0,1]$
and let
$$
\by_t^{(i)} (=\{ y-t: y \in \by^{(i)} \cap [t,t+T) \})=\{ u_j \in
\wtd \by_t^{(i)}: R_j \leq e^{\theta_\bw[g_\bw(u_j)-f(0)]} \}.
%\label{18}
$$

For the analytic approximation, we apply 
Theorem \ref{t1}(a), which gives us
\begin{equation}
P_\bw \{M_a \geq c \} = P_\bw \{ V_c \leq a \} \doteq 1-\exp(-a \zeta_\bw
e^{-T \phi_\bw(c)}).
\label{19}
\end{equation}

\begin{table}
\begin{center}
{\sc Table 1.} Estimates of $P_\bw \{ M_a \geq c \} \pm$ standard error.

\begin{tabular}{c|c|c|c}
\hline
$c$ & Direct MC & Imp. Sampling & Anal. Approx. (\ref{19}) \cr
\hline
0.017 & 0.037$\pm$0.004 & 0.0387$\pm$0.0019 & 0.0383 \cr 
0.018 & 0.024$\pm$0.003 & 0.0237$\pm$0.0012 & 0.0241 \cr
0.019 & 0.016$\pm$0.003 & 0.0158$\pm$0.0008 & 0.0149 \cr
0.020 & 0.009$\pm$0.002 & 0.0095$\pm$0.0005 & 0.0091 \cr
0.021 & 0.005$\pm$0.002 & 0.0054$\pm$0.0003 & 0.0055 \cr
0.022 & 0.003$\pm$0.001 & 0.0033$\pm$0.0002 & 0.0033 \cr
\hline
\end{tabular}
\end{center}
\end{table}

We see from the results summarized in Table 1 that
there is substantial variance reduction when importance sampling 
is used. 
The analytic approximations have also been shown to be quite accurate,
lying within two standard errors of the importance sampling estimate
in all the cases considered.

\subsection{Proofs}

We preface the proofs of Proposition \ref{l1} and Theorem \ref{t1} with the
following preliminary lemmas. We shall let $\lfloor \cdot \rfloor$ denote
the greatest integer function. Let $P_{\theta_\bw,t}$ be the change of measure 
defined in the beginning
of Section 5.2 and let $P_{\theta_\bw} = P_{\theta_\bw,0}$.

\begin{la}
\label{l2}
Let $t \geq 0$ and $c >\mu_\bw$. Then 
$$
P_\bw \{ S_t \geq c \} \sim (2 \pi v_\bw)^{-1/2} \theta_\bw^{-1} T^{-1/2} e^{-T \phi_\bw(c)}
\hspace{0.5cm} \mbox{a.s.\ as $T \rightarrow \infty$}.
$$
\end{la}
{\sc Proof.} Let $E_{\theta,t}$ denote expectation with respect
to the probability measure $P_{\theta,t}$. Let 
$I_T=[z_T,z_T+\varepsilon_T)$ with $\varepsilon_T=o(T^{-1/2})$.
Then by (\ref{5}) and (\ref{14}),
\begin{eqnarray}
P_\bw \{ T^{1/2}(S_t-c) \in I_T \} & = & E_{\theta_\bw,t} \Big[ \frac{dP_{\bw}}{dP_{\theta_\bw,t}}
{\bf 1}_{\{ T^{1/2}(S_t-c) \in I_T \}} \Big] \cr
& = & e^{-T \phi_\bw(c)} E_{\theta_\bw,t} \Big[ e^{T
\theta_\bw (c-S_t)} {\bf 1}_{\{ T^{1/2}(S_t-c) \in I_T \}} \Big] \cr
& \sim & e^{-T \phi_\bw(c) -T^{1/2} \theta_\bw z_T} P_{\theta_\bw,t}
\{ T^{1/2}(S_t-c) \in I_T \}.
\label{21}
\end{eqnarray}
By similar computations, for any $y \in \mathbb{R}$,
\begin{eqnarray}
P_\bw \{ S_t \geq c+y \} 
= e^{-T \phi_\bw(c)} E_{\theta_\bw,t} \Big[ e^{T \theta_\bw(c-S_t)} {\bf 1}_{\{ S_t \geq
c+y \}} \Big] \leq e^{-T \phi_\bw(c)- T \theta_\bw y}.
\label{21a}
\end{eqnarray}
Under $P_{\theta_\bw, t}$, $T S_t^{(i)}$ is compound Poisson with Poisson mean
$\eta_i= \lambda_i \int_0^T e^{\theta_\bw g_\bw^{(i)}(u)} du$ and each summand is identically
distributed as $g_\bw^{(i)}(U_i)$, where $U_i$ is a random variable on $[0,T)$
with density $(\lambda_i/\eta_i) e^{\theta_\bw g_\bw^{(i)}(u)}$. 
We note that
$$
E_{\theta_\bw, t} \Big[ g_\bw^{(i)}(U_i) \Big] = (\lambda_i/\eta_i)
\int_0^T g_\bw^{(i)}(u) e^{\theta_\bw g_\bw^{(i)}(u)} \ du = \frac{d}{d \theta}
\int_0^T \lambda_i e^{\theta g_\bw^{(i)}(u)} \ du \Big|_{\theta=\theta_\bw} \Big/
\eta_i.
%\label{23}
$$
Since $\theta_\bw$ maximizes the right hand side of (\ref{5}), it follows that
\begin{equation}
E_{\theta_\bw, t} [S_t] = T^{-1} \sum_{i=1}^d \eta_i E_{\theta_\bw} [ g_\bw^{(i)}(U_i)]
=T^{-1} \frac{d}{d \theta} \sum_{i=1}^d \lambda_i \int_0^T (e^{\theta g_\bw^{(i)}(u)}-1) \ du
\Big|_{\theta=\theta_\bw} = c.
\label{24}
\end{equation}
Since a compound Poisson $Y =\sum_{j=1}^N Y_j$ has variance ${\rm Var}(Y) = (EN)(EY_1^2)$, 
it follows from (\ref{8}) that 
\begin{equation}
{\rm Var}_{\theta_\bw, t} (S_t) = T^{-2} \sum_{i=1}^d \eta_i \int_0^T [g_\bw^{(i)}(u)]^2
(\lambda_i/\eta_i) e^{\theta_\bw g_\bw^{(i)}(u)} du = T^{-1} v_\bw.
\label{24a}
\end{equation}
By (\ref{24}) and (\ref{24a}), $T^{1/2}(S_t-c)$ is asymptotically normal
with mean 0 and variance $v_\bw$. Hence by equation (5) of Stone (1965),
\begin{equation}
P_{\theta_\bw,t} \{ T^{1/2}(S_t-c) \in I_T \} = (2 \pi v_\bw)^{-1/2} \int_{I_T}
e^{-z^2/(2v_\bw)} dz + o_T(1)(\varepsilon_T+T^{-1/2}) \hspace{0.5cm} {\rm a.s.}
\label{24b}
\end{equation}
as $T \rightarrow \infty$, where $o_T(1)$ is a term not depending on 
$\varepsilon_T$ and $z_T$. 
Let $\varepsilon_T T^{1/2}$ tend to 0 slowly enough such that 
$o_T (1)/(\varepsilon_T T^{1/2}) \rightarrow
0$. Then by (\ref{24b}), 
\begin{equation}
P_{\theta_\bw,t} \{ T^{1/2}(S_t-c) \in I_T \} \sim (2 \pi v_\bw)^{-1/2} \int_{I_T}
e^{-z^2/(2v_\bw)} dz \hspace{0.5cm} \mbox{a.s.\ as $T \rightarrow \infty$},
\label{24c}
\end{equation}
if $I_T$ is uniformly bounded.
Then by (\ref{21}), (\ref{21a}) and (\ref{24c}),
\begin{eqnarray*}
P_\bw \{ S_t \geq c \} & = & \sum_{k=0}^{\lfloor \varepsilon_T^{-1} \rfloor}
P_{\bw} \{ k \varepsilon_T \leq T^{1/2}(S_t-c) < (k+1) \varepsilon_T 
\} + P_{\bw} \{ S_T \geq c+T^{-1/2} \} \cr
& \sim & (2 \pi v_\bw)^{-1/2} e^{-T \phi_\bw(c)} \int_0^\infty  e^{-T^{1/2}
\theta_\bw z-z^2/(2 v_\bw)} dz \hspace{0.5cm}\mbox{a.s.\ as $T\rightarrow\infty$,}
\end{eqnarray*}
and Lemma \ref{l2} holds. \hfill $\Box$

\begin{la}
\label{l3}
Assume (A1)-(A2). There exists $\varepsilon_T =o(T^{-1/2})$ such that 
for all uniformly bounded intervals $I_{1,T},I_{2,T}$ of
length $\varepsilon_T$, 
\begin{eqnarray}
& & P_{\theta_\bw,t} \Big\{ T^{1/2} \Big( S_t-c,\frac{d}{dx} S_x \Big|_{x=t} \Big)
\in I_{1,T} \times I_{2,T} \Big\} \cr
&\sim & (2 \pi)^{-1} (v_\bw \tau_\bw)^{-1/2} \Big(\int_{z_1 \in I_{1,T}} e^{-z_1^2/( 2
v_\bw) } \; dz_1 \Big) \Big( \int_{z_2 \in I_{2,T}} e^{-z_2^2/(2 \tau_\bw)} \ dz_2 \Big),
\label{22}
\end{eqnarray}
almost surely as $T\rightarrow \infty$.
\end{la}
{\sc Proof.} By stationarity, we may assume without loss of generality
that $t=0$. Under $P_{\theta_\bw}$, the vector $(T S_0^{(i)},
T \frac{d}{dx} S_x^{(i)} |_{x=0})'$ is bivariate compound Poisson with
Poisson mean $\eta_i = \lambda_i \int_0^T e^{\theta_\bw g_\bw^{(i)}(u)} \ du$
and with each summand identically distributed as
\begin{displaymath}
(g_\bw^{(i)}(U_i),-\frac{d}{du} g_\bw^{(i)}(u) |_{u=U_i})',
\end{displaymath}
where $U_i$ is a random variable on $[0,T)$ with density $(\lambda_i/\eta_i)
e^{\theta_\bw g_\bw^{(i)}(u)}$. By (A2), $\frac{d}{dx} S_x^{(i)} |_{x=0}$
exists almost surely. 

We shall now compute the means and covariances of 
\begin{displaymath}
(S_0,\frac{d}{dx} S_x |_{x=0})'
= \sum_{i=1}^d (S_0^{(i)}, \frac{d}{dx} S_x^{(i)} |_{x=0})'
\end{displaymath}
under
$P_{\theta_\bw}$. Since
\begin{eqnarray*}
E_{\theta_\bw} \Big[ - \frac{d}{du} g_\bw^{(i)}(u) \Big|_{u=U_i} \Big] & = &
-(\lambda_i/\eta_i) \int_0^T \Big[ \frac{d}{du} g_\bw^{(i)}(u) \Big]
e^{\theta_\bw g_\bw^{(i)}(u)} \ du \cr
& = & -(\lambda_i/\eta_i) \theta_\bw^{-1} (e^{\theta_\bw g_\bw^{(i)}(T)}-
e^{\theta_\bw g_\bw^{(i)}(0)}),
%\label{25}
\end{eqnarray*}
and $g_\bw^{(i)}$ is bounded, it follows that
\begin{equation}
E_{\theta_\bw} \Big[ \frac{d}{dx} S_x \Big|_{x=0} \Big] = T^{-1} \sum_{i=1}^d
\eta_i E_{\theta_\bw} \Big[ - \frac{d}{du} g_\bw^{(i)}(u) \Big|_{u=U^{(i)}} \Big]
=O(T^{-1}).
\label{26}
\end{equation}
The bivariate compound Poisson $(Y,Z)' = \sum_{j=1}^N (Y_j,Z_j)'$, where $(Y_1,Z_1)',
\cdots, (Y_N,Z_N)'$ are independent identically distributed summands conditioned on an independent Poisson
random variable $N$, has covariance
matrix
$$
{\rm Cov} \pmatrix{Y \cr Z} = E(N) \pmatrix{ E(Y^2) & E(YZ) \cr E(YZ) & E(Z^2)}.
%\label{27}
$$
It follows from the relation
\begin{eqnarray*}
E_{\theta_\bw} \Big[ g_\bw^{(i)}(U_i) \frac{d}{du} g_\bw^{(i)}(u) \Big|_{u=U_i} \Big] 
& = & (\lambda_i/\eta_i) \int_0^T \Big[ \frac{d}{du} g_\bw^{(i)}(u) \Big] g_\bw^{(i)}(u) 
e^{\theta_\bw g_\bw^{(i)}(u)} du \cr
& = & (\lambda_i/\eta_i) [\theta_\bw^{-1} g_\bw^{(i)}(u)- \theta_\bw^{-2}] 
e^{\theta g_\bw^{(i)}(u)}
\Big|_{u=0}^{u=T} = O(T^{-1}) \hspace{0.5cm} {\rm a.s.,}
\end{eqnarray*}
that
$$
{\rm Cov}_{\theta_\bw}(S_0, \frac{d}{dx} S_x|_{x=0}) = -T^{-2} \sum_{i=1}^d \eta_i
E_{\theta_\bw} [g_\bw^{(i)}(U_i) \frac{d}{du} g_\bw^{(i)}(u)|_{u=U_i}] = O(T^{-2}) \hspace{0.5cm} {\rm a.s.}
$$
and hence by (\ref{8}),
\begin{equation}
{\rm Cov}_{\theta_\bw} \pmatrix{ S_0 \cr \frac{d}{dx} S_x \big|_{x=0} } \sim T^{-1}
\pmatrix{ v_\bw & 0 \cr 0 & \tau_\bw } \hspace{0.5cm} \mbox{a.s.\
as $T \rightarrow \infty$}.
\label{28}
\end{equation}
It then follows from equation (5) of Stone (1965), (\ref{24}), (\ref{26}) and (\ref{28}) that
\begin{eqnarray*}
& & P_{\theta_\bw,t} \{ T^{1/2} \Big( S_t-c,\frac{d}{dx} S_x|_{x=t} \Big) \in
I_{1,T} \times I_{2,T} \} \cr
& = & 
(2 \pi)^{-1} (v_\bw \tau_\bw)^{-1/2} \Big(\int_{z_1 \in I_{1,T}} e^{-z_1^2/( 2
v_\bw) } \; dz_1 \Big) \Big( \int_{z_2 \in I_{2,T}} e^{-z_2^2/(2 \tau_\bw)} \ dz_2 \Big) \cr
& & \quad + o_T(1)(\varepsilon_T^2+T^{-1}), 
\end{eqnarray*}
where $o_T(1)$ does not depend on $I_{j,T}$, $j=1,2$. Then Lemma \ref{l3} follows by
selecting $\varepsilon_T$ such that $\varepsilon_T 
T^{1/2} \rightarrow 0$ and $o_T(1)/\varepsilon_T^2 T \rightarrow 0$.
\hfill $\Box$

\begin{la}
\label{l4}
Let $\kappa$, $T$, $K$ and $c$ be positive constants. Let
$$
%\label{29}
s(u) = c+ z_1 T^{-1/2} + u z_2 T^{-1/2} -\frac{u^2}{2} K.
$$
Then $\sup_{0 < u < \kappa T^{-1/2}} s(u) \geq \max \{ c,s(0),s(\kappa T^{-1/2}) \}$ 
if and only if 
$\kappa \geq z_2/K \geq 0$ and $z_1 \geq -z_2^2/(2K T^{1/2})$. 
\end{la}
{\sc Proof.} Since the quadratic $s$ has
a unique maximum at $u=z_2/(KT^{1/2})$, it follows that if $\kappa \geq z_2/K \geq 0$, then
$$
\sup_{0 < u < \kappa T^{-1/2}} s(u) = s \Big( \frac{z_2}{ KT^{1/2} } \Big) 
= c+ \frac{z_1}{ T^{1/2}} + \frac{ z_2^2 }{ 2KT}
$$
and Lemma \ref{l4} easily holds. \hfill $\Box$

\begin{la}
\label{l5}
Assume (A1)-(A2). Then for any $\kappa > 0$, $t \geq 0$ and $c >\mu$,
\begin{equation}
P_\bw \Big\{ \sup_{t < u \leq t+\kappa T^{-1/2}} S_u \geq \max (c, S_t, S_{t+\kappa
T^{-1/2}}) \Big\}
\sim \kappa T^{-1/2} \zeta_\bw e^{-T \phi_\bw(c)} \hspace{0.5cm} {\rm a.s.},
\label{31}
\end{equation}
where $\zeta_\bw = (2 \pi)^{-1} (\tau_\bw/v_\bw)^{1/2}$.
\end{la}
{\sc Proof.} Assume without loss of generality $t=0$. Let $H_i$ $(=H_{i,\bw})$
be the set of all $v$ such that a second derivative does not exists
at $g_\bw^{(i)}(v)$. Note that by (A1)-(A2), the number of elements
in $H_i$ is $O(T)$ a.s. for all $i$. Let $0 < u < \kappa T^{-1/2}$ and let
$y \in \by^{(i)}$ be such that
$y-h \not\in (0,u)$ for all $h \in H_i$. Then by (A2) and the mean value
theorem,
\begin{equation}
g_\bw^{(i)}(y-u)-g_\bw^{(i)}(y)+u \frac{d}{dv} g_\bw^{(i)}(v) \Big|_{v=y}
= \frac{u^2}{2} \frac{d^2}{dv^2} g_\bw^{(i)}(v) \Big|_{v=\xi}
\label{32}
\end{equation}
for some $y-u \leq \xi \leq y$. 

If $y \in \by^{(i)}$ is such that $y-h \in (0,u)$ for some $h \in H_i \setminus
\{ 0,T \}$, then
\begin{eqnarray}
\label{33}
&& g_\bw^{(i)}(y-u)-g_\bw^{(i)}(y)+u \frac{d}{dv} g_\bw^{(i)}(v) \Big|_{v=y}
\nonumber \\
& = & \int_{y-u}^y \Big( \frac{d}{dv} g_\bw^{(i)}(v)
\Big|_{v=y} - \frac{d}{d \xi} g_\bw^{(i)}(\xi) \Big) \ d \xi \cr
& = & (h+u-y) \Big( \frac{d}{dv} g_\bw^{(i)}(v)
\Big|_{v \downarrow h} - \frac{d}{d \xi} g_\bw^{(i)}(v) \Big|_{v \uparrow
h} \Big)+o(u^2).  
\end{eqnarray}
Finally, for completeness, we consider $y \in \by^{(i)}$ such that
either $y-T \in (0,u)$ or $y \in (0,u)$. Then we write formally
\begin{equation}
g_\bw^{(i)}(y-u)-g_\bw^{(i)}(y)+u \frac{d}{dv} g_\bw^{(i)}(v) \Big|_{v=y}
= g_\bw^{(i)}(y-u)-g_\bw^{(i)}(y)+u \frac{d}{dv} g_\bw^{(i)}(v) \Big|_{v=y}.
\label{34}
\end{equation}
By adding up (\ref{32})-(\ref{34}) over all $y \in \by^{(i)}$ for $i=1,\cdots,d$
and dividing by $T$, we obtain
\begin{equation}
S_u - S_0 - u \frac{d}{dv} S_v \Big|_{v=0} = - \frac{C_{\bw,u} u^2}{2},
\label{35}
\end{equation}
where $C_{\bw,u}$ is an expression derived from the right hand sides of (\ref{32})-(\ref{34}).
It shall be shown in Appendix B that
\begin{equation}
\lim_{T \rightarrow \infty} \sup_{0 < u < \kappa T^{-1/2}} u^2 T 
\Big| C_{\bw,u}-\theta_\bw \tau_\bw
\Big| \rightarrow 0 \quad {\rm a.s. \ under} \ P_{\theta_\bw}.
\label{36}
\end{equation}
Then by Lemma \ref{l3}, the change of variables
$$
%\label{37}
z_1 = T^{1/2}(S_0-c) \quad {\rm and} \ekp z_2 = T^{1/2} \frac{d}{dx} S_x \Big|_{x=0},
$$
substituting $K = \theta_\bw \tau_\bw$ into Lemma \ref{l4},
(\ref{35}), (\ref{5}) and (\ref{14}),
\begin{eqnarray*}
& & P_\bw \Big\{ \sup_{0 < u < \kappa T^{-1/2}} S_u 
\geq \max(c,S_0,S_{\kappa T^{-1/2}}) 
\Big\} 
\nonumber \\
&=& \quad E_{\theta_\bw} \Big[ \frac{dP_\bw}{dP_{\theta_\bw}}(\by) {\bf 1}_{\{ \sup_{0 < u <
\kappa T^{-1/2}} S_u \geq \max(c,S_0,S_{\kappa T^{-1/2}}) \}} \Big] \cr
&= & e^{-T \phi_\bw(c)} E_{\theta_\bw} \Big[ e^{T \theta_\bw (c-S_0)}
{\bf 1}_{\{ \sup_{0 < u <
\kappa T^{-1/2}} S_u \geq \max(c,S_0,S_{\kappa T^{-1/2}}) \}} \Big] \cr
&\sim & e^{-T \phi_\bw(c)} (2 \pi)^{-1} (v_\bw \tau_\bw )^{-1/2}
\nonumber \\
&&\hspace{0.5cm}\times
\int_0^{\kappa \theta_\bw \tau_\bw} \int^{\infty}_{-z_2^2/(2 \theta_\bw \tau_\bw T^{1/2})} 
e^{-T^{1/2} \theta_\bw z_1-z_1^2/2 v_\bw-z_2^2/
2 \tau_\bw} dz_1 \ dz_2 \cr
&\sim & e^{-T \phi_\bw(c)} (2 \pi)^{-1} (v_\bw \tau_\bw)^{-1/2} 
\nonumber \\
&&\hspace{0.5cm}\times
\int_0^{\kappa \theta_\bw
\tau_\bw} e^{-z_2^2/2 \tau_\bw} (-T^{-1/2} \theta_\bw^{-1} e^{-T^{1/2}
\theta_\bw z_1}) \Big|_{z_1=-z_2^2/(2 \theta_\bw \tau_\bw T^{1/2})}^{z_1=\infty} dz_2,
\end{eqnarray*}
and indeed 
Lemma \ref{l5} holds. \hfill $\Box$

The proof of the next lemma will also be shown in Appendix B.

\begin{la}
\label{l6}
Assume (A1)-(A2). Let 
$$
A_t = \Big\{ \sup_{t < u < t+\kappa T^{-1/2}} S_u \geq \max(c,S_t,S_{t+\kappa
T^{-1/2}}) \Big\}.
%\label{45}
$$
Then there exists $r_\kappa =o(\kappa)$ as $\kappa \rightarrow
\infty$ such that for all $t \geq 0$, with probability 1,
\begin{equation}
\sum_{1 \leq \ell \leq T^{3/2}/\kappa+1} P_\bw(A_t \cap A_{t+\ell \kappa T^{-1/2}})
\leq r_\kappa T^{-1/2} e^{-T \phi_\bw(c)}  
\label{46}
\end{equation}
for all large $T$. 
\end{la}

{\sc Proof of Proposition {\rm \ref{l1}}.} By stationarity, we may assume without
loss of generality $t=0$. By Lemmas \ref{l2}, \ref{l5},
\ref{l6} and the inequalities
\begin{eqnarray*}
& & \sum_{q=0}^{\lfloor \Delta/(\kappa T^{-1/2}) \rfloor-1} 
\Big[ P_\bw(A_{q \kappa T^{-1/2}})
- \sum_{\ell=1}^{\lfloor \Delta/(\kappa T^{-1/2}) \rfloor-1-q} P_\bw
(A_{q \kappa T^{-1/2}} \cap A_{(q+\ell) \kappa T^{-1/2}}) \Big] 
\\
&\leq & P_\bw \Big\{ \sup_{0 < u \leq \Delta} S_u \geq c \Big\} 
\\
&\leq &
\sum_{q=0}^{\lfloor \Delta/(\kappa T^{-1/2}) \rfloor} P_\bw(A_{q \kappa
T^{-1/2}})+\sum_{q=0}^{\lfloor \Delta/(\kappa T^{-1/2}) \rfloor+1} P_\bw \{ S_{q\kappa T^{-1/2}} \geq c \},
%\label{56}
\end{eqnarray*}
it follows that for any $0 < \varepsilon < 1$, there exists $\kappa$ arbitrarily large such that
\begin{eqnarray}
& & (\lfloor \Delta/(\kappa T^{-1/2}) \rfloor-1) \Big[ (1-\varepsilon) \kappa T^{-1/2}
\zeta_\bw e^{-T \phi_\bw(c)} - r_\kappa T^{-1/2} e^{-T \phi_\bw(c)} \Big] \cr
& \leq & P_\bw \Big\{ \sup_{0 < u \leq \Delta} S_u \geq c \Big\} \cr
& \leq & (\lfloor \Delta/(\kappa T^{-1/2}) \rfloor+1)(1+\varepsilon) \Big[ \kappa T^{-1/2} \zeta_\bw
e^{-T \phi_\bw(c)} + (2 \pi v_\bw)^{-1/2} \theta_\bw^{-1} T^{-1/2} e^{-T \phi_\bw(c)} \Big]
\label{46b}
\end{eqnarray}
holds for all large $T$ with probability 1. Select $\kappa$ large enough such that 
(\ref{46b}) and the inequalities
$$
r_\kappa \leq \varepsilon \kappa \zeta_\bw, \hspace{0.5cm} (2 \pi v_\bw)^{-1/2} \theta_\bw^{-1} \leq
\varepsilon \kappa \zeta_\bw,
$$
holds for all large $T$ with probability 1. Then by (\ref{46b}), the inequality 
\begin{equation}
\Big| \frac{P_\bw \{ \sup_{0 < u \leq \Delta} S_u \geq c \}}{\Delta \zeta_\bw e^{-T
\phi_\bw(c)}}-1 \Big| \leq 3 \varepsilon
\label{46c}
\end{equation}
holds for all large $T$ with probability 1 and (\ref{9}) follows from
(\ref{46c}) by selecting $\varepsilon > 0$ arbitrarily small. 
\hfill $\Box$

{\sc Proof of Theorem {\rm \ref{t1}}.} Let $z \in {\mathbb{R}}$ and let $\xi$ ($=\xi_\bw$) 
be such that $\xi/T \rightarrow \infty$ and
$k$ $(=k_\bw):= z e^{T \phi_\bw(c)}/\zeta_\bw \xi$ is a positive integer
tending to infinity almost surely. Define
$B_j = \{ \sup_{(j-1) \xi \leq t < j \xi-T} S_t \geq c \}$ and
$C_j = \{ \sup_{j \xi-T \leq t \leq j \xi} S_t \geq c \}$. Then
\begin{equation}
P_\bw \Big( \bigcup_{j=1}^k B_j \Big) \leq P_\bw \{ \zeta_\bw e^{-T \phi_\bw(c)}
V_c \leq z \} \leq P_\bw \Big( \bigcup_{j=1}^k B_j \Big) + \sum_{j=1}^k
P_\bw(C_j).
\label{57}
\end{equation}
Conditioned on $\bw$ known, the event $B_j$ depends only on the spike
train times of $\by^{(i)}$ lying inside $[(j-1)\xi,j \xi)$. Since these
intervals are disjoint for different $j$, it follows that $B_1,\cdots,B_k$ are 
independent conditioned on $\bw$. 
By Lemmas \ref{l5} and \ref{l6}, it follows that with
probability 1,
\begin{equation}
P_\bw(B_j) \sim (\xi-T) \zeta_\bw e^{-T \phi_\bw(c)} \sim z/k, \quad
P_\bw(C_j) \sim  T \zeta_\bw e^{-T \phi_\bw(c)} \quad {\rm for \ all} \ 
1 \leq j \leq k.
\label{58}
\end{equation}
Since $k \rightarrow \infty$ a.s. as $T \rightarrow \infty$, with 
probability 1,
\begin{equation}
P_\bw \Big( \bigcup_{j=1}^k B_j \Big) = 1 - \prod_{j=1}^k P_\bw(B_j^c)
= 1-(1-z/k)^k+o(1) \rightarrow 1-e^{-z}.
\label{59}
\end{equation}
Moreover because $\xi/T \rightarrow \infty$, it follows from (\ref{58})
that with probability 1, 
\begin{equation}
\sum_{j=1}^k P_\bw(C_j) \sim T k \zeta_\bw e^{-T \phi_\bw(c)} = o(1).
\label{60}
\end{equation}
Theorem \ref{t1}(a) follows from (\ref{57}), (\ref{59}) and (\ref{60}).

To show Theorem \ref{t1}(b), we use the Taylor expansion
\begin{eqnarray}
\phi_\bw \Big( c_\bw + \frac{z+\log \zeta_\bw}{\theta_\bw T} \Big)
& = & \phi_\bw(c_\bw) + \theta_\bw \Big( \frac{z+\log \zeta_\bw}{\theta_\bw
T} \Big) + O(T^{-2}) \cr
& = & T^{-1} [z+\log(a \zeta_\bw)]+O(T^{-2}).
\label{61}
\end{eqnarray}
By the computations in (\ref{57})-(\ref{60}), it follows that 
with probability 1,
\begin{eqnarray}
P_\bw \Big\{ M_a \geq c_\bw + \frac{z+\log \zeta_\bw}{\theta_\bw T}
\Big\} & = & P_\bw \{ V_{c_\bw +[z+\log \zeta_\bw]/(\theta_\bw T)} \leq a \}
\cr
& = & 1 - \exp \Big[-a \zeta_\bw e^{-T \phi_\bw(c_\bw+
[z+\log \zeta_\bw]/(\theta_\bw T) )}+o(1) \Big],
\label{62}
\end{eqnarray}
and Theorem \ref{t1}(b) follows by substituting (\ref{61}) into (\ref{62}).

It remains to show (c). Let $\wtd U_a = \sum_{j=1}^{\lfloor a/\xi \rfloor}
{\bf 1}_{B_j}$. Then by (\ref{58}), $E_\bw[\wtd U_a]-\eta_\bw \rightarrow 0$ a.s. and
(\ref{11}) holds with $U_a$ replaced by $\wtd U_a$, since $B_j$, $1 \leq
j \leq \lfloor a/\xi \rfloor$ are independent events and the limit of
sum of binomial random variables is Poisson. By (\ref{58}) and (\ref{60}),
\begin{equation}
\sum_{j=1}^{\lfloor a/\xi \rfloor} P_\bw(C_j) = o(1) \quad {\rm and} \ 
\sum_{j=1}^{\lfloor a/\xi \rfloor-1} P_\bw(B_j \cap B_{j+1})
= o(1) \hspace{0.5cm} \mbox{a.s.\ as $T \rightarrow \infty$}.
\label{63}
\end{equation}
Moreover, by Lemma \ref{l6}, with probability 1,
\begin{equation}
\sum_{q=0}^{\lfloor a/(\kappa T^{-1/2}) \rfloor} \Big[ \sum_{\ell=\lfloor
(1-\alpha) T^{3/2}/\kappa \rfloor}^{\lfloor T^{3/2}/\kappa+1
\rfloor} P_\bw (A_{q \kappa T^{-1/2}}
\cap A_{(q+\ell)\kappa T^{-1/2}}) \Big] \leq \frac{ \eta_\bw r_\kappa }{ \kappa \zeta_\bw }
\label{64}
\end{equation}
for all large $T$, where $\alpha T$ is the maximal permitted overlap
between two matches. By (\ref{63}) and
(\ref{64}) with $\kappa$ arbitrarily large, 
we can conclude that $\wtd U_a -U_a \rightarrow 0$ in
probability and hence (\ref{11}) holds. \hfill $\Box$

\section{Template matching with kernels containing discontinuities}

In this section, we obtain analogues of Proposition \ref{l1} and Theorem
\ref{t1} when the score function $f$ contains discontinuities. 
A typical example is the box kernel
\begin{equation}
f(x) = \cases{1 & if $x < \varepsilon$, \cr
-\beta & if $x \geq \varepsilon$, } 
\label{65}
\end{equation}
where $\beta,\varepsilon$ are positive real numbers. Instead of (A2), we
have the following regularity condition on $f$.

\smallskip
\noi (A2)$'$ Let $f$ be a discontinuous function and let 
there be a finite set $H$ such that the first derivative
of $f$ exists and is uniformly continuous over any interval within
${\mathbb{R}}^+ \setminus H$. Moreover, (\ref{4}) holds.

\smallskip Under (A2)$'$, the values of $f$ may be concentrated on
$0, \pm q, \pm 2q, \cdots$ for some $q > 0$.

{\sc Definition.} Let
$L(f) = \{ f(x): x \geq 0 \}$ be the range of $f$. We say that $f$ is arithmetic if 
\begin{equation}
L(f) \subseteq q {\mathbb{Z}} \hspace{0.5cm} \mbox{for 
some $q > 0$}.
\label{68}
\end{equation}
Moreover, if $q$ is the largest number satisfying (\ref{68}), then
we say that $f$ is arithmetic with span $q$. If (\ref{68}) is not
satisfied for all $q > 0$, we say that $f$ is nonarithmetic. 

For example,
if $\beta$ in (\ref{65}) is irrational, then $f$ is nonarithmetic while
if $\beta=s/r$ for co-primes $r$ and $s$, then $f$ is arithmetic with span
$q=r^{-1}$. We write for $i=1,\cdots, d$,
\begin{eqnarray*}
g_\bw^{(i)}(u+) &=& \lim_{v \downarrow u} g_\bw^{(i)}(v), 
\nonumber \\
g_\bw^{(i)}(u-) &=& \lim_{v \uparrow u} g_\bw^{(i)}(v), 
\nonumber \\
\delta_i(u) &=& g_\bw^{(i)}(u-)-g_\bw^{(i)}(u+), \cr
D_i &= &\{ u \in (0,T): \delta_i(u) \neq 0 \}, 
%\quad D_i^+ = \{ u \in (0,T): \delta_i(u) > 0 \}.
%\label{67}
\end{eqnarray*}
where $g_\bw^{(i)}$ is as in (\ref{2}).

Let $\phi_\bw$, $\theta_\bw$, $v_\bw$ and $\mu$ be as in Section 5.1.
If $D_i\neq \emptyset$
for some $i$, we can define $h_\bw^*$ to be the probability
mass function 
taking values in $\{ \delta_i(u) \}_{u \in D_i, 1 \leq i \leq d}$
with
\begin{equation}
h_\bw^*(x) = \sum_{i=1}^d \lambda_i \sum_{u \in D_i}
e^{\theta_\bw g_\bw^{(i)}(u-)} {\bf 1}_{\{ \delta_i(u)=x \}} \Big/
\sum_{i=1}^d \lambda_i \sum_{u \in D_i}
e^{\theta_\bw g_\bw^{(i)}(u-)}.
\label{69}
\end{equation}
Let $E_*$ denotes expectation when $X_1,X_2,\cdots$ are 
independent identically distributed random variables with probability mass function
$h_\bw^*$. Define
\begin{equation}
\omega_b = \inf \{ n: X_1+\cdots+X_n \geq b \} \quad {\rm and} \ 
R_b = X_1 + \cdots + X_{\omega_b}.
\label{70}
\end{equation}
Then the overshoot constant
\begin{equation}
\nu_\bw := \lim_{b \rightarrow \infty} E_* e^{-\theta_\bw(R_b-b)},
\label{71}
\end{equation}
where $b$ is taken to be a multiple of $\chi$ if $h_\bw^*$ is
arithmetic with span $\chi$. Note that the statement ``$h_\bw^*$ is
arithmetic with span $\chi$'' implies that $\{ \delta_i(u) \}_{u \in D_i, 1 \leq i \leq d}
\subset \chi \mathbb{Z}$. The constants $\nu_\bw$ have been
well-studied in sequential analysis, see for example
Siegmund (1985) for the existence of the limits in (\ref{71}).
Let us define the asymptotic constant $\zeta_\bw'$ (not depending on {\bf y})
by
\begin{equation}
\zeta_\bw' = (2 \pi T v_\bw)^{-1/2} \nu_\bw K_\bw \sum_{i=1}^d
\lambda_i \sum_{u \in D_i} \delta_i(u) e^{\theta_\bw g_\bw^{(i)}(u-)},
\label{72}
\end{equation}
where
$$
K_\bw = \cases{ 1 & if $h_\bw^*$ is nonarithmetic, \cr
\frac{1}{\theta_\bw \chi} (1-e^{-\theta_\bw \chi})
& if $h_\bw^*$ is arithmetic with span $\chi$, $f$ is nonarithmetic, \cr
\frac{q}{\chi} \left( \frac{1-e^{-\theta_\bw \chi}}{1-e^{-\theta_\bw
q}} \right) & if $h_\bw^*$ is arithmetic with span $\chi$, $f$ is
arithmetic with span $q$.}
%\label{73}
$$
Since we can express each $\delta_i(u)$, $u \in D_i$ in the form $g_1-g_2$ for
$g_1,g_2 \in L(f)$, it follows that if $f$ is arithmetic, then $h_\bw^*$ is arithmetic
and $\chi/q$ is a positive integer.
Analogous to Lemma \ref{l1} and Theorem \ref{t1}, we have the
following asymptotic results.

\begin{pn}
\label{l7}
Assume (A1), (A2)$'$ and let $\Delta > 0$, $t \geq 0$.
If $f$ is nonarithmetic and $c > \mu$, then
\begin{equation}
P_\bw \Big\{ \sup_{t < u \leq t+\Delta} S_u \geq c \Big\}
\sim \Delta \zeta_\bw' e^{-T \phi_\bw(c)} \hspace{0.5cm} \mbox{a.s.\ as $T
\rightarrow \infty$}.
\label{74}
\end{equation}
If $f$ is arithmetic with span $q$, then (\ref{74}) also holds under
the convention that
\begin{equation}
Tc\ (=Tc_T) \in q {\mathbb{Z}} \hspace{0.5cm} \mbox{with $c \rightarrow c'$
as $T \rightarrow \infty$ for some $c' > \mu$}.
\label{75}
\end{equation}
\end{pn}

\begin{tm}
\label{t2}
Assume (A1) and(A2)$'$ and let $f$ be nonarithmetic.

{\rm (a)} Let $c > \mu$. Then the distribution (conditional on $\bw$) of 
$\zeta_\bw' e^{-T \phi_\bw(c)} V_c$
converges to the exponential distribution with mean 1 almost surely as $T\rightarrow \infty$.

{\rm (b)} Let $a \rightarrow \infty$ as $T \rightarrow \infty$ such that
$(\log a)/T$ converges to a positive constant. 
Let $c_\bw > \mu_\bw$ satisfy $\phi_\bw(c_\bw) = (\log a)/T$. Then
for any $z \in \mathbb{R}$, 
$$
P_\bw \{ \theta_\bw T (M_a -c_\bw) - \log \zeta_\bw' \geq z \} 
\rightarrow 1-\exp(-e^{-z}) \hspace{0.5cm} \mbox{a.s.\ as $T \rightarrow \infty$}.
$$

{\rm (c)} Let $a \rightarrow \infty$ as $T \rightarrow \infty$ such that
$(\log a)/T$ converges to a positive constant. Let $c$ ($=c_T$)
be such that $\eta_\bw := a \zeta_\bw' e^{-T \phi_\bw(c)}$ converges to
a constant $\eta > 0$ almost surely. Then
\begin{equation}
P_\bw \{ U_a = k \} - e^{-\eta_\bw} \frac{\eta_\bw^k}{k!}
\rightarrow 0 \hspace{0.5cm} \mbox{a.s. $\forall \ k=0,1,\cdots$}.
\label{10}
\end{equation}

If $f$ is arithmetic with span $q$, then (a) and (c) also hold
under the convention (\ref{75}).
\end{tm}

{\sc Example 2.} We shall conduct here a simulation study similar to Example 1.
The generation of $\bw$ and $\by$ are similar to Example 1 but the
box kernel (\ref{65}) is used instead of the Hamming window function
(\ref{16}) when computing $g_\bw^{(i)}$. We chose parameters $\varepsilon=4$ms and
$\beta=0.3$ in (\ref{65}). Hence $f$ is arithmetic with span $q=0.1$ and $h_\bw^*$
is arithmetic with span $\chi=1.3$. In fact, $h_\bw^*$ is positive only on
the values $-1.3$ and $1.3$ and hence $\nu_\bw=1$. In the template $\bw$ that was
generated, there were a total of 2 $\times$ 59 elements in $\bigcup_i D_i$ with half of 
all $u \in D_i$ satisfying $\delta_i(u)=1.3$ and the other half
satisfying $\delta_i(u)=-1.3$ for each $i$. Hence 
$$
h_\bw^* (-1.3) = e^{-0.3 \theta_\bw}/(e^{\theta_\bw}+
e^{-0.3 \theta_\bw}) \quad {\rm and} \quad h_\bw^* (1.3) = e^{\theta_\bw}/
(e^{\theta_\bw}+e^{-0.3 \theta_\bw}).
$$
These information are then used in the computation of the constant $\zeta_\bw'$ 
in the analytical approximation
\begin{equation}
P_\bw \{ M_a \geq c \} = P_\bw \{ V_c \leq a \} \doteq 1 - \exp(-a
\zeta_\bw' e^{-T \phi_\bw(c)}),
\label{76}
\end{equation} 
an analogue of (\ref{19}) that follows from Theorem \ref{t2}(a).

\begin{table}
\begin{center}
{\sc Table 2.} Estimates of $P_\bw \{ M_a \geq c \} \pm$ standard error with
$a+T=20s.$

\begin{tabular}{c|c|c|c}
\hline
$c$ & Direct MC & Imp. Sampling & Anal. Approx. (\ref{76}) \cr
\hline
0.065 & 0.029$\pm$0.004 & 0.0300$\pm$0.0016 & 0.0289 \cr 
0.066 & 0.019$\pm$0.003 & 0.0218$\pm$0.0012 & 0.0207 \cr
0.067 & 0.012$\pm$0.002 & 0.0140$\pm$0.0008 & 0.0144 \cr
0.068 & 0.008$\pm$0.002 & 0.0103$\pm$0.0006 & 0.0101 \cr
0.069 & 0.005$\pm$0.002 & 0.0067$\pm$0.0004 & 0.0070 \cr
0.070 & 0.003$\pm$0.001 & 0.0051$\pm$0.0003 & 0.0047 \cr
\hline
\end{tabular}
\end{center}
\end{table}

In Table 2, we compare the analytical approximation (\ref{76}) against both
direct Monte Carlo and importance sampling with 2000 simulation runs
being used to obtain each entry. The variance reduction when importance
sampling is used is similar to
that seen in Example 1, and the technique is indeed an effective time saving
device for computing $p$-values especially when the $p$-values are small. 
The analytic approximations are also accurate and
agree with the simulation results that were obtained. 

In addition to the above simulation study, we also conducted a 
similar exercise to check the accuracy of the Poisson approximation of $U_a$
in (\ref{10}), this time with $a+T=200s$ and threshold level $c=0.0614$.
The maximal proportion of overlap between two matches is chosen to be $\alpha=0.8$. The
analytical approximations are compared against 2000 direct Monte Carlo
simulation runs and the results are recorded in
Table 3. Again we see that the analytical approximations are quite accurate
and this indicates the usefulness of using the asymptotic results in 
Theorem \ref{t2} to estimate $p$-values.

\begin{table}
\begin{center}
{\sc Table 3.} Estimates of $P_\bw \{ U_a = k \}$ and $\eta_\bw = E_\bw (U_a)$. 
Standard errors in parentheses.

\begin{tabular}{c|c|c|c|c|c|c|c||c}
\hline
$k$ & 0 & 1 & 2 & 3 & 4 & 5 & $\geq$ 6 & $\eta_\bw$ \cr
\hline
Approx. (\ref{10}) & 0.336 & 0.366 & 0.200 & 0.073 & 0.020 & 0.004 & 0.001 & 1.09 \cr 
%(Theorem \ref{t2}) & & & & & & & \cr
\hline
Direct Monte & 0.328 & 0.363 & 0.195 & 0.084 & 0.024 & 0.005 & 0.001 & 1.13 \cr
Carlo & (0.011) & (0.011) & (0.009) & (0.006) & (0.003) & (0.002) & (0.001) & (0.02) \cr
\hline
\end{tabular}
\end{center}
\end{table}

\medskip
We shall now prove Proposition \ref{l7} and Theorem \ref{t2} via the following
preliminary lemmas. Let 
\begin{equation}
h_\bw(x) = \sum_{i=1}^d \lambda_i \sum_{u \in D_i}
e^{\theta_\bw g_\bw^{(i)}(u+)} {\bf 1}_{\{ \delta_i(u)=x \}} \Big/
\sum_{i=1}^d \lambda_i \sum_{u \in D_i}
e^{\theta_\bw g_\bw^{(i)}(u+)}.
\label{79}
\end{equation}
Then $h_\bw$ and $h_\bw^*$ [see (\ref{69})] are conjugate probability
mass functions in the following
sense.

\begin{la}
\label{l8}
Let $D_i\neq \emptyset$ for some $i$. Then there exists $\gamma_\bw = 
1+O(T^{-1})$ a.s. such that
$$
h_\bw^*(x) = \gamma_\bw e^{\theta_\bw x} h_\bw(x) \hspace{0.5cm} \forall x.
$$
\end{la}
{\sc Proof.} Let $u \in D_i$ with $w_j^{(i)} < u < w_{j+1}^{(i)}$
for adjacent spikes $w_j^{(i)}, w_{j+1}^{(i)} \in \bw^{(i)}$. Then
by the symmetry of $g_\bw^{(i)}$ in the interval $(w_j^{(i)},w_{j+1}^{(i)})$
about its mid-point $(w_j^{(i)}+w_{j+1}^{(i)})/2$, it follows that
$v:=w_{j+1}^{(i)}-(y-w_j^{(i)}) \in D_i$ and that $g_\bw^{(i)}(v-) = 
g_\bw^{(i)}(u+)$. Hence $\gamma_\bw$, which we define here to be the ratio
of the denominators on the right-hand sides of (\ref{69}) and (\ref{79}),
is $1+O(T^{-1})$ a.s. with the $O(T^{-1})$ coming from $u \in D_i$ occuring
before the first spike or after the last spike in $\bw^{(i)}$. Since
$$
e^{\theta_\bw g_\bw^{(i)}(u-)} {\bf 1}_{\{ \delta_i(u)=x \}} = 
e^{\theta_\bw[x+g_\bw^{(i)}(u+)]} {\bf 1}_{\{ \delta_i(u)=x \}}, 
%\label{81}
$$
Lemma \ref{l8} holds. \hfill $\Box$

\begin{la}
\label{l9}
Assume (A1), (A2)$'$ and let $\kappa > 0$. Then for all $\varepsilon > 0$,
there exists $\kappa$ large enough such that for any 
for any $t \geq 0$, the inequality
$$
\Big| \frac{ P_\bw \{ S_t < c, \sup_{t < u \leq t+\kappa T^{-1}} S_u \geq c
\} }{ \kappa T^{-1} \zeta_\bw' e^{-T \phi_\bw(c)} } - 1 \Big| \leq \varepsilon, 
%\label{82}
$$
holds for all large $T$ with probability 1, where the constant $\zeta_\bw'$ is
defined in (\ref{72}) and $c > \mu$ if $f$ is nonarithmetic; $c$
satisfies (\ref{75}) if $f$ is arithmetic with span $q$.
\end{la}
{\sc Proof.} Assume without loss of generality $t=0$ and let
$G_i = \bigcup_{v \in D_i} (v,v+\kappa T^{-1}]$. We can write 
$TS_0 = TS_0' + J_0$, where
$$
S_0' = T^{-1} \sum_{i=1}^d \sum_{y \in \by^{(i)}, y \not\in G_i}
g_\bw^{(i)}(y) \ {\rm and} \ J_0 = \sum_{i=1}^d \sum_{y \in \by^{(i)} \cap G_i} g_\bw^{(i)}(y).
%\label{83}
$$
The random variables $S_0'$ and $J_0$ are independent because they
are functions of the Poisson processes $\by^{(i)}$ over disjoint
subsets of the real line. Let us first consider $f$ arithmetic with span
$q$. Then $TS_0'$ and $J_0$ are both integral multiples of $q$. Since $f$
is constant between jumps, we can 
express $TS_u=TS_0' + J_u$ [see (\ref{1})] where
\begin{equation}
J_u = \sum_{i=1}^d \sum_{y \in \by^{(i)} \cap G_i} g_\bw^{(i)}(y-u) \qquad \forall 
u \in (0,\kappa T^{-1}).
\label{84}
\end{equation}
Hence both $S_0 < c$ and $\sup_{0 < u \leq \kappa T^{-1}} S_u \geq c$ occurs
if and only if
\begin{equation}
\sup_{0 < u \leq \kappa T^{-1}} (J_u-J_0) \geq \ell q \ {\rm and} \ 
S_0' =c-T^{-1} kq \ {\rm for \ some \ integer} \ \ell \geq 1 \ {\rm and} \ k=J_0/q+\ell.
\label{85}
\end{equation}
Since $E_{\theta_\bw}[S_0']=c+O(T^{-1})$ a.s., by the local limit theorem for lattice
random variables [see, for example, Theorem 15.5.3 of Feller (1971)], 
\begin{equation}
P_{\theta_\bw} \{ S_0' = c-T^{-1}qk \} \sim \frac{q}{(2 \pi T
v_\bw)^{1/2}} \ {\rm a.s.}
\label{86}
\end{equation}
for any integer $k$. Since $S_0'$ and $( J_u )_{0 < u \leq
\kappa T^{-1}}$ are independent, it follows from (\ref{85}), (\ref{86}),
(\ref{5}) and the change of measure (\ref{14}) that
\begin{eqnarray}
& & P_\bw \Big\{ S_0 < c, \sup_{0 < u \leq \kappa T^{-1}} S_u \geq c
\Big\} 
\nonumber \\
&=& E_{\theta_\bw} \Big[ \frac{dP_\bw}{dP_{\theta_\bw}}(\by)
{\bf 1}_{\{ S_0 < c, \sup_{0 < u \leq \kappa T^{-1}} S_u \geq c \}}
\Big] 
\nonumber \\
&\sim & \frac{q}{(2 \pi T v_\bw)^{1/2}}
e^{-T \phi_\bw(c)} \sum_{\ell=1}^\infty e^{\theta_\bw \ell q}
P_{\theta_\bw} \Big\{ \sup_{0 < u \leq \kappa T^{-1}}
(J_u-J_0) \geq \ell q \Big\}.
\label{87}
\end{eqnarray}
Since $g_\bw^{(i)}$ is constant between jumps, the graph of
$(J_u-J_0)$ against $u$ is also piecewise constant with jumps at all
$u$ for which $y-u \in D_i$ for some  
$y \in \by^{(i)}$, $1 \leq i \leq d$ [see (\ref{84})]. Let $N_*$ be the total number of spikes in 
$\bigcup_{1 \leq i \leq d} (\by^{(i)} \cap G_i)$. Then there are $N_*$
such jumps and
$$
\sup_{0 < u \leq \kappa T^{-1}} (J_u-J_0) = \sup_{1 \leq j \leq N_*}
(X_1+\cdots+X_j)
$$
where $X_j$ is the $j$th jump and has probability mass
function $h_\bw$.  Moreover, $X_1,X_2,\cdots$ are independent conditioned on 
$N_*$, a Poisson random variable independent of the $X_i$'s 
with mean 
\begin{equation}
EN_* = \kappa T^{-1} \sum_{i=1}^d \lambda_i \sum_{u \in D_i} e^{\theta_\bw
g_\bw^{(i)}(u+)}. 
\label{88}
\end{equation}
If $r \in \{ 0,\cdots,\chi/q-1 \}$
and $s \in {\mathbb{Z}}^+$, then $R_{s \chi-rq} = R_{s \chi}$ [see (\ref{70})]. 
Let $E_*$ and $P_*$ denote the expectation and probability
measure respectively when $X_1,X_2,\cdots$ are
independent identically distributed with probability mass function $h_\bw^*$. 
Then it follows from
a change of measure to $P_*$ and Lemma \ref{l8} that 
\begin{eqnarray}
& & \sum_{\ell=1}^\infty e^{\theta_\bw \ell q}
P_{\theta_\bw} \Big\{ \sup_{0 < u \leq \kappa T^{-1}}
(J_u-J_0) \geq \ell q \Big\} \cr
&= &  \gamma_\bw^{-1} \sum_{\ell=1}^\infty E_* \Big[
e^{-\theta_\bw(R_{\ell q}-\ell q)} {\bf 1}_{\{ \sup_{1 \leq j \leq N_*}
(X_1+\cdots+X_j) \geq \ell q \}} \Big] \cr
&= &  \gamma_\bw^{-1} E_* \Big[ \sum_{r=0}^{\chi/q-1} \sum_{s=1}^\infty
e^{-\theta_\bw[R_{s \chi}-(s \chi-rq)]} {\bf 1}_{\{ \sup_{1 \leq j
\leq N_*} (X_1+\cdots+X_j) \geq s \chi \}} \Big] \cr
&\sim &\chi^{-1} \Big( \sum_{r=0}^{\chi/q-1} e^{-\theta_\bw rq} \Big) \nu_\bw
E_* \Big[ \sup_{1 \leq j \leq N_*} (X_1+\cdots+X_j) \Big].
\label{89}
\end{eqnarray}
Since $X_i$ has positive mean under $P_*$ for all large $T$ and the almost
sure limit of $E_* N$ [see (\ref{88})] is proportional to $\kappa$, it follows that
there exists $\kappa$ large enough such that
\begin{equation}
\Big|  
\frac{ E_* [ \sup_{1 \leq j \leq N_*} (X_1+\cdots+X_j) ] }{ (EN_*) (E_* X_1)} -1 \Big|
< \frac{\varepsilon }{2}
\label{90}
\end{equation}
for all large $T$. Since $(q/\chi) \sum_{r=0}^{\chi/q-1} e^{-\theta_\bw rq} = K_\bw$,
Lemma \ref{l9} then follows from (\ref{69}) and
(\ref{87}) to (\ref{90}). When $f$ is nonarithmetic, 
the local limit result (\ref{22}) with $I_{2,T} = \mathbb{R}$ and $t=0$ 
is used in place of (\ref{86}). \hfill $\Box$. 

The next lemma, needed for the proofs of both Proposition \ref{l7} and Theorem \ref{t2}, 
will be proved in Appendix B.

\begin{la}
\label{l10}
Assume (A1) and (A2)$'$. Let 
$$
A_t = \Big\{ S_t < c, \ \sup_{t < u \leq t+\kappa T^{-1}} S_u \geq c \Big\}.
%\label{91}
$$
Then there exists $r_\kappa =o(\kappa)$ as $\kappa \rightarrow
\infty$ such that for all $t \geq 0$, with probability 1,
$$
\sum_{\ell=1}^{\lfloor T^2/\kappa+1 \rfloor} P_\bw(A_t \cap A_{t+\ell \kappa T^{-1}})
\leq r_\kappa T^{-1/2} e^{-T \phi_\bw(c)}  
%\label{92}
$$
for all large $T$.
\end{la}

{\sc Proof of Proposition {\rm \ref{l7}}.} By stationarity, we may assume without
loss of generality $t=0$. Then (\ref{74}) follows from Lemmas \ref{l2}, \ref{l9}, \ref{l10}
and the inequalities
\begin{eqnarray*}
& & \sum_{q=0}^{\lfloor \Delta/(\kappa T^{-1}) \rfloor-1} 
\Big[ P_\bw(A_{q \kappa T^{-1}})
- \sum_{\ell=1}^{\lfloor \Delta/(\kappa T^{-1}) \rfloor-1-q} P_\bw
(A_{q \kappa T^{-1}} \cap A_{(q+\ell) \kappa T^{-1}}) \Big] \cr
& \leq & P_\bw \Big\{ \sup_{0 < u \leq \Delta} S_u \geq c \Big\} 
\leq 
\sum_{q=0}^{\lfloor \Delta/( \kappa T^{-1}) \rfloor} P_\bw(A_{q \kappa
T^{-1}})+P_\bw \{ S_0 \geq c \},
\end{eqnarray*}
with $\kappa$ arbitrarily large; see for example, the proof of Proposition
\ref{l1}. \hfill $\Box$

{\sc Proof of Theorem {\rm \ref{t2}}.} The proof of Theorem \ref{t2}
proceeds as in the proof of Theorem \ref{t1}. The only
modification needed is the replacement of $\zeta_\bw$ by $\zeta_\bw'$. For the proof
of Theorem \ref{t2}(c), we will also need to replace (\ref{64}) by
$$
\sum_{q=0}^{\lfloor a/(\kappa T^{-1}) \rfloor} \Big[ \sum_{\ell=\lfloor
(1-\alpha) T^2/\kappa \rfloor}^{\lfloor T^2/\kappa+1 
\rfloor} P_\bw (A_{q \kappa T^{-1}}
\cap A_{(q+\ell)\kappa T^{-1}}) \Big] \leq \frac{ \eta_\bw r_\kappa T^{1/2}}{ \kappa \zeta_\bw'}. 
$$
\hfill $\Box$

\section{Acknowledgments}

Wei-Liem Loh would like to thank Professor Yannis Yatracos for the many discussions on metric entropy and minimum
distance estimation and to Professor Zhiyi Chi for introducing to him the field of neuroscience when he visited the
University of Chicago in Spring 2003.

\section{Appendix A}

\begin{la} \label{la:a.3}
Let $\Theta_{\tilde{\kappa}, q, n}^2$ and $\rho_{\Theta_{\tilde{\kappa}, q, n}^2}$ be as in Section 2. 
Then for each $(s_1, r_1)\in \Theta_{\tilde{\kappa}, q, n}^2$,
\begin{eqnarray*}
&& \Big\{ \sum_{j=0}^\infty \int^*_{0<w_1<\cdots < w_j <T} \sup_{ \rho_{\Theta_{\tilde{\kappa}, q, n}^2} ((s_1, r_1), (s_2, r_2))\leq \varepsilon, 
(s_2,r_2)\in \Theta_{\tilde{\kappa}, q, n}^2 } [ p_{s_1,r_1}^{1/2} (\{w_1,\cdots, w_j\})
\nonumber \\
&& \hspace{0.5cm} - p_{s_2, r_2}^{1/2} (\{ w_1,\cdots, w_j\}) ]^2 dw_1 \cdots dw_j \Big\}^{1/2}  \leq  
\varepsilon C_{\tilde{\kappa}},
\end{eqnarray*}
where $C_{\tilde{\kappa}} \geq 1/2$ is a constant depending only on $\tilde{\kappa}$. 
Here $\int^*$ denotes the upper integral [see for example Dudley (1999), page 94].
Consequently,
\begin{equation}
H^B (\varepsilon, {\cal F}_{\tilde{\kappa}, q, n}, \rho_{{\cal F}_{\tilde{\kappa}, q, n}} ) 
\leq H( \frac{ \varepsilon}{ 2 C_{\tilde{\kappa}} }, \Theta_{\tilde{\kappa}, q, n}^2, \rho_{\Theta_{\tilde{\kappa}, q, n}^2})
\leq \frac{ 2^{(q+2)/q} C_{\tilde{\kappa}}^{1/q} C_{\tilde{\kappa}, q } }{\varepsilon^{1/q} }.
\label{eq:4.9}
\end{equation}
\end{la}
{\sc Proof.} 
Let $\bar{\kappa} = \kappa_0\vee 1$.
We observe from (\ref{eq:4.1}) that
\begin{eqnarray*}
&& | p_{s_1,r_1}^{1/2} (\{w_1,\cdots, w_j\})
- p_{s_2, r_2}^{1/2} (\{ w_1,\cdots, w_j\}) |
\nonumber \\
&=& | e^{-\int_0^T s_1(t) r_1(t-w_{\zeta(t)}) dt/2} \prod_{i=1}^j s_1^{1/2} (w_i) r_1^{1/2} (w_i- w_{i-1})
\nonumber \\
&&\hspace{0.5cm}
- e^{-\int_0^T s_2 (t) r_2 (t-w_{\zeta(t)}) dt/2} \prod_{i=1}^j s_2^{1/2} (w_i) r_2^{1/2} (w_i- w_{i-1}) |,
\end{eqnarray*}
where $\zeta(t) = \max\{ k\geq 0: w_k <t\}$.
Since $\rho_{\Theta_{\tilde{\kappa}, q, n}^2} ( (s_1, r_1), (s_2, r_2)) \leq \varepsilon$,
we have
\begin{eqnarray*}
| s_1^{1/2} ( w_i) r_1^{1/2} (w_i -w_{i-1})  - s_2^{1/2} (w_i) r_2^{1/2} (w_i-w_{i-1}) | &\leq & 2 \varepsilon \bar{\kappa}, \hspace{0.5cm}\forall i\geq 1, 
\end{eqnarray*}
and
\begin{eqnarray*}
&& | e^{-\int_0^T s_1(t) r_1(t- w_{\zeta(t)} ) dt/2} - e^{-\int_0^T s_2(t) r_2(t - w_{\zeta(t)} ) dt/2} |
\nonumber \\
&\leq & | 1 - e^{\int_0^T [ s_1(t) r_1(t- w_{\zeta(t)}) - s_2 (t) r_2 (t-w_{\zeta(t)}) ] dt/2} | 
\nonumber \\
&\leq &
2 \varepsilon \bar{\kappa}^3 T \sum_{i=0}^\infty \frac{ ( 2 \varepsilon \bar{\kappa}^3 T )^i}{ (i+1)!}.
\end{eqnarray*}
Consequently
for $j=1, 2, \cdots,$
\begin{eqnarray*}
&& \int^*_{0<w_1<\cdots <w_j <T} \sup_{\rho_{\Theta_{\tilde{\kappa}, q, n}^2}((s_1, r_1), (s_2, r_2)) 
\leq \varepsilon, (s_2, r_2)\in \Theta_{\tilde{\kappa}, q, n}^2 }
[ e^{-\int_0^T s_1(t) r_1(t-w_{\zeta(t)}) dt/2} 
\nonumber \\
&&\hspace{0.5cm}\times
\prod_{i=1}^j s_1^{1/2} (w_i) r_1^{1/2} (w_i- w_{i-1})
\nonumber \\
&&\hspace{0.5cm}
- e^{-\int_0^T s_2 (t) r_2 (t-w_{\zeta(t)}) dt/2} \prod_{i=1}^j s_2^{1/2} (w_i) r_2^{1/2} (w_i- w_{i-1}) ]^2 dw_1 \cdots dw_j
\nonumber \\
&\leq &
[ 2 \varepsilon \bar{\kappa}^{2 j+3} T \sum_{i=0}^\infty \frac{ (2 \varepsilon \bar{\kappa}^3 T )^i}{ (i+1)!}
+ 2 j \varepsilon \bar{\kappa}^{2 j-1} ]^2 \frac{T^j }{ j!},
\end{eqnarray*}
and for each $(s_1, r_1)\in \Theta_{\tilde{\kappa}, q, n}^2 (\varepsilon)$,
\begin{eqnarray*}
&& \Big\{ \sum_{j=0}^\infty \int^*_{0<w_1<\cdots < w_j <T} 
\sup_{ \rho_{\Theta_{\tilde{\kappa}, q, n}^2}((s_1, r_1), (s_2, r_2)) \leq \varepsilon, (s_2, r_2)\in \Theta_{\tilde{\kappa}, q, n}^2 } 
[ p_{s_1,r_1}^{1/2} (\{w_1,\cdots, w_j\})
\nonumber \\
&& \hspace{0.5cm} - p_{s_2, r_2}^{1/2} (\{ w_1,\cdots, w_j\}) ]^2 dw_1 \cdots dw_j \Big\}^{1/2}
\nonumber \\
&\leq & \varepsilon \Big\{ 
\sum_{j= 0}^\infty
[ 2 \bar{\kappa}^{2 j+3} T \sum_{i=0}^\infty \frac{ (2 \varepsilon \bar{\kappa}^3 T )^i}{ (i+1)!}
+ 2 j \bar{\kappa}^{2 j-1} ]^2 \frac{T^j }{ j!}
\Big\}^{1/2}
\nonumber \\
&\leq & \varepsilon C_{\tilde{\kappa} },
\end{eqnarray*}
where $C_{\tilde{\kappa} }$ is a absolute constant depending only on $\tilde{\kappa}$.  
(\ref{eq:4.9}) now follows from Lemma 2.1 of Ossiander (1987) and (\ref{eq:4.20})
since $C_{\tilde{\kappa}} \geq 1/2$. \hfill $\Box$

\begin{la} \label{la:a.6}
With the notation of Section 2, we have
\begin{eqnarray*}
H^B (\varepsilon, \tilde{\cal Z}_{\tilde{\kappa}, q, n}, \rho_{\tilde{\cal Z}_{\tilde{\kappa}, q, n}}) 
& \leq & H^B ( \frac{\varepsilon }{ 2 e^{\tau/2}}, {\cal F}_{\tilde{\kappa}, q, n}, \rho_{{\cal F}_{\tilde{\kappa}, q, n}})
\nonumber \\
&\leq & 2^{(q+2)/q} C_{\tilde{\kappa}}^{1/q} C_{\tilde{\kappa}, q} (\frac{ 2 e^{\tau/2} }{\varepsilon})^{1/q},
\hspace{0.5cm}\forall \varepsilon>0.
\end{eqnarray*}
\end{la}
{\sc Proof.} For $i=1,2$, let $f_i: {\cal N}\rightarrow R$ be a nonnegative function such that
\begin{displaymath}
\sum_{j=0}^\infty \int_{0<w_1<\cdots< w_j<T} f_i (\{w_1,\cdots, w_j\}) dw_1 \cdots dw_j <\infty.
\end{displaymath}
Define for $j=0,1, \cdots,$
\begin{displaymath}
A_{i,j} = \left\{ \{w_1,\cdots, w_j \}: f_i(\{w_1,\cdots, w_j \}) < e^{-\tau} p_{s,r}(\{w_1,\cdots, w_j \}) \right\},
\end{displaymath}
and
\begin{displaymath}
\tilde{f}_i (\{w_1,\cdots, w_j\}) = \left\{ \begin{array}{ll}
f_i (\{ w_1,\cdots, w_j\}), & \mbox{on $A_{i, j}^c$,} \\
e^{-\tau} p_{s,r}(\{w_1,\cdots, w_j\}), & \mbox{on $A_{i, j}$.}
\end{array}
\right.
\end{displaymath}
Letting $\tilde{Z}_{f_i}$ as in (\ref{eq:4.13}), we have
\begin{eqnarray*}
&& \sum_{j=0}^\infty \int_{0<w_1<\cdots<w_j<T} [\tilde{Z}_{f_1} (\{w_1,\cdots, w_j\} ) 
-\tilde{Z}_{f_2} (\{w_1,\cdots, w_j\} ) ]^2
\nonumber \\
&&\hspace{0.5cm}\times
p_{s,r} (\{w_1,\cdots, w_j\}) dw_1\cdots dw_j
\nonumber \\
&=&
4 \sum_{j=0}^\infty \int_{0<w_1<\cdots<w_j<T} \Big[ \log (\frac{ \tilde{f}_1^{1/2}(\{w_1,\cdots, w_j\} ) }{ p_{s,r}^{1/2}
(\{w_1,\cdots, w_j\}) } ) 
-\log(\frac{ \tilde{f}_2^{1/2} (\{w_1,\cdots, w_j\}) }{ p_{s,r}^{1/2} (\{w_1,\cdots, w_j\}) }) \Big]^2
\nonumber \\
&&\hspace{0.5cm}\times
p_{s,r} (\{w_1,\cdots, w_j\}) dw_1\cdots dw_j
\nonumber \\
&\leq & 4 e^\tau 
\sum_{j=0}^\infty \int_{0<w_1<\cdots<w_j<T} [ \tilde{f}_1^{1/2}(\{w_1,\cdots, w_j\} ) 
- \tilde{f}_2^{1/2} (\{w_1,\cdots, w_j\}) ]^2
dw_1\cdots dw_j.
\nonumber \\
\end{eqnarray*}
By dividing the integral into 
four parts: namely $A_{1,j}\cap A_{2, j}, A_{1,j}^c \cap A_{2,j}, A_{1,j}\cap A_{2,j}^c$ and $A_{1,j}^c \cap A_{2,j}^c$, we observe that
\begin{eqnarray*}
&& \sum_{j=0}^\infty \int_{0<w_1<\cdots<w_j<T} [ \tilde{f}_1^{1/2}(\{w_1,\cdots, w_j\} ) 
- \tilde{f}_2^{1/2} (\{w_1,\cdots, w_j\}) ]^2
dw_1\cdots dw_j
\nonumber \\
&\leq &
\sum_{j=0}^\infty \int_{0<w_1<\cdots<w_j<T} [ f_1^{1/2}(\{w_1,\cdots, w_j\} ) 
- f_2^{1/2} (\{w_1,\cdots, w_j\}) ]^2
dw_1\cdots dw_j.
\end{eqnarray*}
Now we conclude that
\begin{eqnarray*}
&& \Big\{ \sum_{j=0}^\infty \int_{0<w_1<\cdots<w_j<T} [ \tilde{Z}_{f_1} (\{w_1,\cdots, w_j\} ) 
- \tilde{Z}_{f_2} (\{w_1,\cdots, w_j\} ) ]^2
\nonumber \\
&&\hspace{0.5cm}\times
p_{s,r} (\{w_1,\cdots, w_j\}) dw_1\cdots dw_j \Big\}^{1/2}
\nonumber \\
&\leq &
2 e^{\tau/2} \Big\{ \sum_{j=0}^\infty \int_{0<w_1<\cdots<w_j<T} [ f_1^{1/2}(\{w_1,\cdots, w_j\} ) 
- f_2^{1/2} (\{w_1,\cdots, w_j\}) ]^2
dw_1\cdots dw_j \Big\}^{1/2}.
\end{eqnarray*}
Lemma \ref{la:a.6} now follows from Lemma \ref{la:a.3}. \hfill $\Box$

The statement of the next lemma can be found in Wong and Shen (1995), page 346, but its proof is not provided there.

\begin{la}[A Bernstein-type inequality] \label{la:4.6}
Let $Z_1, Z_2, \cdots$ be independent identically distributed random variables satisfying
\begin{displaymath}
E (|Z_1|^j) \leq \frac{ j! b^{j-2} \gamma}{2}, \hspace{0.5cm}\forall j\geq 2.
\end{displaymath}
Then 
\begin{displaymath}
P[ \frac{1}{n^{1/2}} \sum_{i=1}^n ( Z_i - EZ_i)  \geq t ] \leq \exp [ - \frac{ t^2}{ 4 (2 \gamma + b t n^{-1/2} ) } ],
\hspace{0.5cm}\forall t>0.
\end{displaymath}
\end{la}
{\sc Proof.}
The following proof is an adaption of the proof given in Bennett (1962), pages 36 to 38. 
Let $ c>0$ be a suitably chosen constant and ${\rm Var}(Z_i) = \sigma^2$. Then
\begin{displaymath}
E ( e^{c ( Z_i - EZ_i) } ) = 1 + c^2 \sum_{j=2}^\infty \frac{ c^{j-2} E[ (Z_i - EZ_i)^j] }{ j!}
= 1 + c^2 F, \hspace{0.5cm} \mbox{say}. 
\end{displaymath}
Since $1+ c^2 F\leq e^{ c^2 F}$, we have
\begin{displaymath}
E \{ \exp[  c\sum_{i=1}^n (Z_i - EZ_i) ] \} \leq e^{ n c^2 F}.
\end{displaymath}
Hence it follows from Markov's inequality that
\begin{eqnarray*}
P[ \sum_{i=1}^n (Z_i - EZ_i) \geq t \sqrt{n} ] &\leq & e^{-c t\sqrt{n} } E e^{ c \sum_{i=1}^n (Z_i - EZ_i) }
\nonumber \\
&\leq & e^{ c^2 n F - c t \sqrt{n}}
\nonumber \\
&\leq & e^{- t^2/(4 F)},
\end{eqnarray*}
by choosing $c$ such that $F = t/(2 c n^{1/2})$.
Now for $j\geq 2$,
\begin{eqnarray*}
| E( Z_1 - EZ_1)^j| &\leq & E( |Z_1 -EZ_1|^j )
\nonumber \\
&\leq &  \sum_{i=0}^j \frac{ j!}{i! (j-i)!} ( E|Z_1|^i) ( E|Z_1|^{j-i} )
\nonumber \\
&\leq & 2^j E(|Z_1|^j )
\nonumber \\
&\leq & j! 2^{j-1} b^{j-2} \gamma.
\end{eqnarray*}
Hence if $2 b c < 1$,
\begin{displaymath}
F \leq 2 \gamma \sum_{j=2}^\infty (2 bc)^{j-2}
= \frac{ 2 \gamma }{ 1 - 2 bc}.
\end{displaymath}
This implies that
\begin{displaymath}
\frac{ t}{2 c n^{1/2} } \leq \frac{ 2 \gamma }{ 1 - 2 b c},
\end{displaymath}
and consequently
\begin{displaymath}
c \geq \frac{ t}{ 4 \gamma n^{1/2} + 2 b t}.
\end{displaymath}
By taking
\begin{displaymath}
c = \frac{ t}{ 4 \gamma n^{1/2} + 2 b t},
\end{displaymath}
we observe that $2 bc <1$ and
\begin{displaymath}
P[ \frac{1}{n^{1/2}} \sum_{i=1}^n ( Z_i - EZ_i)  \geq t ] \leq \exp [ - \frac{ t^2}{ 4 (2 \gamma + b t n^{-1/2} ) } ].
\end{displaymath}
This proves Lemma \ref{la:4.6}.\hfill $\Box$

For $g: {\cal N}\rightarrow R$, let
\begin{equation}
\nu_n (g) = \frac{1}{n^{1/2}} \sum_{i=1}^n [ g(\{w_{i,1},\cdots, w_{i, N_i(T)} \})
- E_{s,r} g(\{w_{i,1},\cdots, w_{i, N_i(T)} \}) ],
\label{eq:3.22}
\end{equation}
assuming that the right hand side exists 
and $\{ w_{i,1}, \cdots, w_{i, N_i(T)} \}$ is as in Section 3.
We observe from Lemma 6 of Wong and Shen (1995) that
\begin{equation}
P_{s,r} [ \nu_n( \tilde{Z}_f ) \geq t ]
\leq
\exp[ - \frac{ t^2}{ 8 (8 c_0 \| f^{1/2} - p_{s,r}^{1/2} \|_2^2 + 2 t n^{-1/2} ) } ],
\hspace{0.5cm}\forall t>0,
\label{eq:3.21}
\end{equation}
where
$\tilde{Z}_f$ is as in (\ref{eq:4.13}),
\begin{displaymath}
c_0= (e^{\tau/2} - 1 -\frac{\tau}{2} )/(1 - e^{-\tau/2} )^2.
\end{displaymath}
The next lemma is motivated by Theorem 3 of Shen and Wong (1994) and Lemma 7 of Wong and Shen (1995). As the proof of the latter lemma is only briefly sketched
in Wong and Shen (1995), page 348, a detailed proof of Lemma \ref{la:4.10} is given below.

\begin{la} \label{la:4.10}
For any $t>0, 0<\gamma<1$ and $M>0$, let
\begin{displaymath}
\psi (M, t^2, n)= \frac{ M^2}{ 16 ( 8c_0 t^2 + M n^{-1/2}) }.
\end{displaymath}
Assume that
\begin{eqnarray}
H^B (\frac{t}{10}, {\cal F}_{\tilde{\kappa}, q, n}, \rho_{{\cal F}_{\tilde{\kappa}, q, n}}) &\leq & \frac{\gamma}{4} \psi (M, t^2, n),
\label{eq:4.88} \\
\int_{\gamma M/(32 n^{1/2})}^{e^{\tau/2}t/5} [ H^B (\frac{ x}{2 e^{\tau/2}}, {\cal F}_{\tilde{\kappa}, q, n}, \rho_{{\cal F}_{\tilde{\kappa}, q, n}} ) ]^{1/2} dx
& \leq & \frac{ M \gamma^{3/2} }{ 2^{10} }.
\label{eq:4.89}
\end{eqnarray}
Then
\begin{displaymath}
P^*_{s,r} [ \sup_{\| p^{1/2}_{s_1, r_1} - p^{1/2}_{s, r} \|_2 \leq t, p_{s_1,r_1} \in {\cal F}_{\tilde{\kappa}, q, n}}
\nu_n(  \tilde{Z}_{p_{s_1, r_1}} ) \geq M]
\leq 3 e^{-(1-\gamma ) \psi( M, t^2, n)},
\end{displaymath}
where $\nu_n(.)$ is as in (\ref{eq:3.22}).
\end{la}
{\sc Proof.}
Without loss of generality, we can assume that
\begin{equation}
3 e^{-(1-\gamma ) \psi( M, t^2, n)} \leq 1,
\label{eq:4.25}
\end{equation}
else Lemma \ref{la:4.10} is trivial.
We observe from Lemma \ref{la:a.6} that 
\begin{displaymath}
H^B (\gamma, \tilde{\cal Z}_{\tilde{\kappa}, q, n}, \rho_{\tilde{\cal Z}_{\tilde{\kappa}, q, n} }) 
\leq H^B (\frac{\gamma }{ 2 e^{\tau/2}}, {\cal F}_{\tilde{\kappa}, q, n},
\rho_{{\cal F}_{\tilde{\kappa}, q, n}}) <\infty.
\end{displaymath}
For any $\delta_{n,0} > \delta_{n,1} > \cdots > \delta_{n,n_0} > 0$, there exist 
${\cal F}_{n,j}, j=0,\cdots, n_0$, with 
\begin{displaymath}
| {\cal F}_{n,j} | = \exp[ H^B (\frac{ \delta_{n,j} }{2 e^{\tau/2}}, {\cal F}_{\tilde{\kappa}, q, n}, \rho_{{\cal F}_{\tilde{\kappa}, q, n}} ) ],
\end{displaymath}
such that
for each $p_{s_1, r_1} \in {\cal F}_{\tilde{\kappa}, q, n}$ one can find
$f_j^L$, $f_j^U \in {\cal F}_{n,j}$ such that
\begin{displaymath}
f_j^L(\{w_1,\cdots, w_{N(T)} \} ) \leq p_{s_1, r_1}(\{w_1,\cdots, w_{N(T)} \} ) 
\leq f_j^U(\{w_1,\cdots, w_{N(T)} \} ), \hspace{0.5cm} \mbox{a. s.},
\end{displaymath}
and
\begin{eqnarray*}
&& \Big\{ \sum_{j_1=0}^\infty \int_{0<w_1<\cdots < w_{j_1} <T} [ \tilde{Z}_{f_j^U}( \{w_1,\cdots, w_{j_1}\})
- \tilde{Z}_{f_j^L}( \{w_1,\cdots, w_{j_1} \}) ]^2 
\nonumber \\
&& \hspace{0.5cm}\times p_{s,r} (\{w_1,\cdots, w_{j_1} \}) dw_1\cdots dw_{j_1} \Big\}^{1/2} 
\nonumber \\
&\leq &
2 e^{\tau/2} \Big\{ \sum_{j_1=0}^\infty \int_{0<w_1<\cdots < w_{j_1}<T} [f^U_j ( \{w_1,\cdots, w_{j_1}\})^{1/2}
\nonumber \\
&&\hspace{0.5cm}
- f_j^L (\{w_1,\cdots, w_{j_1}\})^{1/2} ]^2 dw_1 \cdots dw_{j_1} \Big\}^{1/2}
\nonumber \\
&\leq &  \delta_{n,j}.
\end{eqnarray*}
Here $n_0$ is a nonnegative integer to be suitably chosen later.
Define for $k=0,\cdots, n_0$,
\begin{eqnarray*}
u_k (\{ w_1,\cdots, w_{N(T)}\}) &=& \min_{0\leq j\leq k} f_j^U (\{w_1,\cdots, w_{N(T)}\}),
 \nonumber \\
l_k(\{ w_1,\cdots, w_{N(T)}\}) &=& \max_{0\leq j\leq k} f_j^L (\{w_1,\cdots, w_{N(T)}\}).
 \end{eqnarray*}
Then
$\tilde{Z}_{l_k} \leq \tilde{Z}_{p_{s_1, r_1}} \leq \tilde{Z}_{u_k}$,
$0 \leq \tilde{Z}_{u_{k+1} } 
 - \tilde{Z}_{ l_{k+1} }
 \leq 
\tilde{Z}_{ u_k} - \tilde{Z}_{ l_k}$ a. s.,
and
\begin{displaymath}
[ E_{s,r} (\tilde{Z}_{u_k} - \tilde{Z}_{l_k} )^2 ]^{1/2}
\leq 
[ E_{s,r} (\tilde{Z}_{f_k^U} - \tilde{Z}_{f_k^L} )^2 ]^{1/2}
\leq \delta_{n,k}, \hspace{0.5cm}\forall k=0,\cdots, n_0.
\end{displaymath}

If $n_0=0$, define $B_0= {\cal N}$.
If $n_0\geq 1$, let $a_1> a_2> \cdots > 0$ be a sequence of constants
and define
\begin{eqnarray*}
B_0 &=& \{ \tilde{Z}_{u_0} - \tilde{Z}_{l_0} \geq a_1\},
\nonumber \\
B_k &=& \{ \tilde{Z}_{u_k} - \tilde{Z}_{l_k} \geq a_{k+1},
\tilde{Z}_{u_j} - \tilde{Z}_{l_j} < a_{j+1},
j=0,\cdots, k-1\},
\hspace{0.5cm}\forall k=1,\cdots, n_0-1, 
\nonumber \\
B_{n_0} &=& ( \cup_{k=0}^{n_0-1} B_k )^c.
\end{eqnarray*}
Note that $\{ B_k: k=0,\cdots, n_0\}$ forms a partition of ${\cal N}$.
Consequently writing ${\bf 1}_{B_k}$ to denote the indicator function of $B_k$, we have
\begin{eqnarray*}
\tilde{Z}_{p_{s_1, r_1}}
&=& \tilde{Z}_{u_0} + \sum_{k=0}^{n_0} ( \tilde{Z}_{u_k} {\bf 1}_{B_k}
- \tilde{Z}_{u_0} {\bf 1}_{B_k} )
+ \tilde{Z}_{p_{s_1, r_1}} - \sum_{k=0}^{n_0} \tilde{Z}_{u_k} {\bf 1}_{B_k}
\nonumber \\
&=&
\tilde{Z}_{u_0}
+\sum_{j=1}^{n_0} ( \tilde{Z}_{u_j}
- \tilde{Z}_{u_{j-1}} ) {\bf 1}_{\cup_{j\leq k\leq n_0} B_k} 
+ \sum_{k=0}^{n_0} ( \tilde{Z}_{p_{s_1, r_1}} - \tilde{Z}_{u_k}  ) {\bf 1}_{B_k}.
\end{eqnarray*}
Let $\eta_1, \cdots, \eta_{n_0+1}$ be strictly positive constants such that $2 \eta_1 
+ \cdots + 2 \eta_{n_0} +\eta_{n_0+1} \leq \gamma M/8$.
Then
\begin{eqnarray}
&& P^*_{s,r} [ \sup_{\| p_{s_1, r_1}^{1/2} - p_{s,r}^{1/2} \|_2 \leq t, p_{s_1, r_1} \in {\cal F}_{\tilde{\kappa}, q, n} }
\nu_n (\tilde{Z}_{p_{s_1, r_1}} ) \geq M ]
\nonumber \\
&\leq &
P_{s,r}^* [ \sup_{\| p_{s_1, r_1}^{1/2} - p_{s,r}^{1/2} \|_2 \leq t, p_{s_1, r_1} \in {\cal F}_{\tilde{\kappa}, q, n} }
\nu_n (\tilde{Z}_{u_0} ) \geq M - \frac{\gamma M}{4}  ]
\nonumber \\
&& + P_{s,r}^* [ \sup_{\| p_{s_1, r_1}^{1/2} - p_{s,r}^{1/2} \|_2 \leq t, p_{s_1, r_1} \in {\cal F}_{\tilde{\kappa}, q, n} }
\nu_n (\sum_{j=1}^{n_0} ( \tilde{Z}_{u_j }
- \tilde{Z}_{u_{j-1}} ) {\bf 1}_{\cup_{j\leq k\leq n_0} B_k} ) \geq \sum_{j=1}^{n_0} \eta_j ]
\nonumber \\
&& + P_{s,r}^* [ \sup_{\| p_{s_1, r_1}^{1/2} - p_{s,r}^{1/2} \|_2 \leq t, p_{s_1, r_1} \in {\cal F}_{\tilde{\kappa}, q, n} }
\nu_n (
\sum_{k=0}^{n_0} ( \tilde{Z}_{p_{s_1, r_1}} - \tilde{Z}_{u_k}  ) {\bf 1}_{B_k} ) \geq \frac{\gamma M}{8} +\sum_{k=1}^{n_0+1} \eta_k ]
\nonumber \\
&\leq &
|{\cal F}_{n,0} | \sup_{\| p_{s_1, r_1}^{1/2} - p_{s,r}^{1/2} \|_2 \leq t, p_{s_1, r_1} \in {\cal F}_{\tilde{\kappa}, q, n} }
P_{s,r} [\nu_n (\tilde{Z}_{u_0} ) \geq M - \frac{\gamma M}{4}  ]
\nonumber \\
&& + \sum_{j=1}^{n_0} P_{s,r}^* [ \sup_{\| p_{s_1, r_1}^{1/2} - p_{s,r}^{1/2} \|_2 \leq t, p_{s_1, r_1} \in {\cal F}_{\tilde{\kappa}, q, n} }
\nu_n ( ( \tilde{Z}_{u_j }
- \tilde{Z}_{ u_{j-1}} ) {\bf 1}_{\cup_{j\leq k\leq n_0} B_k} ) \geq \eta_j ]
\nonumber \\
&& + \sum_{k=0}^{n_0-1} P_{s,r}^* [ \sup_{\| p_{s_1, r_1}^{1/2} - p_{s,r}^{1/2} \|_2 \leq t, p_{s_1, r_1} \in {\cal F}_{\tilde{\kappa}, q, n} }
\nu_n (
( \tilde{Z}_{p_{s_1, r_1}} - \tilde{Z}_{u_k}  ) {\bf 1}_{B_k} ) \geq \eta_{k+1} ]
\nonumber \\
&& + P_{s,r}^* [ \sup_{\| p_{s_1, r_1}^{1/2} - p_{s,r}^{1/2} \|_2 \leq t, p_{s_1, r_1} \in {\cal F}_{\tilde{\kappa}, q, n} }
\nu_n (
( \tilde{Z}_{p_{s_1, r_1}} - \tilde{Z}_{u_{n_0}}  ) {\bf 1}_{B_{n_0} } ) \geq \frac{\gamma M}{8} + \eta_{n_0+1} ]
\nonumber \\
&\leq &
|{\cal F}_{n,0}| \sup_{\| p_{s_1, r_1}^{1/2} - p_{s,r}^{1/2} \|_2 \leq t, p_{s_1, r_1} \in {\cal F}_{\tilde{\kappa}, q, n} }
P_{s,r} [\nu_n (\tilde{Z}_{u_0} ) \geq M - \frac{\gamma M}{4}  ]
\nonumber \\
&& + \sum_{j=1}^{n_0} (\prod_{l=0}^j | {\cal F}_{n,l} | ) \sup_{\| p_{s_1, r_1}^{1/2} - p_{s,r}^{1/2} \|_2 \leq t, p_{s_1, r_1} \in {\cal F}_{\tilde{\kappa}, q, n} }
P_{s,r} [ \nu_n ( ( \tilde{Z}_{u_j}
- \tilde{Z}_{u_{j-1}} ) {\bf 1}_{\cup_{j\leq k\leq n_0} B_k} ) \geq \eta_j ]
\nonumber \\
&& + \sum_{k=0}^{n_0-1} P_{s,r}^* [ \sup_{\| p_{s_1, r_1}^{1/2} - p_{s,r}^{1/2} \|_2 \leq t, p_{s_1, r_1} \in {\cal F}_{\tilde{\kappa}, q, n} }
\nu_n (
( \tilde{Z}_{p_{s_1, r_1}} - \tilde{Z}_{u_k}  ) {\bf 1}_{B_k} ) \geq \eta_{k+1} ]
\nonumber \\
&& + P_{s,r}^* [ \sup_{\| p_{s_1, r_1}^{1/2} - p_{s,r}^{1/2} \|_2 \leq t, p_{s_1, r_1} \in {\cal F}_{\tilde{\kappa}, q, n} }
\nu_n (
( \tilde{Z}_{p_{s_1, r_1}} - \tilde{Z}_{u_{n_0}}  ) {\bf 1}_{B_{n_0} } ) \geq \frac{\gamma M}{8} + \eta_{n_0+1} ].
\label{eq:4.66}
\end{eqnarray}
Since $H^B$ may be replaced by a larger continuous function at the expense of an arbitrarily small increase in these values,
we shall follow Alexander (1984), page 1045, and assume for the rest of this proof that
$H^B(., {\cal F}_{\tilde{\kappa}, q, n},\rho_{{\cal F}_{\tilde{\kappa}, q, n} } )$ is continuous and strictly decreasing from $\infty$ to $0$ on $(0, a]$ for some $a$.
We define
\begin{displaymath}
\delta_{n,0} = \inf \{ x: H^B( \frac{x }{ 2 e^{\tau/2}}, {\cal F}_{\tilde{\kappa}, q, n}, \rho_{{\cal F}_{\tilde{\kappa}, q, n} } ) 
\leq \frac{\gamma }{4} \psi( M, t^2, n) \}.
\end{displaymath}

{\sc Case 1.} Suppose that $\delta_{n,0} > \gamma M/( 8 n^{1/2})$.
Then we further define
\begin{eqnarray*}
\delta_{n, j} &=& \max \{ \frac{\gamma M}{ 8 n^{1/2} }, \sup\{ x\leq \frac{\delta_{n,j-1} }{2}: 
H^B (\frac{x }{2 e^{\tau/2}}, {\cal F}_{\tilde{\kappa}, q, n}, \rho_{{\cal F}_{\tilde{\kappa}, q, n} } )
\geq 
\nonumber \\
&& \hspace{0.5cm}
4 H^B (\frac{ \delta_{n,j-1} }{ 2 e^{\tau/2}}, {\cal F}_{\tilde{\kappa}, q, n}, \rho_{{\cal F}_{\tilde{\kappa}, q, n} } ) \} \}, 
\hspace{0.5cm} \forall j= 1, 2,\cdots,
\nonumber \\
n_0 &=& \min\{ j: \delta_{n,j} = \frac{ \gamma M}{ 8 n^{1/2}} \},
\nonumber \\
\eta_j &=& 4 \delta_{n, j-1} [ \frac{ \sum_{i=0}^j H^B (\delta_{n,i} e^{-\tau/2}/2, {\cal F}_{\tilde{\kappa}, q, n}, 
\rho_{{\cal F}_{\tilde{\kappa}, q, n} } ) }{ \gamma } ]^{1/2},
\hspace{0.5cm}\forall j=1,\cdots, n_0+1,
\nonumber \\
a_j &=& \frac{ 8 n^{1/2} \delta_{n, j-1}^2 }{ \eta_j},
\hspace{0.5cm}\forall j=1,\cdots, n_0.
\end{eqnarray*}
Thus 
$\delta_{n, n_0} = \delta_{n, n_0+1} = \gamma M/(8 n^{1/2})$,
\begin{eqnarray*}
H^B ( \frac{\delta_{n,0} }{ 2 e^{\tau/2}}, {\cal F}_{\tilde{\kappa}, q, n}, \rho_{{\cal F}_{\tilde{\kappa}, q, n} } ) &=& \frac{ \gamma \psi( M, t^2, n) }{4},
\nonumber \\
\frac{ \delta_{n,0} }{ 2 e^{\tau/2} } & \leq &  \frac{ t}{10}.
\end{eqnarray*}
We observe from Lemma 3.1 of Alexander (1984) and (\ref{eq:4.89}) that
\begin{eqnarray*}
2 \sum_{j=1}^{n_0+1} \eta_j &= & 
\frac{ 8 }{\gamma^{1/2} } \sum_{j=1}^{n_0+1} \delta_{n, j-1} [ \sum_{i=0}^j H^B (\frac{ \delta_{n,i} }{ 2 e^{\tau/2}}, 
{\cal F}_{\tilde{\kappa}, q, n}, \rho_{{\cal F}_{\tilde{\kappa}, q, n} } ) ]^{1/2}
\nonumber \\
&\leq  &
\frac{ 8 }{\gamma^{1/2} } \sum_{j=1}^{n_0-1} \delta_{n,j-1} [ \frac{4}{3} H^B (\frac{\delta_{n,j} }{ 2 e^{\tau/2}}, 
{\cal F}_{\tilde{\kappa}, q, n}, \rho_{{\cal F}_{\tilde{\kappa}, q, n} } ) ]^{1/2}
\nonumber \\
& &
+ \frac{ 8 }{\gamma^{1/2} } \sum_{j=n_0}^{n_0+1} \delta_{n, j-1} [ 4 H^B (\frac{ \delta_{n,j} }{ 2 e^{\tau/2}}, 
{\cal F}_{\tilde{\kappa}, q, n}, \rho_{{\cal F}_{\tilde{\kappa}, q, n}} ) ]^{1/2}
\nonumber \\
&\leq &
\frac{ 16 }{\gamma^{1/2} } \sum_{j=0}^{n_0} \delta_{n,j} [ H^B (\frac{ \delta_{n,j+1} }{2 e^{\tau/2}}, 
{\cal F}_{\tilde{\kappa}, q, n}, \rho_{{\cal F}_{\tilde{\kappa}, q, n} } ) ]^{1/2}
\nonumber \\
&\leq &
\frac{ 2^7 }{\gamma^{1/2} } \int_{\gamma M/(32 n^{1/2})}^{\delta_{n,0} } [ H^B (\frac{x }{ 2 e^{\tau/2}}, 
{\cal F}_{\tilde{\kappa}, q, n}, \rho_{{\cal F}_{\tilde{\kappa}, q, n} } ) ]^{1/2} dx
\nonumber \\
&\leq & \frac{ \gamma M}{8}.
\end{eqnarray*}

Now we observe from (\ref{eq:3.21}) that
\begin{eqnarray*}
&& |{\cal F}_{n,0} | \sup_{\| p_{s_1, r_1}^{1/2} - p_{s,r}^{1/2} \|_2 \leq t, p_{s_1, r_1} \in {\cal F}_{\tilde{\kappa}, q, n} }
P_{s,r} [\nu_n (\tilde{Z}_{u_0} ) \geq (1 - 2^{-2} \gamma ) M ]
\nonumber \\
&\leq & 
e^{ H^B (\delta_{n,0} e^{-\tau/2}/2, {\cal F}_{\tilde{\kappa}, q, n}, \rho_{{\cal F}_{\tilde{\kappa}, q, n} } )}
\exp\{ - \frac{ (1 - 2^{-2} \gamma )^2 M^2 }{
8 [ 8 c_0 ( 2^{-1} e^{-\tau/2} \delta_{n,0} + t)^2  
+ 2 n^{-1/2} (1- 2^{-2} \gamma ) M ] } \}
\nonumber \\
&\leq &
e^{ H^B (\delta_{n,0} e^{-\tau/2}/2, {\cal F}_{\tilde{\kappa}, q, n}, \rho_{{\cal F}_{\tilde{\kappa}, q, n} } )}
\exp\{ - \frac{ (1 - 2^{-2} \gamma )^2 M^2 }{
16 [ 8 c_0 t^2  
+ n^{-1/2} (1- 2^{-2} \gamma ) M ] } \}
\nonumber \\
&\leq &
\exp\{ [\frac{\gamma }{4} - (1 - \frac{ \gamma }{4} )^2 ]
\psi (M, t^2, n) \}
\nonumber \\
&\leq &
\exp[ - (1 - \frac{ 3 \gamma }{4} ) 
\psi (M, t^2, n) ],
\end{eqnarray*}
since 
\begin{eqnarray*}
\| u_0^{1/2} - p_{s,r}^{1/2} \|_2 
& \leq & \| u_0^{1/2} - p_{s_1, r_1}^{1/2} \|_2 + \| p_{s_1, r_1}^{1/2} - p_{s,r}^{1/2} \|_2
\nonumber \\
&\leq & \frac{ \delta_{n,0} }{ 2 e^{\tau/2} } + t
\nonumber \\
&\leq & \frac{ 11 t}{ 10}.
\end{eqnarray*}

Next we have
\begin{eqnarray*}
{\rm Var}_{s,r} [ (\tilde{Z}_{u_j} - \tilde{Z}_{u_{j-1}} ) {\bf 1}_{\cup_{j\leq k \leq n_0} B_k} ]
&\leq & E_{s,r} [( \tilde{Z}_{u_{j-1}} -\tilde{Z}_{l_{j-1}} )^2 ]
\nonumber\\
&\leq & \delta_{n, j-1}^2, \hspace{0.5cm}\forall j=1,\cdots, n_0,
\end{eqnarray*}
and
$-a_j \leq \tilde{Z}_{l_{j-1}} - \tilde{Z}_{u_{j-1}} \leq \tilde{Z}_{u_j} - \tilde{Z}_{u_{j-1}} \leq 0$
on $\cup_{j\leq k\leq n_0} B_k$.
Hence by the one sided version of Bernstein's inequality [see Bennett (1962), page 38],
\begin{displaymath}
P_{s,r} (\nu_n ( (\tilde{Z}_{u_j} - \tilde{Z}_{u_{j-1}} ) {\bf 1}_{\cup_{j\leq k\leq n_0} B_k} ) \geq \eta_j)
\leq \exp [ -\frac{ \eta_j^2 }{ 2( \delta_{n, j-1}^2 + a_j \eta_j/(3 n^{1/2}) )} ],
\hspace{0.5cm}\forall j=1,\cdots, n_0.
\end{displaymath}
Since 
\begin{displaymath}
\frac{ \eta_j^2 }{ 2( \delta_{n, j-1}^2 + a_j \eta_j/(3 n^{1/2}) )} 
= \frac{ 3 \eta_j^2 }{ 22 \delta_{n, j-1}^2}
\geq \frac{ 2}{ \gamma}  \sum_{i=0}^j H^B (\frac{ \delta_{n,i}}{ 2 e^{\tau/2}}, {\cal F}_{\tilde{\kappa}, q, n}, \rho_{{\cal F}_{\tilde{\kappa}, q, n} } ),
\end{displaymath}
we have
\begin{eqnarray*}
&& \sum_{j=1}^{n_0} (\prod_{l=0}^j | {\cal F}_{n,l} | ) \sup_{\| p_{s_1, r_1}^{1/2} - p_{s,r}^{1/2} \|_2 
\leq t, p_{s_1, r_1} \in {\cal F}_{\tilde{\kappa}, q, n} }
P_{s,r} [ \nu_n ( ( \tilde{Z}_{u_j}
- \tilde{Z}_{u_{j-1}} ) {\bf 1}_{\cup_{j\leq k\leq n_0} B_k} ) \geq \eta_j ]
\nonumber \\
&\leq &
\sum_{j=1}^{n_0} \exp [ \sum_{i =0}^j H^B (\frac{ \delta_{n,i} }{ 2 e^{\tau/2}}, {\cal F}_{\tilde{\kappa}, q, n}, 
\rho_{{\cal F}_{\tilde{\kappa}, q, n} } ) 
-\frac{ \eta_j^2 }{ 2( \delta_{n, j-1}^2 + a_j \eta_j/(3 n^{1/2}) )} ]
\nonumber \\
&\leq &
\sum_{j=1}^{n_0-1} \exp [ -\frac{ 2(1 - \gamma)}{\gamma } \sum_{i=0}^j H^B (\frac{ \delta_{n, i} }{ 2 e^{\tau/2}}, 
{\cal F}_{\tilde{\kappa}, q, n}, \rho_{{\cal F}_{\tilde{\kappa}, q, n} } ) ]
\nonumber \\
&&\hspace{0.5cm}
+ \exp[-\frac{2(1 -\gamma) }{\gamma } \sum_{i=0}^{n_0} H^B (\frac{ \delta_{n,i}}{ 2  e^{\tau/2}}, {\cal F}_{\tilde{\kappa}, q, n}, 
\rho_{{\cal F}_{\tilde{\kappa}, q, n} } )]
\nonumber \\
&\leq &
\sum_{j=1}^{n_0-1} \exp [ -\frac{ 2 (1 - \gamma )( \sum_{i=0}^j 4^i) }{\gamma } H^B (\frac{ \delta_{n,0}}{2 e^{\tau/2}}, 
{\cal F}_{\tilde{\kappa}, q, n}, \rho_{{\cal F}_{\tilde{\kappa}, q, n }} ) ]
\nonumber \\
&&\hspace{0.5cm}
+ \exp[-\frac{ 2 ( 1 -\gamma ) (4^{n_0-1} +\sum_{i=0}^{n_0-1} 4^i ) }{\gamma } H^B (\frac{ \delta_{n,0}}{2 e^{\tau/2}}, 
{\cal F}_{\tilde{\kappa}, q, n}, \rho_{{\cal F}_{\tilde{\kappa}, q, n} } )]
\nonumber \\
&= &
\sum_{j=1}^{n_0-1} \exp [ -\frac{ (1 - \gamma )( \sum_{i=0}^j 4^i) }{2 } \psi (M, t^2, n) ]
\nonumber \\
&&\hspace{0.5cm}
+ \exp[-\frac{ ( 1 -\gamma ) (4^{n_0-1} +\sum_{i=0}^{n_0-1} 4^i ) }{2 } \psi( M, t^2, n) ].
\end{eqnarray*}
Thus it follows from (\ref{eq:4.25}) that
\begin{eqnarray*}
&& \sum_{j=1}^{n_0} (\prod_{l=0}^j | {\cal F}_{n, l} | ) \sup_{\| p_{s_1, r_1}^{1/2} - p_{s,r}^{1/2} \|_2 \leq t, 
p_{s_1, r_1} \in {\cal F}_{\tilde{\kappa}, q, n} }
P_{s,r} [ \nu_n ( ( \tilde{Z}_{u_j}
- \tilde{Z}_{u_{j-1}} ) {\bf 1}_{\cup_{j\leq k\leq n_0} B_k} ) \geq \eta_j ]
\nonumber \\
&\leq &
2 \exp [ - (1 - \gamma ) 
\psi (M, t^2, n) ].
\end{eqnarray*}
Next using Markov's inequality, we observe that for $0\leq k\leq n_0-1$,
\begin{displaymath}
P_{s,r} ( B_k) \leq P_{s,r} ( \tilde{Z}_{u_k} - \tilde{Z}_{l_k} \geq a_{k+1} )
\leq  \frac{ E_{s,r} (\tilde{Z}_{u_k} - \tilde{Z}_{l_k} )^2 }{ a_{k+1}^2 }
\leq \frac{ \delta_{n,k}^2}{ a^2_{k+1}}.
\end{displaymath}
Hence
\begin{eqnarray*}
&& \sup_{\| p_{s_1, r_1}^{1/2} - p_{s,r}^{1/2} \|_2 \leq t, p_{s_1, r_1} \in {\cal F}_{\tilde{\kappa}, q, n} }
\nu_n ((\tilde{Z}_{p_{s_1, r_1} } - \tilde{Z}_{u_k} ) {\bf 1}_{B_k} ) 
\nonumber \\
& \leq &
n^{1/2}
\sup_{\| p_{s_1, r_1}^{1/2} - p_{s,r}^{1/2} \|_2 \leq t, p_{s_1, r_1} \in {\cal F}_{\tilde{\kappa}, q, n} }
E_{s,r} [ ( \tilde{Z}_{u_k} - \tilde{Z}_{l_k} ) {\bf 1}_{B_k} ]
\nonumber \\
&\leq & n^{1/2} \sup_{\| p_{s_1, r_1}^{1/2} - p_{s,r}^{1/2} \|_2 \leq t, p_{s_1, r_1} \in {\cal F}_{\tilde{\kappa}, q, n} }
[ E_{s,r} ( \tilde{Z}_{u_k} - \tilde{Z}_{l_k} )^2 P_{s,r}( B_k ) ]^{1/2}
\nonumber \\
& \leq & \frac{ n^{1/2} \delta^2_{n,k} }{ a_{k+1} }
\nonumber \\
&= & \frac{ \eta_{k+1}}{ 8}, \hspace{0.5cm}\forall k=0,\cdots, n_0-1.
\end{eqnarray*}
This implies that
\begin{displaymath}
\sum_{k=0}^{n_0-1} P_{s,r}^* [ \sup_{\| p_{s_1, r_1}^{1/2} - p_{s,r}^{1/2} \|_2 \leq t, 
p_{s_1, r_1} \in {\cal F}_{\tilde{\kappa}, q, n} }
\nu_n (
( \tilde{Z}_{p_{s_1, r_1}} - \tilde{Z}_{u_k}  ) {\bf 1}_{B_k} ) \geq \eta_{k+1} ]
= 0.
\end{displaymath}
Finally,
\begin{eqnarray*}
&& \sup_{\| p_{s_1, r_1}^{1/2} - p_{s,r}^{1/2} \|_2 \leq t, p_{s_1, r_1} \in {\cal F}_{\tilde{\kappa}, q, n} }
\nu_n (
( \tilde{Z}_{p_{s_1, r_1}} - \tilde{Z}_{u_{n_0}}  ) {\bf 1}_{B_{n_0} } )
\nonumber \\
&\leq &
\sup_{\| p_{s_1, r_1}^{1/2} - p_{s,r}^{1/2} \|_2 \leq t, p_{s_1, r_1} \in {\cal F}_{\tilde{\kappa}, q, n} }
n^{1/2} E_{s, r} | \tilde{Z}_{p_{s_1, r_1}} - \tilde{Z}_{u_{n_0}} | 
\nonumber \\
&\leq & n^{1/2} \delta_{n, n_0}
\nonumber \\
&= &
\frac{\gamma M}{ 8 },
\end{eqnarray*}
and consequently
\begin{displaymath}
P_{s,r}^* [ \sup_{\| p_{s_1, r_1}^{1/2} - p_{s,r}^{1/2} \|_2 \leq t, p_{s_1, r_1} \in {\cal F}_{\tilde{\kappa}, q, n} }
\nu_n (
( \tilde{Z}_{p_{s_1, r_1}} - \tilde{Z}_{u_{n_0}}  ) {\bf 1}_{B_{n_0} } )
\geq \frac{\gamma M}{8} +\eta_{n_0+1} ] = 0.
\end{displaymath}

Thus we conclude from (\ref{eq:4.66}) that
\begin{displaymath}
P_{s,r}^* [ \sup_{\| p^{1/2}_{s_1, r_1} - p^{1/2}_{s, r} \|_2 \leq t, p_{s_1,r_1} \in {\cal F}_{\tilde{\kappa}, q, n} }
\nu_n(  \tilde{Z}_{p_{s_1, r_1}} ) \geq M]
\leq 3 e^{-(1-\gamma ) \psi( M, t^2, n)}.
\end{displaymath}

{\sc Case 2.} Suppose that $\delta_{n,0} \leq \gamma M/( 8 n^{1/2})$.
Then define $n_0=0$ and $\eta_1 = \gamma M/16$.
As in Case 1, we have
\begin{displaymath}
|{\cal F}_{n,0}| \sup_{\| p_{s_1, r_1}^{1/2} - p_{s,r}^{1/2} \|_2 \leq t, p_{s_1, r_1} \in {\cal F}_{\tilde{\kappa}, q, n} }
P_{s,r} [\nu_n (\tilde{Z}_{u_0} ) \geq (1 - \frac{\gamma}{4} ) M ]
\leq
\exp[ - (1 - \frac{ 3\gamma}{4} ) 
\psi (M, t^2, n) ],
\end{displaymath}
and
\begin{displaymath}
P_{s,r}^* [ \sup_{\| p_{s_1, r_1}^{1/2} - p_{s,r}^{1/2} \|_2 \leq t, p_{s_1, r_1} \in {\cal F}_{\tilde{\kappa}, q, n} }
\nu_n (
( \tilde{Z}_{p_{s_1, r_1}} - \tilde{Z}_{u_{n_0}}  ) {\bf 1}_{B_{n_0} } ) \geq \frac{\gamma M}{8} + \eta_{n_0+1} ]
= 0.
\end{displaymath}
This completes the proof of Lemma \ref{la:4.10}.\hfill $\Box$

{\sc Proof of Proposition \ref{pn:4.1}.} 
Without loss of generality, we can assume that
\begin{equation}
4 \exp [ - \frac{ n\varepsilon^2 }{ 2^7 (250) } ] \leq 1.
\label{eq:a.56}
\end{equation}
We observe that
(\ref{eq:4.90}) holds with $\varepsilon$ replaced by any $s$ such that  $\varepsilon \leq s \leq 1$.
Let $\gamma = 1/2$, $e^{\tau/2} = 5$, $t = \sqrt{2} s$ and $M= \gamma n^{1/2} s^2/2$.
Then it follows from Lemma \ref{la:a.6} that
\begin{eqnarray*}
\int_{ \gamma M/(32 n^{1/2}) }^{e^{\tau/2} t/5} [H^B (\frac{x}{2 e^{\tau/2}}, {\cal F}_{\tilde{\kappa}, q, n}, 
\rho_{{\cal F}_{\tilde{\kappa}, q, n} }) ]^{1/2} dx
&\leq &
\int_{s^2/2^8}^{\sqrt{2} s} [ 2^{(q+2)/q}  C_{\tilde{\kappa}}^{1/q} C_{\tilde{\kappa}, q} (\frac{2 e^{\tau/2} }{ x} )^{1/q}]^{1/2} dx
\nonumber \\
& \leq & \frac{ n^{1/2} s^2 }{ 2^{13} \sqrt{2} }
\nonumber \\
&= & \frac{ M \gamma^{3/2} }{ 2^{10} },
\end{eqnarray*}
and since $q>1/2$,
\begin{eqnarray*}
\sqrt{ 2^{(q+2)/q} C_{\tilde{\kappa}}^{1/q} C_{\tilde{\kappa}, q} (2 e^{\tau/2} )^{1/q}  }
&\leq & \frac{ n^{1/2} s^2 }{ 2^{13} \sqrt{2} } (\frac{ 2 q-1}{2 q}) \frac{1}{ (\sqrt{2} s)^{( 2q-1)/( 2q)} - (2^{-3} s^2)^{(2q-1)/(2q)} }
\nonumber \\
&=&
\frac{ M }{ 2^{11} \sqrt{2} }
(\frac{ 2 q-1}{2 q}) \frac{1}{ (\sqrt{2} s)^{( 2q-1)/( 2q)} - (2^{-3} s^2)^{(2q-1)/(2q)} }.
\end{eqnarray*}
Also with $\psi(M, t^2, n)$ as in Lemma \ref{la:4.10}, we have
\begin{eqnarray*}
H^B ( \frac{ t}{10}, {\cal F}_{\tilde{\kappa}, q, n}, \rho_{{\cal F}_{\tilde{\kappa}, q, n} } ) &\leq &
\frac{ 
2^{(q+2)/q} C_{\tilde{\kappa}}^{1/q} C_{\tilde{\kappa}, q} (2 e^{\tau/2} )^{1/q}  
}{ (\sqrt{2} s)^{1/q } }
\nonumber \\
&\leq &
\frac{ M^2 }{ 2^{23} s^{1/q} }
(\frac{ 2q-1}{2q})^2
\frac{1}{ [ (\sqrt{2} s)^{( 2q-1)/( 2q)} - (2^{-3} s^2)^{(2q-1)/(2q)} ]^2 }
\nonumber \\
&\leq & \frac{ M^2}{ 2^{23} s^2} 
\nonumber \\
&< & \frac{ M^2}{ 4 s^2 ( 64 c_0 + 1)}
\nonumber \\
&=& \frac{\gamma}{4} \psi( M, t^2, n).
\end{eqnarray*}
Consequently, we observe from Lemma \ref{la:4.10} that
\begin{displaymath}
P_{s,r}^* [ \sup_{
\| p_{s_1, r_1}^{1/2}- p_{s, r}^{1/2} \|_2^2 \leq 2 s^2, p_{s_1,r_1} \in {\cal F}_{\tilde{\kappa}, q, n} }
\nu_n(  \tilde{Z}_{p_{s_1, r_1}} ) \geq \frac{ n^{1/2} s^2 }{4}]
\leq 3 \exp[ - \frac{ n s^2 }{ 2^7 (2^6 c_0 +1) }].
\end{displaymath}
Let $A = \{ p_{s_1, r_1} \in {\cal F}_{\tilde{\kappa}, q, n}: s^2 \leq \| p_{s_1,r_1}^{1/2}- p_{s,r}^{1/2} \|_2^2 \leq 2 s^2\}$.
By Lemma 4 of Wong and Shen (1995), 
\begin{displaymath}
\sup_A E( \tilde{Z}_{p_{s_1, r_1}} ) \leq -(1-\delta) s^2,
\end{displaymath}
where $\delta = 2 e^{-\tau/2} ( 1- e^{-\tau/2} )^{-2} = 5/8$.
Now
\begin{eqnarray*}
&& P_{s,r}^* \{ \sup_A \prod_{i=1}^n \frac{ p_{s_1, r_1} (\{w_{i,1},\cdots, w_{i, N_i (T)} \}) }{
p_{s, r} (\{w_{i,1},\cdots, w_{i, N_i (T)} \}) }
\geq \exp[ - n s^2 (1 - \delta - \frac{1}{4})] \}
\nonumber \\
&\leq &
P_{s,r}^* \{ \sup_A \nu_n (\tilde{Z}_{p_{s_1, r_1} } ) \geq \frac{ n^{1/2} s^2}{ 4} \},
\end{eqnarray*}
and hence
\begin{eqnarray*}
&& P_{s,r}^* \{ \sup_{ s^2 \leq \|p_{s_1, r_1}^{1/2} - p_{s, r}^{1/2} \|_2^2 \leq 2 s^2, 
p_{s_1, r_1}\in {\cal F}_{\tilde{\kappa}, q, n} } \prod_{i=1}^n \frac{ p_{s_1, r_1} (\{w_{i,1},\cdots, w_{i, N_i (T)} \}) }{
p_{s, r} (\{w_{i,1},\cdots, w_{i, N_i (T)} \}) }
\geq e^{ - n s^2/8} \}
\nonumber \\
&\leq &
3 \exp[ - \frac{ n s^2 }{ 2^7 (2^6 c_0 +1) }].
\end{eqnarray*}
Let $L$ be the smallest integer such that $2^L \varepsilon^2 \geq 4 \geq \max 
\{ \| p_{s_1, r_1}^{1/2} - p_{s, r}^{1/2} \|_2^2:
p_{s_1, r_1}\in {\cal F}_{\tilde{\kappa}, q, n} \}$.
Then
\begin{eqnarray*}
&& P_{s,r}^* \{ \sup_{ \| p_{s_1, r_1}^{1/2} - p_{s, r}^{1/2}\|_2 \geq \varepsilon, 
p_{s_1, r_1}\in {\cal F}_{\tilde{\kappa}, q, n} } \prod_{i=1}^n \frac{ p_{s_1, r_1} (\{w_{i,1},\cdots, w_{i, N_i (T)} \}) }{
p_{s, r} (\{w_{i,1},\cdots, w_{i, N_i (T)} \}) }
\geq e^{ - n \varepsilon^2/8}  \}
\nonumber \\
&\leq & \sum_{j=0}^L
P_{s,r}^* \{ \sup_{2^j \varepsilon^2 \leq  \| p_{s_1, r_1}^{1/2}- p_{s, r}^{1/2} \|_2^2 < 2^{j+1} \varepsilon^2, 
p_{s_1, r_1}\in {\cal F}_{\tilde{\kappa}, q, n} } \prod_{i=1}^n \frac{ p_{s_1, r_1} (\{w_{i,1},\cdots, w_{i, N_i (T)} \}) }{
p_{s, r} (\{w_{i,1},\cdots, w_{i, N_i (T)} \}) }
\geq e^{ - n \varepsilon^2/8} \}
\nonumber \\
&\leq & 3 \sum_{j=0}^L
\exp[ - \frac{ 2^j n \varepsilon^2 }{ 2^7 (2^6 c_0 +1) }].
\end{eqnarray*}
Hence we conclude from (\ref{eq:a.56}) that
\begin{eqnarray*}
&& P_{s,r}^* \{ \sup_{ \| p_{s_1, r_1}^{1/2}- p_{s, r}^{1/2} \|_2 \geq \varepsilon, 
p_{s_1, r_1}\in {\cal F}_{\tilde{\kappa}, q, n} } \prod_{i=1}^n \frac{ p_{s_1, r_1} (\{w_{i,1},\cdots, w_{i, N_i (T)} \}) }{
p_{s, r} (\{w_{i,1},\cdots, w_{i, N_i (T)} \}) }
\geq e^{ - n \varepsilon^2/8} \}
\nonumber \\
&\leq &
4 \exp[ - \frac{ n \varepsilon^2 }{ 2^7 (2^6 c_0 +1) }].
\end{eqnarray*}
This proves Proposition \ref{pn:4.1}. \hfill $\Box$

\begin{la} \label{la:a.43}
Let $s^\dag_n (t)$ and $r^\dag_n (t)$ be as in (\ref{eq:3.23}). Then
\begin{displaymath}
E_{s,r} ( \frac{ p_{s,r} }{p_{s^\dag_n, r^\dag_n } } -1 )  \leq  C_{\tilde{\kappa}, 1} \delta_n,
\end{displaymath}
where $C_{\tilde{\kappa}, 1}$ is a constant depending only on $\tilde{\kappa}$.
\end{la}
{\sc Proof.}
Let $\bar{\kappa} = \kappa_0\vee 1$. Then
\begin{eqnarray*}
&& E_{s,r} (\frac{ p_{s,r} }{p_{s^\dag_n, r^\dag_n }} -1)
\nonumber \\
&=&
\sum_{j=0}^\infty \int_{0<w_1<\cdots < w_j<T}
\frac{ e^{-\int_0^T s(t) r(t- w_{\zeta(t)} ) dt } \prod_{i=1}^j s(w_i) r(w_i - w_{i-1}) }{
e^{-\int_0^T s^\dag_n(t) r^\dag_n (t- w_{\zeta(t)} ) dt } \prod_{i=1}^j s^\dag_n (w_i) r^\dag_n (w_i - w_{i-1}) }
\nonumber \\
&&\hspace{0.5cm}\times
[ e^{-\int_0^T s(t) r(t- w_{\zeta(t)} ) dt } \prod_{i=1}^j s(w_i) r(w_i - w_{i-1}) -
\nonumber \\
&&\hspace{1cm}
e^{-\int_0^T s^\dag_n(t) r^\dag_n (t- w_{\zeta(t)} ) dt } \prod_{i=1}^j s^\dag_n (w_i) r^\dag_n (w_i - w_{i-1}) ]
dw_1 \cdots dw_j
\nonumber \\
&\leq &
e^{ 2 \bar{\kappa}^2 \delta_n (2 \bar{\kappa} +\delta_n) T} \sum_{j=0}^\infty \int_{0<w_1<\cdots < w_j<T}
| e^{-\int_0^T s(t) r(t- w_{\zeta(t)} ) dt } \prod_{i=1}^j s(w_i) r(w_i - w_{i-1}) -
\nonumber \\
&&\hspace{1cm}
e^{-\int_0^T s^\dag_n(t) r^\dag_n (t- w_{\zeta(t)} ) dt } \prod_{i=1}^j s^\dag_n (w_i) r^\dag_n (w_i - w_{i-1}) |
dw_1 \cdots dw_j.
\end{eqnarray*}
Now we observe that
\begin{displaymath}
| s^\dag_n (w_i) r^\dag_n (w_i - w_{i-1}) - s(w_i) r(w_i - w_{i-1}) |
\leq 2 \bar{\kappa}^2 \delta_n ( 2 \bar{\kappa} + \delta_n),
\end{displaymath}
and
\begin{eqnarray*}
&& | e^{ -\int_0^T s^\dag_n(t) r^\dag_n (t- w_{\zeta(t)} ) dt }
- e^{ -\int_0^T s(t) r(t- w_{\zeta(t)} ) dt } |
\nonumber \\
&\leq &
|1 -  e^{ \int_0^T [ s^\dag_n(t) r^\dag_n (t- w_{\zeta(t)} ) -
s(t) r(t- w_{\zeta(t)} ) ] dt } |
\nonumber \\
&\leq &
2 \bar{\kappa}^2 \delta_n (2 \bar{\kappa} + \delta_n) T \sum_{i=0}^\infty \frac{ [ 2 \bar{\kappa}^2 \delta_n (2 \bar{\kappa} + \delta_n) T ]^i }{(i+1)!}.
\end{eqnarray*}
This implies that for $j=1,2, \cdots,$
\begin{eqnarray*}
&& \int_{0<w_1<\cdots < w_j<T}
| e^{-\int_0^T s(t) r(t- w_{\zeta(t)} ) dt } \prod_{i=1}^j s(w_i) r(w_i - w_{i-1}) -
\nonumber \\
&&\hspace{1cm}
e^{-\int_0^T s^\dag_n(t) r^\dag_n (t- w_{\zeta(t)} ) dt } \prod_{i=1}^j s^\dag_n (w_i) r^\dag_n (w_i - w_{i-1}) |
dw_1 \cdots dw_j
\nonumber \\
&\leq &
\delta_n \Big\{ 2 \bar{\kappa}^{4 j + 2} (2 \bar{\kappa}+ \delta_n) T \sum_{i=0}^\infty \frac{ [2 \bar{\kappa}^2 (2 \bar{\kappa} + \delta_n) T ]^i }{(i+1)!}
+ 2 j \bar{\kappa}^{4 j -2} (2 \bar{\kappa} + \delta_n) \Big\}
\frac{ T^j}{j!}.
\end{eqnarray*}
This proves Lemma \ref{la:a.43}. \hfill $\Box$

\begin{la} \label{la:a.1}
Let $N(t), t\in [0, T),$  be a counting process with conditional intensity $\lambda_1(.|.)$ as
in (\ref{eq:1.1}).
Suppose that (\ref{eq:4.86}) holds and
\begin{equation}
\xi (t) := \lim_{\delta \downarrow 0} \frac{1}{\delta} P_{s, r} [ N( t+\delta ) - N(t ) =1 ],
\hspace{0.5cm} \forall t \in [0, T).
\label{eq:a.63}
\end{equation}
Then for $s_1, r_1 \in \Theta_{\tilde{\kappa}, q, n}$,
\begin{eqnarray*}
&&  \sum_{j=0}^{n_\theta } \int_{0<w_1<\cdots < w_j<T} p_{s,r} (\{w_1,\cdots, w_j\})  \log [\frac{p_{s,r} (\{w_1,\cdots, w_j\}) }{
p_{s_1, r_1} (\{w_1,\cdots, w_j\}) } ]  dw_1 \cdots dw_j
\nonumber \\
&=&
\int_0^T  \{ \frac{s_1 (t) }{s(t)} -1 - \log[ \frac{s_1 (t) }{ s(t)} ] \}
s(t) e^{-\int_0^t s(u) du } dt
\nonumber \\
&& + \int_0^T \int_0^t \{ \frac{ s_1 (t) r_1 (u) }{ s(t) r(u) } -1 - \log[ \frac{s_1 (t) r_1 (u)}{ s(t) r(u)} ] \}
\xi( t-u) s(t) r(u) e^{- \int_{t-u }^t s(v) r(v-t+u ) dv } du dt.
\end{eqnarray*}
Also if
\begin{displaymath}
\sum_{j=0}^{n_\theta } \int_{0<w_1<\cdots < w_j<T} p_{s,r} (\{w_1,\cdots, w_j\})  \log [\frac{p_{s,r} (\{w_1,\cdots, w_j\}) }{
p_{s_1, r_1} (\{w_1,\cdots, w_j\}) } ]  dw_1 \cdots dw_j \leq 1,
\end{displaymath}
then
\begin{eqnarray*}
&&  \sum_{j=0}^{n_\theta } \int_{0<w_1<\cdots < w_j<T} p_{s,r} (\{w_1,\cdots, w_j\})  \log [\frac{p_{s,r} (\{w_1,\cdots, w_j\}) }{
p_{s_1, r_1} (\{w_1,\cdots, w_j\}) } ]  dw_1 \cdots dw_j
\nonumber \\
&\geq & \min\{ 
\frac{1}{20 \int_0^T s(t) e^{-\int_0^t s(u) du} dt} , \frac{1}{200} \}
[ \int_0^T  | s_1 (t) - s(t) | 
e^{-\int_0^t s(u) du } dt ]^2,
\end{eqnarray*}
and
\begin{eqnarray*}
&&  \sum_{j=0}^{n_\theta } \int_{0<w_1<\cdots < w_j<T} p_{s,r} (\{w_1,\cdots, w_j\})  \log [\frac{p_{s,r} (\{w_1,\cdots, w_j\}) }{
p_{s_1, r_1} (\{w_1,\cdots, w_j\}) } ]  dw_1 \cdots dw_j
\nonumber \\
&\geq & \min\{ 
\frac{1}{20 \int_0^T \int_0^t \xi( t-u) s(t) r(u) e^{-\int_{t-u}^t s(v) r( v-t+u) dv} du dt} , \frac{1}{200} \}
\nonumber \\
&& \hspace{0.5cm}\times
[ \int_0^T \int_0^t | s_1 (t) r_1 (u) - s(t) r(u) |
\xi( t-u) e^{- \int_{t-u }^t s(v) r(v-t+u ) dv } du dt ]^2.
\end{eqnarray*}
\end{la}
{\sc Proof.}
Writing the expectation as a Lebesgue-Stieltjes integral [cf.\ Aalen (1978) and Karr (1987);
see also Miller (1985), page 1455, for a different approach], we observe that
\begin{eqnarray*}
&&  \sum_{j=0}^{n_\theta } \int_{0<w_1<\cdots < w_j<T} p_{s,r} (\{w_1,\cdots, w_j\})  \log [\frac{p_{s,r} (\{w_1,\cdots, w_j\}) }{
p_{s_1, r_1} (\{w_1,\cdots, w_j\}) } ]  dw_1 \cdots dw_j
\nonumber \\
&=&
E_{s,r} \log [\frac{p_{s,r} (\{w_1,\cdots, w_{N(T)} \}) }{
p_{s_1, r_1} (\{w_1,\cdots, w_{N(T)} \}) } ]  
\nonumber \\
&=&
E_{s,r} \{ - \int_0^T s(t) r(t- w_{N(t)} ) dt + \sum_{j=1}^{N(T)} \log [s(w_j) r(w_j - w_{j-1} ) ]
\nonumber \\
&&
+ \int_0^T s_1(t) r_1(t- w_{N(t)} ) dt - \sum_{j=1}^{N(T)} \log [s_1 (w_j) r_1 (w_j - w_{j-1} ) ]
\}
\nonumber \\
&=& 
-\int_0^T  s (t ) 
e^{-\int_0^t s(u) du } dt
+ \int_0^T \log [s (t) ] e^{-\int_0^t s(u) du} s(t) dt
\nonumber \\
&& - \int_0^T \int_0^t s (t) r (u) 
\xi( t-u) e^{- \int_{t-u }^t s(v) r(v-t+u ) dv } du dt
\nonumber \\
&& + \int_0^T \int_0^t   \log[ s (t) r( u) ]
\xi( t-u) e^{- \int_{t-u }^t s(v) r(v-t+u ) dv }
s(t) r( u) du dt
\nonumber \\
&& +\int_0^T  s_1 (t ) 
e^{-\int_0^t s(u) du } dt
- \int_0^T \log [s_1 (t) ] e^{-\int_0^t s(u) du} s(t) dt
\nonumber \\
&& + \int_0^T \int_0^t s_1 (t) r_1 (u) 
\xi( t-u) e^{- \int_{t-u }^t s(v) r(v-t+u ) dv } du dt
\nonumber \\
&& - \int_0^T \int_0^t   \log[ s_1 (t) r_1 ( u) ]
\xi( t-u) e^{- \int_{t-u }^t s(v) r(v-t+u ) dv }
s(t) r( u) du dt
\nonumber \\
&=&
\int_0^T  \{ \frac{s_1 (t) }{s(t)} -1 - \log[ \frac{s_1 (t) }{ s(t)} ] \}
s(t) e^{-\int_0^t s(u) du } dt
\nonumber \\
&& + \int_0^T \int_0^t \{ \frac{ s_1 (t) r_1 (u) }{ s(t) r(u) } -1 - \log[ \frac{s_1 (t) r_1 (u)}{ s(t) r(u)} ] \}
\xi( t-u) s(t) r(u) e^{- \int_{t-u }^t s(v) r(v-t+u ) dv } du dt.
\end{eqnarray*}
Next suppose that
\begin{displaymath}
\sum_{j=0}^{n_\theta } \int_{0<w_1<\cdots < w_j<T} p_{s,r} (\{w_1,\cdots, w_j\})  \log [\frac{p_{s,r} (\{w_1,\cdots, w_j\}) }{
p_{s_1, r_1} (\{w_1,\cdots, w_j\}) } ]  dw_1 \cdots dw_j
\leq 1.
\end{displaymath}
Since $y-1-\log(y) \geq 0$ for all $y \in [0, \infty)$ with equality only if $y=1$, 
we conclude that
\begin{equation}
\int_0^T  \{ \frac{s_1 (t) }{s(t)} -1 - \log[ \frac{s_1 (t) }{ s(t)} ] \}
s(t) e^{-\int_0^t s(u) du } dt \leq 1,
\label{eq:a.31}
\end{equation}
and
\begin{displaymath}
\int_0^T \int_0^t \{ \frac{ s_1 (t) r_1 (u) }{ s(t) r(u) } -1 - \log[ \frac{s_1 (t) r_1 (u)}{ s(t) r(u)} ] \}
\xi( t-u) s(t) r(u) e^{- \int_{t-u }^t s(v) r(v-t+u ) dv } du dt
\leq 1.
\end{displaymath}
Observing that
\begin{displaymath}
y-1 -\log (y) \geq \frac{ (y-1)^2}{10} {\bf 1}\{y\in (0, 6)\}
+\frac{ y-1}{10} {\bf 1}\{y\geq 6\},
\end{displaymath}
it follows from (\ref{eq:a.31}) that
\begin{eqnarray*}
&& \int_0^T  \{ \frac{s_1 (t) }{s(t)} -1 - \log[ \frac{s_1 (t) }{ s(t)} ] \}
s(t) e^{-\int_0^t s(u) du } dt 
\nonumber \\
&\geq &
\frac{1}{10} \int_0^T  ( \frac{s_1 (t) }{s(t)} -1 )^2 {\bf 1}\{ \frac{ s_1(t) }{ s(t)} \in (0, 6)\}
s(t) e^{-\int_0^t s(u) du } dt 
\nonumber \\
&& + \frac{1}{10} \int_0^T  ( \frac{s_1 (t) }{s(t)} -1 ) {\bf 1}\{ \frac{ s_1(t) }{ s(t)} \geq 6\}
s(t) e^{-\int_0^t s(u) du } dt 
\nonumber \\
&\geq &
\frac{1}{10 \int_0^T s(t) e^{-\int_0^t s(u) du} dt} [ \int_0^T  | \frac{s_1 (t) }{s(t)} -1 | {\bf 1}\{ \frac{ s_1(t) }{ s(t)} \in (0, 6)\}
s(t) e^{-\int_0^t s(u) du } dt ]^2
\nonumber \\
&& + \frac{1}{100} [ \int_0^T  | \frac{s_1 (t) }{s(t)} -1 | {\bf 1}\{ \frac{ s_1(t) }{ s(t)} \geq 6\}
s(t) e^{-\int_0^t s(u) du } dt ]^2
\nonumber \\
&\geq & \min\{ 
\frac{1}{20 \int_0^T s(t) e^{-\int_0^t s(u) du} dt} , \frac{1}{200} \}
[ \int_0^T  | \frac{s_1 (t) }{s(t)} -1 | 
s(t) e^{-\int_0^t s(u) du } dt ]^2.
\end{eqnarray*}
In a similar manner, we have
\begin{eqnarray*}
&& \int_0^T \int_0^t \{ \frac{ s_1 (t) r_1 (u) }{ s(t) r(u) } -1 - \log[ \frac{s_1 (t) r_1 (u)}{ s(t) r(u)} ] \}
\xi( t-u) s(t) r(u) e^{- \int_{t-u }^t s(v) r(v-t+u ) dv } du dt
\nonumber \\
&\geq & \min\{ 
\frac{1}{20 \int_0^T \int_0^t \xi( t-u) s(t) r(u) e^{-\int_{t-u}^t s(v) r( v-t+u) dv} du dt} , \frac{1}{200} \}
\nonumber \\
&& \hspace{0.5cm}\times
[ \int_0^T \int_0^t | \frac{ s_1 (t) r_1 (u) }{ s(t) r(u) } -1 |
\xi( t-u) s(t) r(u) e^{- \int_{t-u }^t s(v) r(v-t+u ) dv } du dt ]^2.
\end{eqnarray*}
This proves Lemma \ref{la:a.1}.
\hfill $\Box$

{\sc Proof of Lemma \ref{la:3.1}.}
Following Yatracos (1988), page 1183, we observe that
\begin{eqnarray}
&& \sup \{ E_{s, r} [ \int_0^T | \tilde{s}_n (t) - s(t) | dt]: s \in \Theta_{\tilde{\kappa}, q},
r \in \Theta_{\theta, \tilde{\kappa}, q}  \} 
\nonumber \\
&\geq &
\sup \{ E_{s_1, r_1} [ \int_0^T | \tilde{s}_n (t) - s_1 (t) | dt]: s_1 \in \tilde{\Theta}_{\tilde{\kappa}, q, n} \}
\nonumber \\
&\geq &
\frac{1}{{\rm card}(\tilde{\Theta}_{\tilde{\kappa}, q, n} ) }
\sum_{ s_1 \in \tilde{\Theta}_{\tilde{\kappa}, q, n} }
E_{s_1, r_1} [ \int_0^T | \tilde{s}_n (t) - s_1 (t) | dt],
\label{eq:3.2}
\end{eqnarray}
for any $r_1\in \Theta_{\theta, \tilde{\kappa}, q}$.
Define $\tilde{s}_n^* \in \tilde{\Theta}_{\tilde{\kappa}, q, n}$ such that
\begin{displaymath}
\int_0^T | \tilde{s}_n (t) - \tilde{s}_n^* (t) | dt
= \inf \{ \int_0^T | \tilde{s}_n (t) - s_1(t) | dt: s_1 \in \tilde{\Theta}_{\tilde{\kappa}, q, n} \}.
\end{displaymath}
Then we have for $s_1 \in \tilde{\Theta}_{\tilde{\kappa}, q, n}$, 
\begin{eqnarray*}
\int_0^T |\tilde{s}_n^* (t) - s_1(t) | dt
&\leq & 
\int_0^T | \tilde{s}_n^* (t) - \tilde{s}_n (t) | dt
+ \int_0^T | \tilde{s}_n (t) - s_1(t)| dt
\nonumber \\
&\leq & 2 \int_0^T | \tilde{s}_n (t) - s_1(t) | dt.
\end{eqnarray*}
So
\begin{eqnarray}
&&
\frac{1}{{\rm card}(\tilde{\Theta}_{\tilde{\kappa}, q, n} ) }
\sum_{ s_1 \in \tilde{\Theta}_{\tilde{\kappa}, q, n} }
E_{s_1, r_1} [ \int_0^T | \tilde{s}_n (t) - s_1 (t) | dt]
\nonumber \\
&\geq &
\frac{1}{2 {\rm card}(\tilde{\Theta}_{\tilde{\kappa}, q, n} ) }
\sum_{ s_1 \in \tilde{\Theta}_{\tilde{\kappa}, q, n} }
E_{s_1, r_1} [ \int_0^T | \tilde{s}_n^* (t) - s_1 (t) | dt]
\nonumber \\
&\geq &
\inf \{ \int_0^T |s_1(t) - s_2(t) | dt: s_1\neq s_2,
s_1, s_2 \in \tilde{\Theta}_{\tilde{\kappa}, q, n} \}
\nonumber \\
&&\hspace{0.5cm}\times
\frac{1}{ 2 {\rm card}(\tilde{\Theta}_{\tilde{\kappa}, q, n} ) }
\sum_{ s_1 \in \tilde{\Theta}_{\tilde{\kappa}, q, n} }
P_{s_1, r_1} ( \tilde{s}_n^* \neq s_1).
\label{eq:3.3}
\end{eqnarray}
We observe from Fano's lemma [cf.\ Ibragimov and Has'minskii (1981), pages 323 to 325, or
Yatracos (1988), page 1182] that
\begin{eqnarray}
&& \frac{1}{ {\rm card}(\tilde{\Theta}_{\tilde{\kappa}, q, n} ) }
\sum_{ s_1 \in \tilde{\Theta}_{\tilde{\kappa}, q, n} }
P_{s_1, r_1} ( \tilde{s}_n^* \neq s_1 )
\nonumber \\
&\geq &
1 - \frac{1}{ \log [{\rm card}(\tilde{\Theta}_{\tilde{\kappa}, q, n}) -1 ] }
\Big\{ \log 2 +
\nonumber \\
&&\hspace{0.5cm}
+ \frac{1}{ [{\rm card} (\tilde{\Theta}_{\tilde{\kappa}, q, n})]^2 }
\sum_{s_1, s_2 \in \tilde{\Theta}_{\tilde{\kappa}, q, n}} E_{s_1, r_1} \log [ \prod_{i=1}^n \frac{ p_{s_1, r_1} (\{ w_{i,1},\cdots, w_{i, N_i(T)}\} )
}{
p_{s_2, r_1} (\{ w_{i,1},\cdots, w_{i, N_i(T)}\} ) } ] \Big\}.
\label{eq:3.4}
\end{eqnarray}
(\ref{eq:3.5}) now follows from (\ref{eq:3.2}), (\ref{eq:3.3}) and (\ref{eq:3.4}).
(\ref{eq:3.6}) is proved in a similar manner. \hfill $\Box$

\section{Appendix B}

{\sc Proof of {\rm (\ref{36})}.}
Let $u=\delta T^{-1/2}$ for some
$0 < \delta < \kappa$. The
contribution to $C_{\bw,u}$ from (\ref{32}) is equal to (under $P_{\theta_\bw}$)
\begin{eqnarray}
& & -T^{-1} \sum_{i=1}^d \sum_{y \in \by^{(i)}} \frac{d^2}{dv^2} g_\bw^{(i)}(v) \Big|_{v=y}
+ o(1) \cr
&= & -T^{-1} \sum_{i=1}^d \lambda_i \int_0^T \Big[ \frac{d^2}{dv^2} g_\bw^{(i)}(v)
\Big] e^{\theta_\bw g_\bw^{(i)}(v)} dv + o(1) \hspace{0.5cm} {\rm a.s. \ as } \ T \rightarrow \infty.
\label{39}
\end{eqnarray}
If $\frac{d^2}{dv^2} g_\bw^{(i)}(v)$ exists for all $v$ in an interval $(v_0,v_1)$,
then by integration by parts,
\begin{eqnarray}
& & -\int_{v_0}^{v_1} \Big[ \frac{d^2}{dv^2} g_\bw^{(i)}(v)
\Big] e^{\theta_\bw g_\bw^{(i)}(v)} dv \cr
&= & - \frac{d}{dv} g_\bw^{(i)}(v)
e^{\theta_\bw g_\bw^{(i)}(v)} \Big|_{v=v_0}^{v=v_1} + \theta_\bw
\int_{v_0}^{v_1} \Big[ \frac{d}{dv} g_\bw^{(i)}(v)
\Big]^2 e^{\theta_\bw g_\bw^{(i)}(v)} dv.
\label{40}
\end{eqnarray}
By letting $v_1 \uparrow h_{j+1}$ and $v_0 \downarrow h_j$, where $h_j$
and $h_{j+1}$ are adjacent points in $H_i$, it follows from (\ref{40}) that
(\ref{39}) is equal to 
\begin{equation}
T^{-1} \sum_{i=1}^d \lambda_i \sum_{h \in H_i} \Big( \frac{d}{dv} g_\bw^{(i)}
(v) \Big|_{v \downarrow h} - \frac{d}{dv} g_\bw^{(i)}
(v) \Big|_{v \uparrow h} \Big) e^{\theta_\bw g_\bw^{(i)}(h)} + \theta_\bw
\tau_\bw +o(1).
\label{41}
\end{equation}
By the law of large numbers, the contribution to $C_{\bw,u}$ from (\ref{33}) is equal to
\begin{eqnarray}
&& \frac{2}{u^2 T} \sum_{i=1}^d \lambda_i \sum_{h \in H_i} \Big[ \int_h^{h+u}
(h+u-y) \; dy \Big] \Big( \frac{d}{dv} g_\bw^{(i)}(v) \Big|_{v 
\uparrow h} - \frac{d}{dv} g_\bw^{(i)}
(v) \Big|_{v \downarrow h} \Big) e^{\theta_\bw g_\bw^{(i)}(h)}+o(1) 
\label{42} \\
&=& T^{-1}  \sum_{i=1}^d \lambda_i \sum_{h \in H_i} 
\Big( \frac{d}{dv} g_\bw^{(i)}(v) \Big|_{v 
\uparrow h} - \frac{d}{dv} g_\bw^{(i)}
(v) \Big|_{v \downarrow h} \Big) e^{\theta_\bw g_\bw^{(i)}(h)}+o(1) 
\nonumber
\end{eqnarray}
almost surely as $T \rightarrow \infty$. Since the contribution from
(\ref{34}) to $C_{\bw,u}$ is asymptotically negligible, it
follows from adding up (\ref{41}) and (\ref{42}) that 
\begin{equation}
\lim_{T \rightarrow \infty} [C_{\bw,u}-\theta_\bw \tau_\bw] \rightarrow 0 \hspace{0.5cm} \mbox{a.s.\
under $P_{\theta_\bw}$}.
\label{43}
\end{equation}
Since the first and second derivatives of $g_\bw^{(i)}$ are bounded and
continuous by (A2), it follows from (\ref{35})
that there exists $\beta_s \rightarrow 0$ as $s \rightarrow 0$ such that
\begin{equation}
\sup_{u \leq x \leq u+sT^{-1/2}} T | u^2C_{\bw,u}-x^2 C_{\bw,x}| \leq \beta_s
\label{44}
\end{equation}
for all large $T$ with probability 1. We can conclude (\ref{36}) from (\ref{43})
and (\ref{44}). \hfill $\Box$

{\sc Proof of Lemma {\rm \ref{l6}}.} By stationarity, we may assume without loss of generality 
$t=0$. Let $\ell \geq 1 $ and let us denote by $Q_\theta$ $(=Q_{\theta,\ell})$ 
the probability measure under which
$\by^{(i)}$ is generated as a Poisson point process on $[0,T +
\ell \kappa T^{-1/2})$ with intensity
$$
\eta_i(u) = \lambda_i \exp[\theta \wtd g_\bw^{(i)}(u)] \ {\rm for \ all} \ 
1 \leq i \leq d,
%\label{47}
$$
where $\wtd g_\bw^{(i)}(u) = g_\bw^{(i)}(u)+g_\bw^{(i)}(u-\ell \kappa T^{-1/2})$. 
Let
\begin{equation}
\wtd \phi_\bw(c) = \sup_{\theta > 0} \Big[ 2 \theta c - T^{-1} \sum_{i=1}^d
\lambda_i \int_0^{T+\ell \kappa T^{-1/2}} (e^{\theta \tilde g_\bw^{(i)}(u)}-1) \ du
\Big]
\label{48}
\end{equation}
and let $\wtd \theta_\bw > 0$ attain the supremum on the right hand side
of (\ref{48}). Define $\wtd S_x = S_x + S_{x+\ell \kappa T^{-1/2}}$. 
It follows from the arguments in (\ref{24}), (\ref{26}) and (\ref{28}) that
$$
E_{\tilde \theta_\bw} [\wtd S_0] = 2c, \quad
E_{\tilde \theta_\bw} \Big[ \frac{d}{dx} \wtd S_x \Big|_{x=0} \Big] = O(T^{-1}),
%\label{49}
$$
where $E_{\tilde \theta_\bw}$ denotes expectation with respect to
$Q_{\tilde \theta_\bw}$, and
$$
{\rm Cov}_{\tilde \theta_\bw} \pmatrix{ \wtd S_0 \cr
\frac{d}{dx} \wtd S_x \big|_{x=0} } \sim T^{-1}
\pmatrix{ \wtd v_\bw & 0 \cr 0 & \wtd \tau_\bw }, 
%\label{50}
$$
where $\wtd v_\bw$ and $\wtd \tau_\bw$ are defined as in (\ref{8}) but with
$\int_0^{T+\ell \kappa T^{-1/2}}$ replacing $\int_0^T$, $\wtd g_\bw^{(i)}$
replacing $g_\bw^{(i)}$ and $\wtd \theta_\bw$ replacing $\theta_\bw$. 
Hence for intervals $I_{1,T}$, $I_{2,T}$ satisfying the conditions of Lemma 3,
\begin{eqnarray*}
& & Q_{\tilde \theta_\bw} \Big\{ T^{1/2} \Big( \wtd S_0-2c,
\frac{d}{dx} \wtd S_x \Big|_{x=0} \Big) \in  I_{1,T} \times
I_{2,T} \Big\} 
\nonumber \\
&\sim & (2 \pi)^{-1} (\wtd v_\bw \wtd \tau_\bw
)^{-1/2} \Big( \int_{z_1 \in I_{1,T}} e^{-z_1^2/(2
\tilde v_\bw)} \ dz_1 \Big) \Big(
\int_{z_2 \in I_{2,T}} e^{-z_2^2/(2 \tilde \tau_\bw)} \ dz_2 \Big).
%\label{51}
\end{eqnarray*}
By the arguments in the proof of Lemma \ref{l5}, it follows that analogous to
(\ref{31}), 
\begin{eqnarray}
P_\bw(A_0 \cap A_{\ell \kappa T^{-1/2}}) & \leq & 
P_\bw \Big\{ \sup_{0 < u < 
\kappa T^{-1/2}} \wtd S_u \geq \max(2c,\wtd S_0,\wtd S_{\kappa T^{-1/2}}) \Big\} \cr
& \sim & \kappa T^{-1/2} \wtd \zeta_\bw e^{-T \tilde \phi_\bw
(c)},
\label{52}
\end{eqnarray}
for all $\ell$ and large $T$, where $\wtd \zeta_\bw = (2 \pi)^{-1} (\wtd \tau_\bw/
\wtd v_\bw)^{1/2}$. It remains for us to show that there exists
a constant $\gamma > 0$ such that with probability 1,
\begin{eqnarray}
\wtd \phi_\bw(c) & \geq & \theta_\bw c - T^{-1} \sum_{i=1}^d \lambda_i
\int_0^{\ell \kappa T^{-1/2}+T} (e^{\theta_\bw \tilde g_\bw^{(i)}(u)/2}-1) \ du \cr
& \geq & \phi_\bw(c) + \gamma \min \{ \ell^2 \kappa^2 T^{-1},1 \}
\label{53}
\end{eqnarray}
for all $1 \leq \ell \leq T^{3/2}/\kappa+1$ and $T$ large, so that (\ref{46})
follows by adding up (\ref{52}) over $1 \leq \ell \leq
T^{3/2}/\kappa+1$. 

The first inequality in (\ref{53}) follows directly from letting $\theta = 
\theta_\bw/2$ in the right hand side of (\ref{48}). By the identity
\begin{eqnarray*}
& & 2 e^{\theta_\bw[g_\bw^{(i)}(u) + g_\bw^{(i)}(u-\ell \kappa T^{-1/2} )]/2} \cr
&= & 
e^{\theta_\bw g_\bw^{(i)}(u)} +e^{\theta_\bw
g_\bw^{(i)}(u-\ell \kappa T^{-1/2})} - (e^{\theta_\bw g_\bw^{(i)}(u)/2} -e^{\theta_\bw
g_\bw^{(i)}(u-\ell \kappa T^{-1/2})/2})^2,
\end{eqnarray*}
it follows from (\ref{5}) that
\begin{eqnarray}
& & \theta_\bw c- \frac{1}{T} \sum_{i=1}^d \lambda_i \int_0^{\ell \kappa T^{-1/2}+T}
(e^{\theta_\bw \tilde g_\bw^{(i)}(u)/2}-1) \ du \cr
&= & \phi_\bw(c) + \frac{1}{ 2 T}
\sum_{i=1}^d \lambda_i \int_0^{\ell \kappa T^{-1/2}+T} 
(e^{\theta_\bw g_\bw^{(i)}(u)/2} -e^{\theta_\bw
g_\bw^{(i)}(u-\ell \kappa T^{-1/2})/2})^2 \ du
\label{55}
\end{eqnarray}
and the second inequality of (\ref{53}) indeed holds for all large $T$ with probability 1.
\hfill $\Box$

{\sc Proof of Lemma {\rm \ref{l10}}.} The proof of
Lemma \ref{l10} uses arguments similar to the proof of Lemma \ref{l6}. 
The main changes are in replacing $\kappa T^{-1/2}$ by $\kappa T^{-1}$.
Analogous to (\ref{52}), there exists a constant $C > 0$ such that
with probability 1,
\begin{eqnarray}
P_\bw(A_t \cap A_{t+\ell \kappa T^{-1}}) & \leq & P_\bw \Big\{ S_t+
S_{t+\ell \kappa T^{-1}} < 2c, \sup_{t < u \leq t+\kappa T^{-1}}
(S_u+S_{u+\ell \kappa T^{-1}}) \geq 2c \Big\} \cr
& \leq & C \kappa T^{-1/2}
e^{-T \tilde \phi_\bw(c)}
\label{93}
\end{eqnarray}
for all large $T$, if $\kappa$ is chosen large enough. 
By (\ref{55}) (with $\kappa T^{-1/2}$ replaced by $\kappa T^{-1}$), the
following analogue to (\ref{53}),
\begin{equation}
\wtd \phi_\bw(c) \geq \phi_\bw(c) + \gamma \min \{ \ell \kappa T^{-1},1 \}
\label{94}
\end{equation}
holds with probability 1 for some $\gamma > 0$ and hence Lemma \ref{l10} follows from
(\ref{93}) and (\ref{94}). \hfill $\Box$

\end{document}